\documentclass[mnsc,nonblindrev]{informs4}

\DoubleSpacedXI

\usepackage{amsmath,mathtools}
\usepackage{appendix}
\usepackage{graphicx}
\usepackage{float} 
\usepackage{tikz, tikz-3dplot}
\usepackage{caption}
\usepackage{subcaption}
\usepackage{mathtools} 
\usepackage{comment}
\usepackage{algorithm2e}
\RestyleAlgo{ruled}
\definecolor{links}{RGB}{204,36,29}
\usepackage[colorlinks=true,breaklinks=true,bookmarks=true,urlcolor=links,citecolor=links,linkcolor=links,bookmarksopen=false,draft=false]{hyperref}
\usepackage{makecell}
\usepackage{threeparttable}
\usepackage{authblk}
\usepackage[figuresright]{rotating}
\RequirePackage{pdflscape}
\usepackage{multirow, booktabs}
\usepackage{pgfplots}
\pgfplotsset{compat=newest}

\usepackage{natbib}
 \bibpunct[, ]{(}{)}{,}{a}{}{,}%
 \def\bibfont{\small}%

\DeclareMathOperator{\mibp}{MIBP}

\DeclareMathOperator{\proj}{proj}

\DeclareMathOperator{\conv}{conv}

\DeclareMathOperator{\vertex}{vert}

\DeclareMathOperator{\pr}{Pr}

\DeclareMathOperator{\R}{\mathbb{R}}
\DeclareMathOperator{\Z}{\mathbb{Z}}

\newcommand{\yun}{\textcolor{black}} 

\newcommand{\our}{\textsf{CH}}
\newcommand{\ourip}{\textsf{CH-2}}
\newcommand{\conic}{\textsf{Conic}}
\newcommand{\milp}{\textsf{MILP}}

\newcommand{\ourC}{\textsf{CH-Chain}}
\newcommand{\ouripC}{\textsf{CH-Chain-2}}
\newcommand{\conicC}{\textsf{Conic-Chain}}
\newcommand{\milpC}{\textsf{MILP-Chain}}

\newcommand{\revision}{\textcolor{black}}
\newcommand{\rvtwo}{\textcolor{black}}

\TheoremsNumberedThrough     
\ECRepeatTheorems

\EquationsNumberedThrough    

\renewcommand{\theARTICLETOP}{}

\begin{document}


\RUNAUTHOR{Chen, He, Rong and Wang}

\RUNTITLE{IP Approach for Quick-Commerce Assortment Planning}

\TITLE{An Integer Programming Approach for Quick-Commerce Assortment Planning}


\ARTICLEAUTHORS{%
\AUTHOR{
Yajing Chen \qquad 
Taotao He \qquad 
Ying Rong \qquad 
Yunlong Wang}
\AFF{Antai College of Economics and Management, Shanghai Jiao Tong University, Shanghai 200030, China. \\
\EMAIL{yajingchen@sjtu.edu.cn} (YC), 
\EMAIL{hetaotao@sjtu.edu.cn}, \URL{https://taotaoohe.github.io/} (TH),    
\EMAIL{yrong@sjtu.edu.cn} (YR), 
\EMAIL{wylwork\_sjtu@sjtu.edu.cn} (YW)
\URL{}}
} 

\ABSTRACT{In this paper, we explore the challenge of assortment planning in the context of quick-commerce, a rapidly-growing business model that aims to deliver time-sensitive products. In order to achieve quick delivery to satisfy the immediate demands of online customers in close proximity, personalized online assortments need to be included in brick-and-mortar store offerings. With the presence of this physical linkage requirement and distinct multinomial logit choice models for online consumer segments, the firm seeks to maximize overall revenue by selecting an optimal assortment of products for local stores and by tailoring a personalized assortment for each online consumer segment. We employ an integer programming approach to solve this NP-hard problem to global optimality.  \revision{In particular, we derive convex hull results to represent the consumer choice of each online segment under a general class of operational constraints, and to characterize the relation between assortment decisions and choice probabilities of products.}
Our convex hull results, coupled with a modified choice probability ordered separation algorithm, yield formulations that provide a significant computational advantage over existing methods.  \revision{Finally, we illustrate how our convex hull results can be used to address other assortment optimization problems}. 
}




\KEYWORDS{quick commerce; assortment optimization; multinomial logit model; mixed-integer nonlinear programming; convexification; cutting-plane}
\HISTORY{This version: \today}

\maketitle
%


\section{Introduction}
Over recent years, brick-and-mortar stores have witnessed slower growth compared to the online sector. Nevertheless, there exist certain product categories, including meals, fresh food, groceries, and flowers, which are not well-suited for next-day delivery, which is the fastest service offered by many online retail giants. Consumers demand these time-sensitive products to be delivered within an hour to satiate their immediate demands. With the advent of third-party delivery platforms such as DoorDash, Delivery Hero, and Grab, brick-and-mortar stores are capitalizing on their proximity to local consumers to offer their products online and deliver them within a short timeframe. This strategy has given birth to the rapidly-evolving business model of quick-commerce, \revision{which is expected to produce a global market volume of 265 billion US dollars in 2029 \citep{statista2024}}.

As brick-and-mortar stores are increasingly teaming up with third-party delivery platforms to offer quick delivery services, some online retail giants are taking an alternative approach by establishing their own local stores that serve both online and offline consumers. Alibaba's Hema Fresh is a prime example of such an approach. The store guarantees delivery within 30 minutes to online consumers residing within three kilometers of their local stores. By leveraging in-store staff to handpick products for online consumers, Hema Fresh is able to enhance its sales per square foot and even to venture into opening massive stores in prime downtown areas. Similarly, Amazon chose to acquire Whole Foods and to roll out speedy delivery to the consumers in the selected cities.

In traditional omni-channels, the constraints of physical space often lead to a smaller assortment of products available in brick-and-mortar stores compared to that of the online channel. \revision{These offline outlets typically complement the online channel by allowing consumers to touch, feel, and sample products in person \citep{DzyaburaJagabathula2018,LoTopaloglu2022}, thereby influencing their attractiveness in the online channel. However, the quick-commerce model generally focuses on providing everyday items that consumers are already familiar with.} The advantage of quick-commerce lies in leveraging the in-store product assortment to facilitate prompt delivery and to satisfy the immediate demands of local consumers.

In this paper, we explore the challenge of product assortment in the context of quick commerce. To meet the need for speedy delivery, it is essential for all products, whether offered online or offline, to be available at brick-and-mortar stores. Offline consumers can choose their desired products in the conventional manner, whereas the delivery platform or online retailers can leverage extensive consumer data, such as browsing and purchase behavior, to offer personalized assortments to online consumers. 

Specifically, we assume that the choice behavior of both offline consumers and each of the online consumer segments follow separate multinomial logit (MNL) models. The linkage between the offline and online channels in the quick commerce setting is established by ensuring that the personalized assortments for each online segment are carried at a local brick-and-mortar store. The goal is to maximize the expected revenue by finding the optimal assortment to carry at the local store and the personalized assortment for every online consumer segment, subject to the linkage constraint and the other operational constraints. We refer to this problem as the quick-commerce assortment problem.

In this paper, we formally define the quick-commerce assortment problem as in~(\ref{eq:omni-mnl}). We show that this problem is NP-hard even when there are no operational constraints and there is one single online consumer type. Nevertheless, we develop integer programming techniques for solving~\eqref{eq:omni-mnl} to its global optimality. In particular, we propose convexification techniques to address the  combinatorial and nonconvex nature of the optimization problem. Our methodology yields formulations of~\eqref{eq:omni-mnl} that are not only capable of providing provably optimal solutions at a large scale but also offer flexibility in accommodating complex operational considerations.

\subsection{Contribution}
We summarize the contribution of our study as follows:

\begin{enumerate}
    \item 
    We prove in Theorem~\ref{theorem:partial-convexification} that modeling each online consumer type’s choice probabilities via their convex hull, \textit{i.e.}, tightest convex relaxation, yields an exact formulation. In Theorem~\ref{them:extended-formulation}, we further show that, by leveraging the geometry underlying Luce’s choice axiom, this convex hull admits a compact polyhedral representation, provided the convex hull of operational constraints does.   \revision{As an illustrative example, we specialize our results to the case where online consumers make choices using the two-stage Luce choice model introduced in~\cite{EcheniqueSaito2019}, which satisfactorily handles zero-probability choices. In Section~\ref{section:MNL}, we leverage this convex hull result to derive linear programming (LP) formulation for constrained assortment planning under MNL. When specializing to totally unimodular constraints, our result recovers  LP formulations in~\cite{avadhanula2016tightness} and~\cite{sumida2021revenue}}.
    \item We conduct a convex hull study on the relation between offered products and the resulting choice probabilities. While a convex hull description is generally intractable for the full product set, we provide a tractable characterization for the \textit{single-product} case (Theorems~\ref{them:staircase}–\ref{them:convexhull}). In Section~\ref{section:poly}, we leverage the geometry of this convex hull to design a polynomial-time randomized rounding algorithm for a variant of~\eqref{eq:omni-mnl}, where each online consumer type follows an independent demand model and the firm faces a specific class of precedence constraints (see Theorem~\ref{them:polyLP}). This generalizes the LP formulation for assortment optimization under mixed independent and MNL demand models presented in~\cite{cao2023revenue}.
    

    \item We apply the convex hull results to obtain formulations for the quick-commerce assortment problems. Our formulations, coupled with a modified choice probability ordered cutting-plane algorithm (refer to Algorithm~\ref{alg:cutting}), outperform  the formulations based on the state-of-the-art conic integer optimization approach~\citep{sen2018conic}. \revision{In Section~\ref{section:frontwarehouse}, we illustrate that our formulations can be modified to handle quick-commerce assortment planning problems with additional front warehouses, which subsumes the problem setting studied in~\cite{HousniTopaloglu2022}.}
    
    
\end{enumerate}

\subsection{Literature Review}
Our research is related to the field of optimizing product assortments by using random utility choice models, specifically the widely used MNL model \citep{Luce1959, McFadden1977}. In this case, without restrictions on offered assortments, \cite{GallegoIyengarPhillipsDubey2004_techreport} and \cite{TalluriVanRyzin2004} show that the optimal assortment under the MNL model is revenue-ordered, and~\cite{avadhanula2016tightness,sumida2021revenue} obtain an LP formulation for assortment optimization under MNL model with additional totally unimodular constraints. \revision{Various optimization techniques have been proposed to solve assortment planning under other choice models, including the mixture of MNL models~\citep{DesirGoyalZhang2022}, the mixture of MNL and linear demand model~\citep{cao2023revenue}  the nested logit model~\citep{davis2014assortment}, the Markov chain based choice model~\citep{blanchet2016markov}, and random consideration set model~\citep{gallego2024random}.}

The majority of research on assortment optimization problems focuses on a single channel context. However, with the rise of omni-channel retailing, assortment planning in omni-channel systems has emerged as a relatively new and important topic for both academics and practitioners. \cite{DzyaburaJagabathula2018} study the problem of determining the subset of products from the retailer's online channel to offer in the offline channel to maximize the aggregate revenue. It is assumed that each product is defined by a set of attributes, and there is a utility associated with each attribute that depends on whether the product is offered in the offline channel. They incorporate the impact of physical evaluation on preferences into the consumer demand model. Under this model, they demonstrate that the decision problem is NP-hard and propose approximation algorithms with theoretical guarantees. \cite{LoTopaloglu2022} introduce a novel features-tree structure to organize products by features in an omni-channel setting for the assortment optimization problem of a retailer that operates a physical store and an online store. 
They show that the assortment optimization problem in this setup is NP-hard, and leverage the features tree structure to obtain approximately optimal assortments. \cite{chen2022offline} study the optimal offline store locations and location-dependent assortments decision problem for an omni-channel retailer in the presence of an online channel that carries all products. They develop a tractable mixed integer second-order conic programming reformulation and explore the structural properties of the reformulation to derive strengthening cuts in closed form. 

Our paper most closely related to the growing literature on developing integer programming techniques to solve optimization problems involving customer choice model.~\cite{sen2018conic} propose a conic integer optimization formulation for assortment problem under MMNL model.~\cite{BertsimasMisic2019} propose a new mixed integer optimization model of the product line design problem under a first-choice rule, and solve the model using benders decomposition.~\cite{ChenMisic2021} develop a mixed-integer optimization methodology for solving the assortment optimization problem when the choice model is a decision forest model.~\cite{chen2022offline} develop a tractable mixed integer second-order conic programming reformulation for location-dependent offline-channel assortment planning in omni-channel retailing.\revision{~\cite{he2024discrete} provide strong mixed-integer linear programming formulations for solving revenue maximization problems with discrete decision variables.} \rvtwo{~\cite{linexpress} propose mixed-integer linear programming formulations for solving assortment planning with multi-choice rank list model in e-commerce}.   

Last, we comment on the relation of our paper to the optimization literature. First, our paper provides new insights into the Charnes-Cooper transformation which is proposed in~\cite{charnes1962programming} to solve linear fractional programs. Second, our convex hull characterization of choice probabilities under MNL model with operational constraints is related to~\cite{Megiddo1979}, which shows that optimizing a rational function over a combinatorial set is polynomial time solvable if optimizing a linear function over the feasible region is polynomial time solvable. Last, our polyhedral study on the relation between  assortment decision and choice probabilities under MNL is related to the literature on mixed-integer bilinear programming. We review recent developments in this topic as follows. Convex hulls of bilinear functions are studied in~\citep{tawarmalani2013explicit,gupte2020extended}.~\cite{gupte2013solving} present a mixed-integer linear programming formulation for mixed integer bilinear problems. Convex hulls of various mixed integer bilinear sets are studied in~\citep{tawarmalani2010strong,fampa2021convexification,tawarmalani2024new}. \revision{In the context of QAP, the bilinear structure is special in the sense that the function domain is given by the MNL models, and such structure has not been exploited in the literature.}

\subsection{Organization}
The remainder of this paper is organized as follows. 
In Section \ref{section:model}, we formulate our quick-commerce assortment planning problem \eqref{eq:omni-mnl} and show its NP-hardness. 
In Section \ref{section:reformulation}, we give the relaxed model \eqref{eq:omni-rlx} and show its exactness. 
\revision{In Section \ref{section:convexhull}, we present convex hull results.}
\revision{In Section~\ref{section:CH}, we present explicit formulations and cutting-plane implementations.}
In Section \ref{section:computation}, we present computational experiments. 
\revision{In Section \ref{section:applications}, we discuss applications of our methods.}
Finally, conclusions are given in Section \ref{section:conclusion}.


\def \p {\boldsymbol{p}}
\def \q {\boldsymbol{q}}
\def \Y {\boldsymbol{Y}}
\def \X {\mathcal{F}}
\def \calZ {\mathcal{Z}}
\def \x {\boldsymbol{x}}
\def \B {\mathcal{B}}
\def \r {\boldsymbol{r}}
\def \y {\boldsymbol{y}}
\def \u {\boldsymbol{u}}
\def \z {\boldsymbol{z}}

\def \v {\boldsymbol{v}}
\def \V {\mathcal{V}}
\def \c {\phi}
    
\section{Model}\label{section:model}

\subsection{Formulation}
\def \N {\mathcal{N}}
\def \M {\mathcal{M}}
We consider a set of products $N:=\{1, \ldots,n\}$. Due to enormous data available from online consumers, we assume that the firm can segment the online consumers into $m$ types.  Let $M:=\{1, \ldots, m\}$ denote the set of consumer types in the online channel and let $\{0\}$ denote the only consumer segment in the offline channel. Let $M^+ = \{0\} \cup M$ denote the set of all consumer types, which include $\{0\}$ and $M$. For each consumer type $i \in M^+$ and each product $j \in N$, let $r_{ij}$ denote the revenue obtained from selling the product $j$ to the type $i$ consumer, which allows for modeling personalized pricing such as discounts for VIP consumers or personalized discount offered by mobile application \citep{ElmachtoubGuptaHamilton2021}. The probability of a consumer of type $i$ arriving to the system is $\alpha_i$, where $\sum_{i \in M^+} \alpha_i = 1$. 

\revision{As noted by \citet{chen2022offline} in the context of omni-channel retailing, a consumer's decision to visit an offline store is influenced by the store's overall attractiveness and the distance required to travel. This overall attractiveness is not impacted by differences in the assortment between online and offline channels. This assumption also applies to the quick-commerce business model. For example, Alibaba's Hema Fresh enhances its offline appeal by offering cooking and dining services for the fresh food sold in its stores. Moreover, these stores are typically located in large shopping malls, allowing consumers to address other needs beyond making purchases. Additionally, the size of consumer type $i$ for $i > 0$ is determined by the firm's consumer segmentation strategy, which considers factors such as price sensitivity and preferences for various product types, regardless of their availability. Therefore, in this paper, we assume that $\alpha_i$ is exogenously given and independent of the assortment offered to each consumer type.}

For each type $i \in M^+$, let $\x_i = (x_{i1}, \ldots, x_{in} ) \in \{0,1\}^n$ be a binary decision variable to model the subset of products offered to consumer type $i$, namely, $x_{ij}=1$  if and only if  product $j$ is offered to consumer type $i$. In the online channel, $\x_i$ can be different as the firm can utilize the personalized assortment to enhance revenue \citep{HousniTopaloglu2022}. In addition, due to the prompt delivery requirement, the personalized online assortment needs to be included in the physical store, \textit{i.e.}, $\x_0 \geq \x_i$.

The choice behavior of each consumer type, from either offline or online channels, follows the multinomial logit (MNL) model. For each consumer type $i \in M^+$, we use $u_{ij} \geq 0$ to denote their preference weight on  product $j,$ and $u_{i0}$ to denote their preference weight on the no-purchase option. If we offer $\x_i \in \{0,1\}^n$ to a consumer of type $i$, the consumer purchases product $j$ with probability $u_{ij}x_{ij}/(u_{i0} + \sum_{j \in N}u_{ij}x_{ij})$. Thus, the expected revenue obtained from consumer type $i$ is given by 
\[
R_i^{\text{MNL}}(\x_i): = \frac{\sum_{j \in N}r_{ij}u_{ij}x_{ij}}{u_{i0} + \sum_{j \in N}u_{ij}x_{ij}}.
\]

\def \portionoff {\lambda^{\text{off}}} 
\def \portionon {\lambda^{\text{on}}} 


For a consumer of type $i \in M$, we assume that products recommended to her must be stocked in the physical store. Our goal is to jointly find feasible offline assortment and personalized assortments to online consumers to maximize the expected revenue over all consumer types. More specifically, we are interested in solving the following nonlinear discrete optimization problem:
\begin{equation}\label{eq:omni-mnl}
	\begin{aligned}
 \max\quad & \sum_{i \in M^+} \alpha_{i}R_i^{\text{MNL}}(\x_i)\\
\text{s.t.}\quad &\x_0 \geq \x_i  && \text{for } i \in M \\
& \x_i \in \X_i \subseteq \{0,1\}^n   && \text{for } i \in M^+, \\
\end{aligned}\tag{\textsc{QAP}}
\end{equation}
where the first constraint  models that products recommended to online consumers must be available in the offline physical  store, and $\X_i$ is the set of feasible assortments that we can offer to consumer type $i \in M^+$. \revision{Our methodology merely requires the offline feasible region $\X_0$  to be an arbitrary binary polyhedral set, that is, $\X_0$ is the intersection of $\{0,1\}^n$ and a system of linear inequalities. This allows us to model complex business rules for offline operations. For each online feasible region $\X_i$, we assume that its tightest convex relaxation (convex hull) can be constructed. This class of feasible regions includes various constraints studied in the assortment optimization literature, such as matroid constraints~\citep{ghuge2022constrained,udwani2023submodular}, totally unimodular constraints~\citep{avadhanula2016tightness,sumida2021revenue}, and permutahedron~\citep{tserenjigmid2021order}. In particular, we will focus on exemplifying  the usefulness of our methods with the two-stage Luce model in~\cite{echenique2018perception}, as we detailed as follows.}

\subsubsection{\revision{Two-stage Luce model}}\label{section:2SLM} In reality, an online consumer segment may never purchase certain products within the same category. According to the dataset of Taobao, \cite{chen2022data} find that the consumer-item matrices are very sparse (only $1.25\%$ nonzero entries), indicating the zero-probability scenario of consumer behavior in real-world scenarios. Many utility based discrete choice models are not able to capture zero-probability scenarios. For example, if a person always prefers Coke to Pepsi, then the probability of choosing Pepsi is always zero if Coke is in the assortment. However, if Coke is unavailable, this person has a positive probability  purchasing Pepsi. Such behavior cannot be modeled by the MNL model. To address this,~\cite{EcheniqueSaito2019} propose a two-stage Luce model (2SLM) to accommodate the dominance relationship among some products. Under the 2SLM, given an assortment, all the \textit{dominated products} are removed by the consumers in the first stage, and then, at the second stage, they make the choice according to the MNL model based on the remaining undominated products.

To use the 2SLM to model choice behavior in our context, we need to model all undominated products, \rvtwo{\textit{i.e.}, products not removed in the first stage.} 
For $i \in M$, let $\mathcal{P}_i$ be a partial order set specifying consumers' preference on the set of $n$ products, where the partial order is denoted by $\preceq_i$. A strict order $a_1 \prec_i a_2$ means that $a_1 \preceq_i a_2$ and $a_1 \neq a_2$.  
A \textit{chain} in $\mathcal{P}_i$ is a collection of elements that are totally ordered. \rvtwo{For $i \in M$, given an assortment $S \subseteq N$, the undominated products in $S$, proposed in~\cite[Definition 2]{EcheniqueSaito2019}, are defined as 
\[
c_i(S) := \bigl\{j \in S \bigm| \not\exists k \in S \text{ s.t. } j \prec_i k \bigr\}.
\]
Assume that there is no constraint on online assortments. Then, the set of all undominated products is $\bigl\{c_i(S) \bigm | S \subseteq N \bigr\}$, which has the following binary representation:}
\begin{equation}\label{eq:undominated}
\Bigl\{\x_i\in \{0,1\}^n \Bigm| 0\leq x_{ia_1} + \cdots + x_{ia_k} \leq 1 \text{ for chains } a_1 \prec_i \cdots \prec_i a_k \text{ in } \mathcal{P}_i \Bigr\}. \tag{\textsc{Undominated}}
\end{equation}
\revision{Henceforth, if  a consumer of type $i$ chooses according 2SLM then our objective is to maximize the revenue function $R_i^{\text{MNL}}(\x_i)$ over the undominated set~\eqref{eq:undominated}.}

\subsection{Two step revenue-ordered policy} \label{subsec:2S-RO}
\rvtwo{Before exactly solving~\eqref{eq:omni-mnl}, we first propose a heuristic approach that yields near-optimal solutions.} When there is no constraint on assortments, that is, $\X_i = \{0,1\}^n$ for $i \in M^+$, a natural idea for solving (\ref{eq:omni-mnl}) is to use a two step revenue-ordered (RO) policy. In the first step, we determine the offline assortment by maximizing revenue from the offline channel alone, \rvtwo{temporarily ignoring the contribution from the online channel}. Since this corresponds to solving a standard unconstrained MNL assortment optimization problem, the resulting  offline assortment follows a revenue-ordered structure. In the second step, the personalized online assortment for each consumer type is chosen from the offline assortment obtained in the first step. Similarly, the personalized online assortment obtained in the second step also possesses the revenue-order property. Formally, the RO policy can be obtained by solving the following sequential problems:
\begin{equation}\label{eq:2SNR}
\begin{aligned}
\x^{\textsc{RO}}_0 &\in \argmax \bigl\{R_0^{\text{MNL}}(\x_0) \bigm| \x_0 \in \{0,1\}^n \bigr\} \\
\x^{\textsc{RO}}_i  &\in \argmax \bigl\{R_i^{\text{MNL}}(\x_i) \bigm| \x_i \in \{0,1\}^n,\ \x^{\textsc{RO}}_0 \geq \x_i \bigr\} \quad \text{ for } i \in M,
\end{aligned}\tag{\textsc{2-Step-RO}}
\end{equation}
Unfortunately, this strategy fails to solve~(\ref{eq:omni-mnl}) to global optimality, as we illustrate in the next example.   
\begin{example}\label{ex:RO-is-bad}
\revision{Here, we consider a toy example with three products and two consumer types in the online channel. The arrival rates of the three types are $\alpha = (0.4, 0.4, 0.2)$, and the revenues of products are $(r_{01}, r_{02}, r_{03}) =(r_{11}, r_{12}, r_{13}) = (r_{21}, r_{22}, r_{23}) =(10, 9, 8)$. The preference weight  of products are $(u_{01}, u_{02}, u_{03}) = (100, 100, 1)$, $(u_{11}, u_{12}, u_{13}) = (1, 1, 100)$ and $(u_{21}, u_{22}, u_{23}) =(100, 1, 1)$, and the preference weights of the no-purchase option is $u_{00}=u_{10}=u_{20}=1$. As shown in Table~\ref{table:2stageRO},~(\ref{eq:2SNR}) offers product $\{1\}$ to all three consumer types, while~(\ref{eq:omni-mnl}) offers $\{1,3\}$ to consumer type 0 and type 1 and $\{1\}$ to consumer type 2. By adding product $\{3\}$ to offline assortment $\{1\}$, the increases of online revenue offsets the decrease of the offline revenue decreases. \hfill  \Halmos
\begin{table}[H]
\fontsize{9pt}{10pt}\selectfont 
\renewcommand\tabcolsep{2pt}
\renewcommand{\arraystretch}{1.1}
  \centering
  \caption{Assortment and revenue of two methods}
  \resizebox{\textwidth}{!}{
    \begin{tabular}{cccccccrccc}
    \toprule
    \multirow{2}[4]{*}{\textbf{Channel}} & \multirow{2}[4]{*}{\textbf{Type}} & \multirow{2}[4]{*}{\textbf{Rate}} &   & \multicolumn{3}{c}{\textbf{\ref{eq:2SNR}}} &  & \multicolumn{3}{c}{\textbf{\ref{eq:omni-mnl}}} \\
    \cmidrule{5-7}\cmidrule{9-11}  &  &  &  & \textbf{Assort.} & \textbf{Exp.Rev} & \textbf{Total Rev.} &  & \textbf{Assort.} & \textbf{Exp.Rev} & \textbf{Total Rev.} \\
    \midrule
    Offline & 0 & 0.4 &  & \{1\} & 9.90  & \multicolumn{1}{c}{\multirow{3}[1]{*}{\makecell{0.4*9.90+0.4*5.00 \\ +0.2*9.90=7.94}}} &  & \{1,3\} & 9.88  & \multicolumn{1}{c}{\multirow{3}[2]{*}{\makecell{0.4*9.88+0.4*7.94\\+0.2*9.90=9.11}}} \\
    Online & 1  & 0.4   &   & \{1\} & 5.00  &  &  & \{1,3\} & 7.94  &  \\
     & 2   & 0.2   &   & \{1\} & 9.90  & &  & \{1\}   & 9.90  &  \\
    \bottomrule
    \end{tabular}}
  \label{table:2stageRO}
\end{table}}
\end{example}
The key insight behind this example is that the RO policy ignores the needs of online consumers completely when determining the assortment for the offline channel. Online consumers who do not appreciate the high priced-products often leave the firm's online channel without making a purchase. To provide a better experience to those consumers, it is necessary  to expand the assortments. \rvtwo{To better assess the performance of the revenue-ordered policy, Section~\ref{section:comparisonwithRO} introduces an improved revenue-ordered heuristic and presents numerical evidence showing that its optimality gap can be significant across a range of scenarios.}



\subsection{The Connection between~(\ref{eq:omni-mnl}) and Other Assortment Problems}\label{sec:connection_literature}

The assortment under quick commerce is closely related to two other assortment problems. The first is the assortment under the mixture of MNL models \citep{bront2009column,rusmevichientong2014assortment,sen2018conic,DesirGoyalZhang2022}. The difference between (\ref{eq:omni-mnl}) and this stream of literature is that the assortment for each consumer type can be different. If we modify the constraint in \eqref{eq:omni-mnl} to require $\x_0 = \x_i$ instead of $\x_0 \geq \x_i$, then it becomes assortment optimization under the mixture of MNL models. Thus, \eqref{eq:omni-mnl} generates more revenue by allowing personalized assortments for different online consumer types.  

The second closely related problem is to analyze the value of personalized assortment with each consumer type following the MNL model \citep{HousniTopaloglu2022}. In the case where the offline channel is disabled (i.e., $\alpha_0 = 0$), there is a cardinality constraint on the assortment from the offline channel (i.e., $\X_0 = \{\x_0\mid \sum_{j\in N} x_{0j} \leq k\}$), and the prices of all products are the same across different consumer types (i.e., $r_{ij} = r_{i'j}$ for $i \neq i'$),  then (\ref{eq:omni-mnl}) reduces to the problem studied in \citet{HousniTopaloglu2022}.

Similar to the aforementioned assortment problems, the following proposition also shows that (\ref{eq:omni-mnl}) is an NP-hard problem. 
\begin{proposition}\label{prop:complexity}
The maximization problem $(\ref{eq:omni-mnl})$ is NP-hard even when there is only one online consumer type, $\X_0=\X_1 = \{0,1\}^n$, and $\r_0 = \r_1$.
\end{proposition}
There are two different directions to tackle difficult assortment problems. The first is to develop approximation algorithms \citep{rusmevichientong2014assortment, DesirGoyalZhang2022, HousniTopaloglu2022}. The other direction is to aim for global optimality in solving such problems \citep{sen2018conic,BertsimasMisic2019,ChenMisic2021}. In this paper, we take the latter approach by deriving integer programming formulation techniques for~(\ref{eq:omni-mnl}). 

It is also noteworthy that solving~(\ref{eq:omni-mnl}) to global optimality not only generates more revenue but also yields more selections to offline stores than the revenue-ordered heuristic  proposed in Section~\ref{subsec:2S-RO}. 
\begin{proposition}\label{prop:more-offering}
Assume that the price ranks of products are consistent across different consumer types, that is, $r_{i1} \geq r_{i2} \geq \cdots \geq  r_{in}$ for each $i \in M^+$.  Let $\x^{\textsc{RO}}:=(\x^{\textsc{RO}}_0, \ldots, \x^{\textsc{RO}}_m)$ be a solution given by~$(\ref{eq:2SNR})$. Let $\x^{\textsc{QAP}}:=(\x^{\textsc{QAP}}_0, \ldots, \x^{\textsc{QAP}}_m)$ be the optimal solution of~$(\ref{eq:omni-mnl})$. When there exist multiple optimal solutions of~$(\ref{eq:omni-mnl})$, $\x^{\textsc{QAP}}$ is the optimal solution with the maximal number of activated products in all channels. Then we have $\x^{\textsc{RO}} \leq \x^{\textsc{QAP}}$.
\end{proposition}
To gauge consumer satisfaction in the assortment planning setting, the total utility of the offered assortment is utilized by \citet{sumida2021revenue}. 
Higher total utility of the offered assortment corresponds to greater consumer satisfaction. 
Proposition~\ref{prop:more-offering} suggests that  optimizing assortments jointly for offline and online channels, as opposed to sequentially, can lead to simultaneous increases in overall revenue and consumer satisfaction for all consumer types simultaneously. This is due to the higher total utility of the optimal assortment achieved under joint optimization.


\subsection{\revision{Preliminaries}}
 \revision{Before presenting our formulation of (\ref{eq:omni-mnl}), we review some basic concepts that we use throughout this paper}. First, we review a prevalent transformation which is proposed in~\cite{charnes1962programming} to solve linear fractional programming. For each  $i \in M^+$, the Charnes-Cooper (C-C) transformation, denoted as $\Pi_i$, maps an assortment  offering $\x_i \in \{0,1\}^n$ to a vector $(y_{i0},\y_i)$ defined as follows:
\begin{equation*}
y_{i0} = \frac{1}{u_{i0}+ \sum_{j \in N}u_{ij}x_{ij}} \qquad \text{and} \qquad y_{ij} =  \frac{x_{ij}}{u_{i0}+ \sum_{j \in N}u_{ij}x_{ij}}\quad \text{ for } j \in N.
\end{equation*}
The C-C transformation  has a few properties that are useful in our development. Algebraically,  for a given set of assortments $\X_i \subseteq \{0,1\}^n$, it is a one-to-one mapping between $\X_i$ and the image set
\begin{equation*}\label{eq:vertex-CC}
\Pi_i(\X_i):= \Biggl\{(y_{i0},\y_i) \Biggm| u_{i0} y_{i0} + \sum_{j \in  N} u_{ij} y_{ij} = 1,\ \y_i \in y_{i0} \cdot \X_i,\ y_{i0} > 0 \Biggr\},
\end{equation*}
where we recall that  for a set $S$ in $\R^n$ and a real scalar $\lambda>0$, the scalar multiple $\lambda \cdot S$ is defined as $\{\lambda x\mid x \in S\}$. Geometrically, the image set $\Pi_i(\X_i)$ can be viewed as the intersection of a hyperplane and a cone, given as follows:
\[
\biggl\{(y_{i0},\y_i) \biggm| u_{i0} y_{i0} + \sum_{j \in  N} u_{ij} y_{ij} = 1 \biggr\} \quad\text{ and } \quad \biggl\{(y_{i0},\y_i) \biggm| \y_i \in y_{i0} \cdot \X_i ,\ y_{i0} > 0 \biggr\},
\]
respectively. This geometrical interpretation will be used to prove the main result of Section~\ref{section:formulation-MNL}.  In addition, the vector $(y_{i0},\y_i)$ can be used to generate the choice probability in (\ref{eq:omni-mnl}). Specifically, for each  fixed point $(y_{i0}, \y_i)$ in $\Pi_i(\X_i)$, we may interpret  $u_{ij}y_{ij}$ as the probability that a consumer of type $i$ chooses item $j$ and $u_{i0}y_{i0}$ as the probability that such consumer chooses the outside option. In other words, the set $\Pi_i(\X_i)$ captures the choice probabilities among the products in the given assortment and the no-purchase alternative. In Section~\ref{section:conv-bilinear}, we use this interpretation to devise a (choice) probability-ordered algorithm for generating cutting planes. Henceforth, we will refer to $\Pi_i(\X_i)$ as \textit{choice probability set} over a given set of assortments, and we will refer to $\x_i$ as \textit{assortment variable} and $\y_i \in \Pi_i(\X_i)$ as \textit{choice probability variable}.

Next, we review some concepts from integer programming. Consider a set defined by linear inequalities and continuous and integer variables as follows
\[
E = \bigl\{(x,y,\lambda) \in \R^n \times \R^{s} \times \Z^t \bigm| \BFA x + \BFB y + \BFC \lambda \leq \BFb \bigr\},
\]
where $\BFA \in \R^{m\times n}$, $\BFB \in \R^{m\times s}$, $\BFC \in \R^{m\times t}$, and $\BFb \in \R^m$. We say that $E$ is a \textit{mixed integer programming} (MIP) \textit{formulation} of a set $S \subseteq \R^n$ if the projection of $E$ onto the space of $x$ variables is $S$, that is $S = \{x \in \R^n \mid \exists (y,\lambda) \text{ s.t. } (x,y,\lambda) \in E\}$. The polyhedron obtained by dropping all integrality requirements in $E$ is called the \textit{continuous relaxation} or \textit{linear programming} (LP) \textit{relaxation} of $E$. One of the factors that has a strong impact on the performance of an MIP formulation is the strength of the continuous relaxation~\citep{vielma2015mixed}. 
For two MIP formulations $E_1$ and $E_2$ of a set $S \subseteq \R^n$, we say $E_1$ is a \textit{tighter} than $E_2$ if the LP relaxation of $E_1$ is contained in that of $E_2$.  

When MIP formulations already have a number of constraints that is exponential in the data size of the problem, solving the corresponding LP relaxations is not straightforward. In this case,  we would like to solve the \textit{separation problem} of these linear programs, namely, given a polyhedron $P \subseteq \R^n$ and a point $\bar{\x} \in \R^n$, either show that $\bar{\x} \in P$ or give a valid inequality $ \langle \alpha, \x \rangle \leq \beta $ for $P$ such that $ \langle \alpha, \bar{\x} \rangle > \beta$. The celebrated theorem of~\cite{grotschel1981ellipsoid} states that optimizing a linear function over $P$ can be solved in polynomial time if and only if the separation problem over $P$ can be solved in polynomial time. 

One of the most successful techniques for deriving MIP formulations is \textit{convexification}, see~\citep{schrijver2003combinatorial,bertsimas2005optimization,conforti2014integer}.
The basic idea in convexification is to identify substructures that capture the essence in the problem at hand and then derive the convex hull of substructures.  Given a set $S \subseteq \R^n$, the \textit{convex hull} of $S$, denoted as $\conv(S)$, is the inclusionwise minimal convex set containing $S$. To prove a convex hull result, it is often useful to invoke its dual definition, that is, the set of all possible convex combinations of points in $S$.

\section{Reformulation framework}\label{section:reformulation}

In this section, we introduce a framework to reformulate~(\ref{eq:omni-mnl}) based on a convex relaxation approach. Namely, we relax the optimization problem of each online consumer type to its tightest convex relaxation while keeping the offline optimization problem the same.  We show that such relaxation is exact although~(\ref{eq:omni-mnl}) is non-separable due to the presence of linkage constraints $\x_0 \geq \x_i$ for each $i \in M$. 
This idea reduces the number of binary variables in from $mn+n$ to $n$. 
Moreover, it  yields a framework for constructing MIP formulations of~(\ref{eq:omni-mnl}), \revision{paving the way for our theoretical developments  in Section~\ref{section:convexhull}}.






By using the algebraic property of the C-C transformation, we obtain an mixed integer bilinear reformulation of~(\ref{eq:omni-mnl}), that is,
\begin{equation*}~\label{eq:CC-OAP}
\begin{aligned}
\max \quad &  \sum_{i \in M^+} \alpha_i\biggl( \sum_{j \in N} r_{ij} u_{ij}y_{ij}  \biggr)\\
	\text{s.t.} \quad 	 &  (y_{i0},\y_i) \in \Pi_i(\X_i) \quad\text{for } i \in \revision{M} \\
 &\y_i =  y_{i0}\x_i \quad \text{and} \quad \x_i \in \X_i  \quad \text{for } i \in M^+  \\	
	 	 &\x_0 \geq \x_i  \quad \text{for } i \in M. \\
	\end{aligned}
\end{equation*}
Since the objective function only involves the choice probability variable $(\y_i)_{i\in M^+}$, we use the relation $\x_i = \y_i / y_{i0}$ to project the online assortment variable $\x_i$ out of the feasible region for $i \in M$, and obtain an equivalent formulation: 
\begin{align}
\max \quad &   \sum_{i \in M^+} \alpha_i\biggl( \sum_{j \in N} r_{ij} u_{ij}y_{ij}  \biggr) \notag\\
	\text{s.t.} \quad &(y_{i0},\y_i) \in \Pi_i(\X_i)  \quad \text{for } i \in \revision{M} \tag{\textsc{Choice}}\label{eq:on-mnl} \\
 & \y_{0} =  y_{00}\x_0  \quad \text{and} \quad (y_{00},\y_0)  \in \Pi_0(\X_0)   \tag{\textsc{Off-Bi}}\label{eq:off-bi} \\	 
	 	 &  \y_i \leq y_{i0}\x_0 \quad \text{for } i \in M.\tag{\textsc{Link-Bi}} \label{eq:link-bi}
\end{align}
The constraint~(\ref{eq:on-mnl}) describes  \revision{online  consumers' choice behavior} using the choice probability set.
The bilinear constraint~(\ref{eq:off-bi}) models offline consumers' choice behavior using a bilinear relation between the assortment variable and the choice probability variable. Last, the bilinear inequality~(\ref{eq:link-bi}) models that if a product does not appear in the offline store then a consumer in an online channel purchases such a product with a probability of zero.

%

Now, we are ready to construct our relaxation for (\ref{eq:omni-mnl}). We relax each choice probability set in~(\ref{eq:on-mnl}) to its convex hull and obtain:
\begin{equation}
\begin{aligned}
\max \quad &  \sum_{i \in M^+} \alpha_i\biggl( \sum_{j \in N} r_{ij} u_{ij}y_{ij}  \biggr) \notag\\
	\text{s.t.} \quad 	 &(y_{i0},\y_i) \in \conv\bigl(\Pi_i(\X_i)\bigr)  \quad \text{for } i \in \revision{M} \\
	&(\ref{eq:off-bi})  \text{ and } 	 	 (\ref{eq:link-bi}).
\end{aligned}\tag{\textsc{QAP-Rlx}}\label{eq:omni-rlx}
\end{equation}
Theorem~\ref{theorem:partial-convexification} shows that such relaxation is an exact reformulation of~(\ref{eq:omni-mnl}). \rvtwo{The main step in the proof is to invoke a decomposition result, described as in Lemma~\ref{lemma:facial decomposition} and visualized in Figure~\ref{fig:facial-decom}. This lemma has been used in studying disjunctive programming~\citep[e.g. Lemma 5.1 in][]{balas1998disjunctive} and the lift-and-project rank for 0-1 linear programs~\citep[e.g. Theorem 5.22 in][]{conforti2014integer}. }
\begin{lemma}\label{lemma:facial decomposition}
	Consider a set $S \subseteq \R^n$ and a hyperplane $H:=\{\x \mid \alpha^\intercal \x =\beta \}$ such that $\alpha^\intercal \x \leq \beta$ for $\x \in S$. Then, $\conv(S) \cap H = \conv(S \cap H)$.
\end{lemma} 
\begin{figure}
\centering
\begin{subfigure}{.45\textwidth}
  \centering
    \captionsetup{font=normalsize}
  \begin{tikzpicture}[scale=0.6]
\draw[black] (0,-1) -- (3,-1);
\draw[black] (2,2) -- (3,0);
\draw[black] (0,-1) -- (0,3);
\draw[black] (2,2) -- (0,3);
\draw[black] (3,0) -- (3,-1);
\filldraw[black] (0,-1) circle (1.5pt); 
\filldraw[black] (1,-1) circle (1.5pt); 
\filldraw[black] (2,-1) circle (1.5pt); 
\filldraw[black] (3,-1) circle (1.5pt); 
\filldraw[black] (0,0) circle (1.5pt); 
\filldraw[black] (0,1) circle (1.5pt); 
\filldraw[black] (0,2) circle (1.5pt); 
\filldraw[black] (0,3) circle (1.5pt);
\filldraw[black] (1,0) circle (1.5pt); 
\filldraw[black] (1,1) circle (1.5pt); 
\filldraw[black] (1,2) circle (1.5pt); 
\filldraw[black] (2,0) circle (1.5pt); 
\filldraw[black] (2,1) circle (1.5pt); 
\draw[black] (1,4) -- (2,2);
\draw[blue,ultra thick] (2,2) -- (3,0);
\draw[black] (3,0) -- (4,-2);
\filldraw[red] (2,2) circle (2pt); 
\filldraw[red] (3,0) circle (2pt); 
\end{tikzpicture}
\end{subfigure}
\begin{subfigure}{.45\textwidth}
  \centering
    \captionsetup{font=normalsize}
  \begin{tikzpicture}[scale=0.6]
\draw[black] (0,-1) -- (3,-1);
\draw[black] (2,2) -- (3,0);
\draw[black] (0,-1) -- (0,3);
\draw[black] (2,2) -- (0,3);
\draw[black] (3,0) -- (3,-1);
\filldraw[black] (0,-1) circle (1.5pt); 
\filldraw[black] (1,-1) circle (1.5pt); 
\filldraw[black] (2,-1) circle (1.5pt); 
\filldraw[black] (3,-1) circle (1.5pt); 
\filldraw[black] (0,0) circle (1.5pt); 
\filldraw[black] (0,1) circle (1.5pt); 
\filldraw[black] (0,2) circle (1.5pt); 
\filldraw[black] (0,3) circle (1.5pt);
\filldraw[black] (1,0) circle (1.5pt); 
\filldraw[black] (1,1) circle (1.5pt); 
\filldraw[black] (1,2) circle (1.5pt); 
\filldraw[black] (2,0) circle (1.5pt); 
\filldraw[black] (2,1) circle (1.5pt); 
\filldraw[black] (2,2) circle (2pt); 
\filldraw[black] (3,0) circle (2pt); 
\draw[black] (0,4) -- (0.65,2.7);
\draw[blue,ultra thick] (0.65,2.7) -- (2.5,-1);
\draw[black] (2.5,-1) -- (3,-2);
\filldraw[red] (1,2) circle (2pt); 
\filldraw[red] (2,0) circle (2pt); 
\end{tikzpicture}
\end{subfigure}
\caption{\textnormal{\rvtwo{Illustration of the decomposition in Lemma~\ref{lemma:facial decomposition}. The line segment connecting two red dots represents $\conv(S \cap H)$, and the blue line segment represents $\conv(S) \cap H$. When the plane $H$ defines a face of set $S$ (left), the line segment connecting two red points is exactly equal to the blue line segment. However, when $H$ does \textit{not} define a face of set $S$ (right),  the line segment connecting two red points is strictly contained in the blue one.}}}\label{fig:facial-decom}
\end{figure}

Using Lemma~\ref{lemma:facial decomposition}, we derive the main result of Section~\ref{section:reformulation}.
\begin{theorem}~\label{theorem:partial-convexification}
Problem $(\ref{eq:omni-mnl})$ has the same optimal objective value as~$(\ref{eq:omni-rlx})$.
\end{theorem}
\rvtwo{As a consequence, to obtain an MIP formulation for~(\ref{eq:omni-mnl}), it suffices to derive the convex hull of $\Pi_i(\X_i)$ and MIP formulations for~(\ref{eq:off-bi}) and~(\ref{eq:link-bi}). 
In the next section, we will use a convexification approach to tackle these two challenges. }

\section{\revision{Convexification results}}\label{section:convexhull}

\subsection{Convex hull characterization of~\eqref{eq:on-mnl}}\label{section:formulation-MNL}
To utilize Theorem~\ref{theorem:partial-convexification}, an important steppingstone is to specify the condition that the convex hull of the choice probability set (i.e. $\Pi_i(\X_i)$) for each consumer type $i$ can be found efficiently.
The early works by \cite{avadhanula2016tightness} and \cite{sumida2021revenue} imply that  the convex hull of the choice probability set coincides with its linear programming (LP) relaxation when the feasible assortment set $\X_i$ is characterized by a totally unimodular matrix. \revision{This important finding provides LP formulations for the assortment planning under the MNL model with various constraints, including cardinality constraints, display location effects, discrete price menus, price ladder constraints and product precedence constraints. In this section, we use a geometric approach to extend this result. More specifically, we will show:}
\[
\textit{$\conv(\Pi_i(\X_i))$ admits a computationally tractable polyhedral description if $\conv(\X_i)$ does.}
\]
\revision{
Built upon this theoretical result, in Remark~\ref{rmk:2slm-formulation}, we are able to provide a tractable formulation to deal with the case where online consumers make choice decision according to the two-stage Luce model introduced in Section~\ref{section:2SLM}. }

Before we proceed with the main results in this section, we remove the subscript $i$ to streamline the presentation.  More specifically, we study the convex hull of the choice probability set $\Pi_i(\X)$ over a given set of feasible assortments $\X$, where
\[
\Pi(\X):= \Biggl\{(y_0,\y) \in \R^{n+1} \Biggm|u_0y_0  + \sum_{j  \in N}u_jy_j = 1,\ \y \in y_0 \cdot \X,\ y_0 > 0 \Biggr\},
\]
where $\X$ is a subset of $\{0,1\}^n$ and models constraints on available products.

The main idea is to use the geometry of the choice probability set to obtain its convex hull description. For a given nonempty subset $\X$ of $\{0,1\}^n$, the choice probability set $\Pi(\X)$ is the intersection of a hyperplane $H$ and a cone $K$, where,
\[
H:=\biggl\{(y_0,\y) \biggm| u_0y_0 + \sum_{j  \in N}u_jy_j = 1 \biggr\}\qquad \text{and} \qquad K:=\biggl\{(y_0,\y)\biggm| \y \in y_0 \cdot \X,\ y_0 > 0 \biggr\}. 
\]
Recall that the expression $u_i y_i$ can be interpreted as the probability of purchasing product $i$ and that $u_0 y_0$ can be interpreted as the nonpurchase probability. The hyperplane $H$ is used to enforce probability normalization, which requires that the sum of those probabilities equals one, \revision{one of the two axioms behind the MNL/Luce choice model~\citep{Luce1959}}. Hence, $H$ is referred to as the \textit{probability normalization} hyperplane. \revision{Another axiom behind the MNL/Luce model is independence of irrelevant alternatives (IIA)~\citep{Luce1959}. It essentially requires that the relative ratio of purchase probabilities of two \textit{offered} alternatives is independent of others}. \rvtwo{Consider a point $(y_0,\y)$ in the cone $K$ such that $y_a$ and $y_b$ are non-zero, that is, products $a$ and $b$ are offered to the consumer. It follows that  }
\[
\frac{u_a y_a}{u_b y_b} = \frac{u_ay_0}{u_by_0} = \frac{u_a}{u_b}.
\]
That is, every solution to the cone $K$ satisfies IIA. Thus, $K$ is referred to as the \textit{IIA-assortment} cone. The following example illustrates the geometry of $\Pi(\X)$ and its role in our study. 
\begin{example}\label{ex:geometry}
Let us consider the case with two products where $(u_0,u_1,u_2) = (1,1,2)$. We consider two structures of feasible assortments. One is  $\X = \{0,1\}^2$ in Figure~\ref{fig:decom-a}  and the other is $\X = \{(0,0), (1,0),(0,1)\}$ in Figure~\ref{fig:decom-b}. In both settings, the  probability normalization hyperplane $H$ is the same as the fixed value of preference weights, which is visualized as the light blue plane. The IIA-assortment cone $K$ is depicted as ultra-thick dark rays. The intersection between $H$ and $K$ are blue points, which correspond to the $\y$ variable generated by all feasible assortments specified by $\X$. Finally, one can visualize the red polygon, which depicts $\conv(H \cap K)$, as the intersection of the normalization hyperplane and the gray shadow, which depicts $\conv(K)$. \hfill \Halmos
\end{example}

\begin{figure}
\centering
\begin{subfigure}{.5\textwidth}
  \centering
    \captionsetup{font=normalsize}
\begin{tikzpicture}[scale=0.8, color={lightgray}]
    \definecolor{point_color}{rgb}{ 0,0,0 }
    \tikzstyle{point_style} = [fill=point_color]

    \coordinate (q00) at (0, 1*2.2, 0);
    \coordinate (q10) at (1/2*1.5, 1/2*2.2, 0);
    \coordinate (q01) at (0, 1/3*2.2, 1/3*1.8);
    \coordinate (q11) at (1/4*1.5, 1/4*2.2, 1/4*1.8);
    \coordinate (p1) at (-1*1.5, 1*2.2, 0.5*1.8);
    \coordinate (p2) at (1*1.5, 1*2.2, -0.6*1.8);
    \coordinate (p3) at (1.45*1.5, 1/4*2.2, -0.35*1.8);
    \coordinate (p4) at (-0.65*1.5, 1/4*2.2, 0.65*1.8);
    \coordinate (x00) at (0, 1.2*2.2, 0);
    \coordinate (x10) at (1.2*1.5, 1.2*2.2, 0);
    \coordinate (x01) at (0, 1.2*2.2, 1.2*1.8);
    \coordinate (x11) at (1.2*1.5, 1.2*2.2, 1.2*1.8);

    

\draw[black,->] (0,0,0) -- (1.2*1.5,0,0) node [anchor=west]{$y_{1}$};
\draw[black,->] (0,0,0) -- (0,3,0) node [anchor=south west, yshift=-0.5ex]{$y_{0}$};
\draw[black,->] (0,0,0) -- (0,0,1.4*1.8) node [left]{$y_{2}$};


    \definecolor{edge_color}{RGB}{100,100,100}

   \definecolor{plane_color}{RGB}{249, 6, 33}
    \definecolor{shadow_color}{RGB}{100,100,100}

    \tikzstyle{plane_style} = [fill=plane_color, fill opacity=0.8, draw=plane_color, line width=1 pt, line cap=round, line join=round] 
    \tikzstyle{shadow_style} = [fill=shadow_color, fill opacity=0.3, line width=0.01 pt, line cap=round, line join=round]

     \draw[fill=blue,fill opacity=0.05,draw= white,line width=0.001pt] (p1) -- (p2) -- (p3) -- (p4) -- (p1) -- cycle;

    \draw[shadow_style] (0.0,0.0,0.0) -- (0,1.2*2.2,1.2*1.8)-- (0,1.2*2.2,0) -- (0.0,0.0,0.0) -- cycle;
    \draw[shadow_style] (0.0,0.0,0.0) -- (0,1.2*2.2,0) -- (1.2*1.5, 1.2*2.2, 0)-- (0.0,0.0,0.0)-- cycle;
    \draw[shadow_style] (0.0,0.0,0.0) -- (1.2*1.5, 1.2*2.2, 0) --(1.2*1.5,1.2*2.2,1.2*1.8) -- (0.0,0.0,0.0)-- cycle;
    \draw[shadow_style] (0.0,0.0,0.0) --(1.2*1.5,1.2*2.2,1.2*1.8) -- (0,1.2*2.2,1.2*1.8)-- (0.0,0.0,0.0)-- cycle;

    \draw[plane_style] (q00) -- (q01) -- (q11) -- (q10) -- (q00) -- cycle;
    


\draw[black,ultra thick] (0,0,0) -- (1.2*1.5, 1.2*2.2, 0);
\draw[black,ultra thick] (0,0,0) -- (1.2*1.5,1.2*2.2,1.2*1.8);

\draw[black,ultra thick] (0,0,0) -- (0,1.2*2.2,1.2*1.8);
    
\draw[black, ultra thick] (0,0,0) -- (0,1.2*2.2,0);

\filldraw[blue] (q10) circle (1.5pt) ; 
\filldraw[blue] (q11) circle (1.5pt); 
\filldraw[blue] (q01) circle (1.5pt) ; 
\filldraw[blue] (q00) circle (1.5pt); 


\filldraw[gray] (x00) circle (1.5pt) ; 
\filldraw[gray] (x10) circle (1.5pt); 
\filldraw[gray] (x01) circle (1.5pt) ; 
\filldraw[gray] (x11) circle (1.5pt); 

%
 
  \end{tikzpicture}  \caption{$\X = \{0,1\}^2$}
  \label{fig:decom-a}
\end{subfigure}%
\begin{subfigure}{.5\textwidth}
  \centering
    \captionsetup{font=normalsize}
\begin{tikzpicture}[scale=0.8, color={lightgray}]
    \definecolor{point_color}{rgb}{ 0,0,0 }
    \tikzstyle{point_style} = [fill=point_color]

    \coordinate (q00) at (0, 1*2.2, 0);
    \coordinate (q10) at (1/2*1.5, 1/2*2.2, 0);
    \coordinate (q01) at (0, 1/3*2.2, 1/3*1.8);
    \coordinate (p1) at (-1*1.5, 1*2.2, 0.5*1.8);
    \coordinate (p2) at (1*1.5, 1*2.2, -0.6*1.8);
    \coordinate (p3) at (1.45*1.5, 1/4*2.2, -0.35*1.8);
    \coordinate (p4) at (-0.65*1.5, 1/4*2.2, 0.65*1.8);

\draw[black,->] (0,0,0) -- (1.2*1.5,0,0) node [anchor=west]{$y_{1}$};
\draw[black,->] (0,0,0) -- (0,3,0) node [anchor=south west, yshift=-0.5ex]{$y_{0}$};
\draw[black,->] (0,0,0) -- (0,0,1.4*1.8) node [left]{$y_{2}$};

    \definecolor{edge_color}{RGB}{100,100,100}

   \definecolor{plane_color}{RGB}{249, 6, 33}
    \definecolor{shadow_color}{RGB}{100,100,100}

    \tikzstyle{plane_style} = [fill=plane_color, fill opacity=0.8, draw=plane_color, line width=1 pt, line cap=round, line join=round] 
    \tikzstyle{shadow_style} = [fill=shadow_color, fill opacity=0.3, line width=0.01 pt, line cap=round, line join=round]

     \draw[fill=blue,fill opacity=0.05,draw= white,line width=0.001pt] (p1) -- (p2) -- (p3) -- (p4) -- (p1) -- cycle;
    \draw[shadow_style] (0.0,0.0,0.0) -- (0,1.2*2.2,1.2*1.8)-- (0,1.2*2.2,0) -- (0.0,0.0,0.0) -- cycle;
    \draw[shadow_style] (0.0,0.0,0.0) -- (0,1.2*2.2,0) -- (1.2*1.5, 1.2*2.2, 0)-- (0.0,0.0,0.0)-- cycle;
    \draw[shadow_style] (0.0,0.0,0.0) -- (1.2*1.5, 1.2*2.2, 0) --(0,1.2*2.2,1.2*1.8) -- (0.0,0.0,0.0)-- cycle;
    \draw[black, ultra thick] (0,0,0) -- (0,1.2*2.2,0);

    \draw[plane_style] (q00) -- (q01) -- (q10) -- (q00) -- cycle;


\draw[black,ultra thick] (0,0,0) -- (1.2*1.5, 1.2*2.2, 0);

\draw[black,ultra thick] (0,0,0) -- (0,1.2*2.2,1.2*1.8);
\filldraw[blue] (q10) circle (1.5pt) ; 
\filldraw[blue] (q01) circle (1.5pt) ; 
\filldraw[blue] (q00) circle (1.5pt); 

  \end{tikzpicture}  \caption{$\X=\{0,1\}^2 \setminus (1,1)$}
  \label{fig:decom-b}
\end{subfigure}
\caption{Convex Hull Characterization via Hyperplane-Cone Decomposition.}
\label{fig:decom}
\end{figure}

In Example~\ref{ex:geometry}, we obtain the convex hull of a choice probability set by taking the intersection of the normalization hyperplane and the convex hull of IIA-assortment cone. 
The following lemma shows that such an observation is not limited to the special structure of $\X$ in Example~\ref{ex:geometry}.
\begin{lemma}\label{lemma:fractional-LP}
\revision{For any set $\X \subseteq \{0,1\}^n$, $\conv\bigl(\Pi(\X)\bigr) = H \cap \conv(K)$ and $\conv(K) = \bigl\{(y_0,\y)\bigm| \y \in y_0 \cdot \conv(\X),\ y_0 > 0 \bigr\}$.}
\end{lemma}

In the remainder of this subsection, we formally discuss the technical implication  of Lemma~\ref{lemma:fractional-LP} when a compact extended formulation of $\conv(\X)$ is given. Here, we say a system of linear inequalities that defines a polyhedron $Q:=\{(\x,\z) \in \R^{n+q} \mid \BFA\x + \BFB \z \leq \BFb\}$ is an \textit{extended formulation} of $\conv(\X)$ if $ \conv(\X)=\proj_{\x}(Q):= \bigl\{ \x \bigm| \exists \z \in \R^q \text{ s.t. } (\x,\z) \in Q\bigr\}$.
Moreover, we say that the extended formulation is \textit{compact} if both $q$ and the number of inequalities defining $Q$ is polynomial in $n$.

\begin{theorem}\label{them:extended-formulation}
    Let $\BFA \x +\BFB\boldsymbol{z} \leq \BFb$ be an extended formulation for  $\conv(\X)$. Then, an  extended formulation for $\conv\bigl(\Pi(\X)\bigr)$ is $\bigl\{(y_0,\y,\boldsymbol{z}) \bigm|  \BFA\y + \BFB \boldsymbol{z} \leq \BFb y_0,\ (y_0,\y) \in H,\ y_0 \geq 0\bigr\}$. Moreover, if $\conv(\X)$ admits a compact extended formulation then so does $\conv\bigl(\Pi(\X)\bigr)$. 
\end{theorem}
 \revision{When the feasible assortment set is defined by a totally unimodular matrix, its convex hull coincides with its LP relaxation~\citep{conforti2014integer}, recovering the results in~\cite{avadhanula2016tightness} and \cite{sumida2021revenue}.  Relying on the progress in polyhedral combinatorics~\citep{schrijver2003combinatorial,conforti2014integer}, one can explore the structure of $\X$ beyond the totally unimodular matrix. In the following remark, we utilize Theorem~\ref{them:extended-formulation} and a compact extended formulation of the chain polytope to handle the two-stage Luce model (2SLM) introduced in Section~\ref{section:2SLM}}. 
\begin{remark}\label{rmk:2slm-formulation}
	\revision{In this remark, we discuss the case where an online consumer of type $i$ chooses according to 2SLM. In this case, we need to describe the convex hull of $\Pi_i(\X_i)$, where $\X_i$ is the set of undominated products given as in~\eqref{eq:undominated}. It is known that the convex hull of all undominated products $\X_i$ is called the \textit{chain polytope}~\citep{stanley1986two}. The continuous relaxation of~\eqref{eq:undominated} yields a polyhedral description for the chain polytope, but its size is exponential.}
    
    Recently, Corollary 2.9 in~\cite{fawzi2022lifting} gives an extended formulation for the chain polytope with $\mathcal{O}(n^2)$ inequalities  and $n$  additional variables. 
This formulation, together with Theorem~\ref{them:extended-formulation}, gives a compact extended formulation for the convex hull of the choice probability set under 2SLM. To provide an explicit description, we introduce few notions on a partial order. Recall that $\mathcal{P}_i$ is a partial order set specifying consumers' preference on the set of $n$ products, where the partial order is denoted by $\preceq_i$. We say that $a_2$ \textit{covers} $a_1$ if $a_1 \prec_i a_2 $ and there is no $b \in \mathcal{P}_i$ with $a_1 \prec_i b \prec_i a_2$, and we call this a \textit{cover relation}. 
Then, we obtain that $(y_{i0},\y_i) \in \R^n$ belongs to  the convex hull of $\Pi_i(\X_i)$ if there exists $\z_i \in \R^n$ such that
\begin{equation}\label{eq:2stageluce}
\begin{aligned}
	0 \leq z_{ia_1} \leq z_{ia_2} \leq  y_{i0} &\qquad \text{for all } a_2\succeq_i a_1 \text{ in } \mathcal{P}_i \\
	0 \leq y_{ia_2} \leq z_{ia_2} - z_{ia_1}& \qquad \text{for all cover relation } a_2 \succ_i a_1 \\
	y_{ia} = z_{ia}& \qquad \text{for all minimal element } a \in \mathcal{P}_i \\
 u_0y_{i0} + \sum_{j \in N} u_jy_{ij} = 1& \quad \text{and} \quad y_{i0} \geq 0. 
\end{aligned}\tag{\textsc{2-Stage-Luce}}
\end{equation}
This formulation reduces to the formulation of the chain polytope given by Corollary 2.9 of~\cite{fawzi2022lifting} when $y_{i0}=1$ and the normalization plane is removed.  \hfill \Halmos
\end{remark}


To conclude this subsection, we note that~\cite{Megiddo1979} constructs a fast combinatorial algorithm for connecting the complexity of optimizing a linear function over $\X$ and over $\Pi(\X)$. Specifically, maximizing a linear function over $\Pi(\X)$ is solvable in $\mathcal{O}\bigl(p(n)\cdot(p(n)+q(n))\bigr)$ if maximizing a linear function over $\X$ is solvable within $\mathcal{O}\bigl(p(n)\bigr)$ comparisons and $\mathcal{O}\bigl(q(n)\bigr)$ additions. While this combinatorial algorithm is sufficient to solve the stand-alone assortment problem with a single consumer type, it cannot be directly used to solve quick commerce assortment planning. In contrast, as demonstrated in Theorems~\ref{theorem:partial-convexification} and~\ref{them:extended-formulation}, the compact extended formulation  for $\conv(\X)$ can be easily incorporated into the MIP formulation of \eqref{eq:omni-rlx}, which can be directly solved by modern commercial solvers. \revision{For example, although, leveraging the Megiddo's algorithm and the polynomial-time solvability of the maximum weighted independent set over comparability graphs,~\cite{flores2017assortment} derive a polynomial-time algorithm for assortment optimization under 2SLM, it is not evident that this algorithm can be adopted to solve~\eqref{eq:omni-mnl}. }

\subsection{Convex hull results on~\eqref{eq:off-bi} and~\eqref{eq:link-bi}}\label{section:conv-bilinear}
After establishing tractable formulations for the convex hull of the choice probability set, we next derive an MIP formulation for~\eqref{eq:off-bi} and~\eqref{eq:link-bi}. To do so, we need to devise an efficient method to deal with the constraints related to bilinear terms. A prevalent method is to use the McCormick inequalities~\citep{mccormick1976computability}, which, in the context of quick commerce assortment planning, reformulates~\eqref{eq:off-bi} and~\eqref{eq:link-bi} as follows:
\begin{equation}\label{eq:mccormick}
\begin{aligned}
y_{ij} &\leq  y_{i0}^U \cdot x_j \qquad && y_{ij} \leq  y_{i0}^L \cdot x_j + y_{i0} - y_{i0}^L \quad \text{ for } i \in M^+ \text{ and } j \in N \\
    y_{0j} &\geq y_{00}^L \cdot x_j \qquad  && y_{0j} \geq y_{00}^U \cdot x_j + y_0 - y_{i0}^U  \quad \text{ for } j \in N,
\end{aligned}\tag{\textsc{McCormick}}
\end{equation}
where $y_{i0}^L$ (resp. $y_{i0}^U$) is a constant lower (resp. upper) bound on the no-purchase probability variable $y_{i0}$ of consumer type $i$. However,~\eqref{eq:mccormick} ignores the dependency between offline assortment decisions 
and choice probabilities of offline/online consumer segments. It is precisely this interdependence that motivates our convex hull study. In Section~\ref{section:bilinearsetup}, we formally define the nonconvex sets that are studied, and, at the end of Section~\ref{section:bilinearsetup}, we overview our theoretical findings.

\subsubsection{Motivation: tractable and intractable results}\label{section:bilinearsetup}
To streamline the presentation of Section~\ref{section:conv-bilinear}, we remove the subscript $i$. Moreover, we assume that there are no constraints on the feasible assortments, that is $\X = \{0,1\}^n$. Thus, we can use a compact notation $Y$ to denote the convex hull of a choice probability set, that is,  
\[
Y = \conv\bigl( \Pi(\{0,1\}^n) \bigr) = \biggl\{(y_0,\y)\biggm| u_0y_0 + \sum_{j \in N}u_jy_j=1,\ 0\leq y_j\leq y_0 \, \text{ for all } j \in N\biggr\}. 
\]
Now, (\ref{eq:off-bi}) and~(\ref{eq:link-bi}) lead us to study two mixed integer bilinear sets
\[
\B:= \biggl\{(\x, y_0, \y) \in \{0,1 \}^n \times Y \biggm| y_j = y_0x_j \text{ for all } j \in N \biggr\}
\]
and 
\[
\B^\leq := \biggl\{(\x, y_0, \y) \in \{0,1\}^n \times Y \biggm|  y_j \leq y_0x_j \text{ for all } j \in N \biggr\}, \\
\]
respectively. However, the following complexity result indicates that it is impossible to obtain a tractable description for $\conv(\B)$ and $\conv(\B^\leq)$.
\begin{proposition}\label{prop:hullcomplexity}
	The separation problems of $\conv(\B)$ and $\conv(\B^\leq)$ are NP-hard.
\end{proposition}

Thus, instead of studying $\conv(\B)$ and $\conv(\B^\leq)$, we opt to focus on a mixed integer bilinear set that is specifically associated with a single product $j$, that is
\[
\B_j :=  \bigl\{(x_j, y_0, \y) \in \{0,1\} \times Y \bigm|  y_j = y_0x_j\bigr\}.
\]
We also consider a decomposition of $\B_j$, that is $\B_j^\leq \cap \B_j^\geq$, where
\[
\B^\leq_j := \bigl\{(x_j, y_0, \y) \in  \{0,1 \} \times Y \bigm|  y_j \leq y_0x_j\bigr\} \quad \text{and} \quad \B^\geq_j := \bigl\{(x_j, y_0, \y) \in  \{0,1 \} \times Y\bigm|  y_j \geq y_0x_j\bigr\}.
\]
These sets are relevant to our study since their convex hulls provide a partial description of the convex hull of $\B$ and $\B^\leq$. In the following two remarks, we discuss relevant convexification results on single product bilinear sets.
\begin{remark}\label{rmk:disjunctive}
The convex hull of $\B_j$  can be obtained by using disjunctive programming~\citep{balas1998disjunctive} 
 since $\B_j$ is expressible as the union of two polytopes: 
\[
\B_j \cap \{(x_j,y_0,\y) \mid x_j =0 \} \qquad \text{and} \qquad \B_j \cap \{(x_j,y_0,\y) \mid x_j = 1 \}.
\]
This approach requires a convex multiplier and a copy of $(x_i,y_0,\y)$ variables for each disjunction, thus introducing $2(n+2) + 1$ additional continuous variables.
Hence, using the disjunctive programming approach needs to introduce $(m n +n)\cdot \bigl(2(n+2) +1\bigr)$ continuous variables since there are $mn+n$ single product bilinear terms in~(\ref{eq:off-bi}) and~(\ref{eq:link-bi}). \hfill \Halmos
\end{remark}
\begin{remark}\label{rmk:conic}
 \cite{sen2018conic} propose a novel mixed-integer conic formulation for assortment planning under the MMNL model. This approach yields a convex relaxation for  $\B_j^\geq$ given as follows: 
\begin{equation*}\label{eq:conic}
    \Bigl\{ (x_j,y_0,\y) \Bigm| wy_j  \geq  x^2_j\; \text{ for } j \in N ,\  w y_0 \geq 1 ,\ w = u_0 + \sum_{j\in N}u_jx_j,\ (y_0,\y) \in Y ,\ \x \in [0,1]^n \Bigr\},
\end{equation*}
by exploiting the relation $x^2_j = x_j$ for binary variable $x_j$. 
However, this relaxation strictly contains the convex hull of $\B_j^{\geq}$. To see this, we consider the following example where $n=2$, and $(u_0,u_1,u_2) = (1,1,2)$. To depict its convex hull in 3D, we consider $\B_1^\geq$ in the space of $(x_1,y_1,y_2)$ variables:
\[
\B_1^\geq = \Bigl\{(x_1,y_1,y_2) \Bigm| 0\leq y_1 \leq 1-y_1-2y_2,\ 0\leq y_2 \leq 1-y_1-2y_2,\ y_1 \geq x_1(1-y_1-2y_2),\ x_1 \in \{0,1\} \Bigr\}.
\]
The convex hull of $\B^\geq_1$, depicted as a gray polytope in Figure~\ref{fig:convexhull},  is given by 
\[
\begin{aligned}
\Bigl\{(x_1,y_1,y_2) \Bigm| 0\leq y_1 \leq 1-y_1-2y_2,\ 0\leq y_2 &\leq 1-y_1-2y_2,\ 0\leq  x_1 \leq 1,\ \\
&2y_1+2y_2-x_1\geq 0,\ 4y_1 -x_1 \geq 0  \Bigr\},
\end{aligned}
\]
where the last two inequalities are derived using~(\ref{eq:stair-under}) in Theorem~\ref{them:staircase},  and, in Theorem~\ref{them:convexhull}, we show that they indeed yield a convex hull description for $\B_1^\geq$.  In contrast, the convex relaxation from~\cite{sen2018conic} is given as follows:
\[
\begin{aligned}
\bigl\{(x_1,y_1,y_2) \bigm| wy_1 &\geq x_1^2,\  wy_2 \geq x_2^2,\ wy_0 \geq 1 ,\ w = 1+x_1+2x_2,\ \\
&0\leq y_1 \leq y_0,\ 0\leq y_2 \leq y_0,\ y_0 + y_1+2y_2 = 1,\  \x \in [0,1]^2 \bigr\}.    
\end{aligned}
\]
This relaxation contains the red point $(\bar{x}_1,\bar{y}_1,\bar{y}_2) = (0.95,0.25,0.2)$ in Figure~\ref{fig:convexhull} since it can be extended to a point $(\bar{y}_0,\bar{y}_1,\bar{y}_2) = (0.35,0.25,0.2)$, $\bar{w} = 3.61$, $(\bar{x}_1, \bar{x}_2) = (0.95, 0.83)$ which satisfies constraints in the relaxation. However, the red point does not belong to the gray polytope as it violates one of the defining inequalities $2y_1+2y_2-x_1 \geq 0$ and its projection of the plane $2y_1+2y_2-x_1 = 0$ is depicted as the blue point in Figure~\ref{fig:convexhull}. \hfill \Halmos
\end{remark} 

\begin{figure}
\centering
  \tdplotsetmaincoords{80}{330}

  \begin{tikzpicture}[scale=2, color={lightgray},tdplot_main_coords]
    \definecolor{point_color}{rgb}{ 0,0,0 }
    \tikzstyle{point_style} = [fill=point_color]

    \coordinate (p00) at (0, 0, 0);
    \coordinate (p01) at (0, 1/3*2, 0);
     \coordinate (p10) at (1/2*2, 0, 0);
     \coordinate (p11) at (1/4*2, 1/4*2,0);
    \coordinate (q10) at (1/2*2, 0, 1);
    \coordinate (q11) at (1/4*2, 1/4*2,1);
    \coordinate (w) at (1/4*2,1/5*2,0.95);
    \coordinate (p) at (1/4*2.2,1/5*2.2,0.8);



 \filldraw[blue] (p00) circle (0.1pt) ;  \filldraw[black] (p01) circle (0.1pt) ; 
  \filldraw[black] (p10) circle (0.1pt) ;  \filldraw[black] (p11) circle (0.1pt) ; 
 \filldraw[blue] (q10) circle (0.1pt) ; 
 \filldraw[blue] (q11) circle (0.1pt) ; 





    \definecolor{edge_color}{RGB}{100,100,100}

    \definecolor{domain_color}{RGB}{100,100,100}

    \tikzstyle{domain_style} = [fill=domain_color, fill opacity=0.3, draw=domain_color, line width=1 pt, line cap=round, line join=round]
    \tikzstyle{domain+_style} = [fill=domain_color, fill opacity=0.2, draw=domain_color, line width=1 pt, line cap=round, line join=round]

    
   \definecolor{conv_color}{RGB}{255,217,0}

    \tikzstyle{conv_style} = [fill=conv_color, fill opacity=0.3, draw=conv_color, line width=1 pt, line cap=round, line join=round] 


   \definecolor{conc_color}{RGB}{0,95,255}

    \tikzstyle{conc_style} = [fill=lightgray, fill opacity=0.3, draw=black, line width=0.1 pt, line cap=round, line join=round] 

    \draw[fill=lightgray, fill opacity=0.3, draw=lightgray,dashed, line width=0.1 pt] (p00) --(p01) -- (p11) -- (p10) -- (p00) -- cycle;
    \draw[fill=lightgray, fill opacity=0.3, draw=lightgray,dashed, line width=0.1 pt] (q11) --(q10) -- (p10) -- (p11) -- (q11) -- cycle;
    \draw[fill=lightgray, fill opacity=0.3, draw=lightgray,dashed, line width=0.1 pt] (q11) -- (p11) -- (p01) -- (q11) -- cycle;
    \draw[conc_style] (p00) -- (p10) -- (q10) -- (p00) -- cycle;
    \draw[conc_style] (p00) -- (p01) -- (q11) -- (p00) -- cycle;
    \draw[conc_style] (p00) -- (q10) -- (q11) -- (p00) -- cycle;

\draw[black,->] (-0.2*2,-0.3*2,0) -- (0.6*2,-0.3*2,0) node [anchor=north]{$y_{1}$};
\draw[black,->] (-0.2*2,-0.3*2,0) -- (-0.2*2,0.6*2,0) node [anchor=north]{$y_{2}$};
\draw[black,->] (0.6*2,-0.3*2,0) -- (0.6*2,-0.3*2,1.5) node [anchor=west]{$x_{1}$};

 \filldraw[blue] (q10) circle (0.1pt) ; 
 \filldraw[blue] (q11) circle (0.1pt) ; 
\filldraw[blue] (p00) circle (0.1pt) ;
\draw[black,thick] (p) -- (w);
     \filldraw[red] (w) circle (0.8pt) ; 
     \filldraw[blue] (p) circle (0.8pt) ; 




  \end{tikzpicture}
    \caption{The conic quadratic relaxation in~\cite{sen2018conic} fails to describe $\conv(\B_1^\geq)$}
    \label{fig:convexhull}
\end{figure}
On the one hand, disjunctive programming provides a convex hull description of $\B_j$ using $2(n+2)+1$ additional continuous variables.  On the other hand, while the state-of-the-art conic formulation uses few additional variables to obtain a strong relaxation for $\B^\geq_j$, it fails to describe the convex hull for $\B_j$. 
It seems that there is a trade-off between size and tightness in describing the convex hull of $\B_j$.  
In fact, the remainder of Section~\ref{section:conv-bilinear} aims to demonstrate the existence of explicit  convex hull description of $\B_j$ without the introduction of new variables. In particular, in Section~\ref{section:ineq}, we present a family of valid linear inequalities,~(\ref{eq:stair-over}) and~(\ref{eq:stair-under}), for the single product set $\B_j$. Although the number of inequalities is exponential in the number of products, we develop a modified choice probability-ordered policy in Section~\ref{section:sep} to solve its separation problem. Then, we leverage the policy to establish that 
\[
\textit{$(\ref{eq:stair-over})$ and~$(\ref{eq:stair-under})$ describe $\conv(\B_j)$ without using additional variables.}
\]



\subsubsection{Single product relaxations}\label{section:ineq}
\def \ext {\text{ext}}

To streamline the presentation,  we define 
\[
U(S) := u_0+\sum_{i \in S}u_i \quad \text{and} \quad \alpha(S) := \frac{1}{U(S)} \quad \text{ for every } S \subseteq N.
\]
In Theorem~\ref{them:staircase}, we use a constructive procedure to derive valid linear inequalities,~(\ref{eq:stair-over}) and~(\ref{eq:stair-under}), for the single product set $\B_j$. Our rationale of providing a constructive proof is that it may be useful in generating valid inequalities for similar structures. 
\revision{Before proceeding to the formal description, we illustrate the basis idea by deriving one of inequalities defining the polytope in Figure~\ref{fig:convexhull}. A similar argument can be used to derive the convex hull of $\B_j^\leq$.}
\begin{example}\label{ex:inequalities}
In this example, we consider the case where $n=2$ and $(u_0,u_1,u_2) = (1,1,2)$ again, and derive one of the inequalities defining the polytope in Figure~\ref{fig:convexhull}. Recall that we are interested in deriving the convex hull of $\B_1^\geq= \bigl\{(x_1, y_0, y_1,y_2) \in \{0,1\} \times  Y \bigm| y_1 \geq x_1  y_0\bigr\}$, where $Y=\{(y_0,y_1,y_2) \mid y_0 + y_1 +2y_2 = 1,\ 0 \leq y_1 \leq y_0,\ 0\leq y_2 \leq y_0\}$. Our procedure exploits structures in vertices of $Y$, which are enumerated as follows:
\[
\ext(\emptyset) = (1,0, 0) \quad \ext(\{1\}) = \Bigl(\frac{1}{2},\frac{1}{2},0\Bigr) \quad  \ext(\{2\}) = \Bigl(\frac{1}{3},0,\frac{1}{3}\Bigr) \text{ and } \quad \ext(\{1,2\}) = \Bigl(\frac{1}{4},\frac{1}{4},\frac{1}{4}\Bigr),
\]
where, in addition to satisfying the normalization condition, the $j^{\text{th}}$ coordinate of an extreme point is either zero or equal to the $0^{\text{th}}$ coordinate, for example, $\ext(\{1\})$ means that only the first coordinate is nonzero.  Our procedure starts by using an extreme point, for example $\ext(\{1\})$, to define a truncated function $\tau: Y \to \R$, that is, $\tau(y_0,y_1,y_2) := \max \{y_0,\frac{1}{2}\}$, which can be interpreted as the range of $y_0$ provided that products in $\{1\}$ are available to be selected. Next, we use the truncated function to obtain an underestimating function  of  $x_1y_0$ over $\{0,1\}\times Y$ as follows:
\[
\begin{aligned}
    y_1 &\geq x_1 y_0  = x_1 \cdot\tau(y_0,y_1,y_2)  + x_1 \bigl(y_0 - \tau(y_0,y_1,y_2)\bigr) \geq x_1 \cdot \frac{1}{2} + y_0 - \tau(y_0,y_1,y_2) \geq x_1 \cdot \frac{1}{2} - y_2,
\end{aligned}
\]
where the equality expands $x_1 y_0 $ using the truncated function as an intermediate step,  the second inequality holds since the two terms are relaxed individually. For the first term, the truncated function is relaxed to its lower bound $\frac{1}{2}$, and for the second term, $x_1$ is relaxed to its upper bound $1$. The last inequality holds since $-y_2$ is a linear underestimating function of $y_0 - \tau(y_0,y_1,y_2)$. This relaxation is  performed by linearly interpolating the concave function $(1-y_1 - 2y_2) -\tau(y_0,y_1,y_2)$ over the extreme point $\ext(\{1\})$ and its neighborhood $\ext(\{1,2\})$. \hfill \Halmos
\end{example}

\begin{theorem}\label{them:staircase}
Fix  $j \in N$. For $S \subseteq N \setminus \{j\}$, the following are valid for $\B_j$
\begin{align}
        y_j    &\geq \alpha(S\cup j) \cdot x_j  - \sum_{t \in  N \setminus (S \cup \{j\})} u_{t} \cdot  \alpha(S\cup j) \cdot y_t\tag{\textsc{Under}}\label{eq:stair-under} \\
	y_j &\leq    \alpha(S\cup j) \cdot x_j + \bigl(1-(u_0+u_j) \cdot \alpha(S\cup j) \bigr) \cdot y_0 - \sum_{t \in S } u_t \cdot \alpha(S\cup j)  \cdot y_t \tag{\textsc{Over}} \label{eq:stair-over}.
\end{align}
\end{theorem}
Next, we specialize \eqref{eq:stair-over} and \eqref{eq:stair-under} to the case when $S \in \{\emptyset, N\setminus \{j\}\}$. The resulting four inequalities imply McCormick inequalities, thus yielding an MIP formulation of~\eqref{eq:off-bi} and~\eqref{eq:link-bi} which is tighter than that based on~\eqref{eq:mccormick}.

\begin{corollary}\label{cor:mccormick}
    Let $j \in N$ be a given product. For $S \in \{\emptyset, N\setminus \{j\}\}$,~\eqref{eq:stair-under} and \eqref{eq:stair-over} yield
    \begin{equation}\label{eq:improvedmccormick}
    \begin{aligned}
    y_j &\geq \alpha(N) \cdot x_j \qquad  &&   y_j \geq \alpha(\emptyset) \cdot x_j  + y_0 - \alpha(\emptyset) \\
           y_j &\leq  \alpha(j) \cdot x_j \qquad && y_j \leq  \alpha(N \setminus j ) \cdot x_j + y_0  -  \alpha(N\setminus j ),
\end{aligned}
\tag{\textsc{McCormickPlus}}
\end{equation}
which implies McCormick inequalities for the bilinear set $\B_j$. 
\end{corollary}

\subsubsection{\revision{Fast separation algorithm and convex hull characterization}}\label{section:sep}
\rvtwo{Note that the number of inequalities is exponential in the number of products.} In order to leverage these inequalities to speed up the computation time for \eqref{eq:omni-rlx}, we need to design an efficient algorithm to find an inequality, out of the exponential many, to separate an infeasible solution. \revision{In particular, we propose a modified choice probability-order policy to solve the  separation problem with a complexity of $\mathcal{O}(n\log n)$. More interestingly, such policy will be used to show that the single product relaxation defined by~(\ref{eq:stair-under}) and~(\ref{eq:stair-over}) is theoretically optimal in a sense that it coincide with the convex hull of $\mathcal{B}_j$.}

Given a vector  $(\bar{y}_0,\bar{\y}) \in Y$ and $\bar{x} \in [0,1]$, finding a violated inequality in~(\ref{eq:stair-under}) and~(\ref{eq:stair-over}) is equivalent to solving the following combinatorial optimization problems
\begin{align}
	&\max_{S \subseteq N \setminus \{j\}}\biggl\{\frac{\bar{x}_j}{U(S \cup j)}  - \frac{\sum_{t \notin  (S \cup j)} u_{t} \cdot   \bar{y}_t}{U(S \cup j)}  \biggr\} \tag{\textsc{Sep-Under}}\label{eq:sep-cover} ,\\
	 &\min_{S \subseteq N \setminus \{j\}} \biggl\{   \frac{\bar{x}_j}{U(S \cup j)} + \frac{\sum_{t \in S }  u_t \cdot (\bar{y}_0 - \bar{y}_t)}{U(S \cup j)} \biggr\} \tag{\textsc{Sep-Over}}\label{eq:sep-pack},
\end{align}
respectively. 
\rvtwo{The optimization problem defined in \eqref{eq:sep-cover} and \eqref{eq:sep-pack} can be interpreted as follows. We assume that a set $S \cup \{j\}$ is the ``imaginary'' assortment offered to consumers. If $S \cup \{j\}$ is indeed the actual offered assortment (i.e. $y_t$ is equal to $y_0$ for $t \in S \cup \{j\}$ and $0$ otherwise), then $\bar{y}_j$ should be exactly equal to  $\frac{\bar{x}_j}{U(S \cup j)}$. However, the actual assortment induced by a feasible solution $(\bar{x}_j, \bar{y}_0, \ldots, \bar{y}_j,\ldots ,\bar{y}_n)$ to $\B_j$ may not be $S \cup \{j\}$. Therefore, the value of $\bar{y}_t$ can deviate from $\frac{\bar{x}_j}{U(S \cup j)}$. If $\bar{y}_t$ is positive for $t \notin S \cup j$, it indicates that consumers choose products outside the imaginary assortment. Thus, the probability of choosing product $j$ needs to be adjusted downward. On the other hand, if $\bar{y}_t$ is strictly less than $\bar{y}_0$ for $t \in S$, it means that some products inside of the imaginary assortment are chosen less than expected. As a result, the probability of choosing product $j$ needs to be adjusted upward. Both adjustments are conducted according to the preference weights. The optimization problem of \eqref{eq:sep-cover} and \eqref{eq:sep-pack} is to seek the best imaginary assortment that provides the largest lower bound and the smallest upper bound.} 

Next, we present Algorithm~\ref{alg:sep}, 
a modified choice probability-order policy, to obtain the optimal assortments for \eqref{eq:sep-cover} and \eqref{eq:sep-pack}. \rvtwo{To highlight its connection with the prevalent revenue-ordered argument in assortment optimization, see~\citep{talluri2004revenue}, we prove the correctness of Algorithm~\ref{alg:sep} using a similar argument as in the proof of Theorem 5.1 of~\cite{gallego2019revenue}}.

\begin{algorithm}
\fontsize{9pt}{10pt}\selectfont 
\renewcommand{\baselinestretch}{1.3}\selectfont
\caption{Separation for~(\ref{eq:stair-under}) and~(\ref{eq:stair-over})}\label{alg:sep}
\KwData{$\bar{x}_j \in [0,1]$ and $(\bar{y}_0,\bar{\y}) \in Y$}
\KwResult{either $(\bar{x}_j,\bar{y}_0,\bar{\y})$ satisfies~(\ref{eq:stair-under}) and~(\ref{eq:stair-over}) or a violated inequality}
find a permutation $\sigma$ of $N \setminus \{j\}$ so that $\bar{y}_{\sigma(1)} \geq \cdots \geq \bar{y}_{\sigma(n-1)}$\;
$S^* \gets \bigl\{\sigma(1), \ldots,\sigma(k) \bigm| \bar{y}_{\sigma(k)} \geq \bar{y}_j \bigr\}$\;
$T^* \gets \bigl\{\sigma(1), \ldots,\sigma(k) \bigm| \bar{y}_{\sigma(k)} \geq \bar{y}_0-\bar{y}_j\bigl\}$\;
\eIf{$(\bar{x}_i,\bar{y}_0,\bar{\y})$ satisfies~$(\ref{eq:stair-under})$ with $S=S^*$ and~$(\ref{eq:stair-under})$ with $S = T^*$ }{
	$(\bar{x}_i,\bar{y}_0,\bar{\y})$ satisfies~$(\ref{eq:stair-under})$ and~$(\ref{eq:stair-over})$\; 
  }{
	either~$(\ref{eq:stair-under})$ with $S=S^*$ or~$(\ref{eq:stair-over})$ with $S=T^*$ 
violates at point	$(\bar{x}_j,\bar{y}_0,\bar{\y})$\;
  } 
\end{algorithm}



\begin{theorem}\label{them:poly-sep-stair}
\revision{Algorithm~\ref{alg:sep} separates~$(\ref{eq:stair-under})$ and~$(\ref{eq:stair-over})$ in $\mathcal{O}(n\log n)$. }
\end{theorem}


\revision{Last, we investigate the theoretical strength of our proposed inequalities. In particular, we show that~(\ref{eq:stair-under}) and~(\ref{eq:stair-over}) yield a convex hull description for a single product bilinear set $\B_j$. The key step of the proof is to show that for a given point satisfying $\B_j$, we need to decompose into a convex combination of points in $\B_j$. This decomposition is found by using the modified choice probability-ordered separation algorithm.}

\begin{theorem}\label{them:convexhull} 
\revision{For a given product $j$, ~$(\ref{eq:stair-under})$ (resp. ~$(\ref{eq:stair-over})$) describes the convex hull of $\mathcal{B}_j^\geq$ (resp. $\mathcal{B}_j^\leq$). Moreover,~$(\ref{eq:stair-under})$ and~$(\ref{eq:stair-over})$ describe the convex hull of $\B_j$.}  
\end{theorem}

\section{\revision{Explicit formulations and a cutting-plane implementation}}\label{section:CH}
\revision{In this section, we leverage the convex hull descriptions derived in Section~\ref{section:convexhull} to present explicit reformulations for~\eqref{eq:omni-mnl}. First, we present a formulation for constrained~\eqref{eq:omni-mnl}, referred to as~\textsf{CH} \revision{(convex hull)}, and then, we specialize this formulation to the case where the choice behavior of online consumers follows the two-stage Luce model, referred to as~\textsf{CH-Chain}.} 
Last,  we delve into a discussion on the implementation aspect, wherein we propose a cutting-plane algorithm. \revision{With the help of the separation algorithm, Algorithm~\ref{alg:sep}, the cutting-plane approach is easy to implement and computationally efficient.}

We start by presenting formulation \textsf{CH} for \revision{constrained}~\eqref{eq:omni-mnl} with the MNL model. Recall that  Theorem~\ref{theorem:partial-convexification} provides an alternative formulation~\eqref{eq:omni-rlx} of~\eqref{eq:omni-mnl}.  This reformulation consists of three types of nonconvex constraints, \eqref{eq:on-mnl}, \eqref{eq:off-bi} and \eqref{eq:link-bi}, which can be convexified using techniques in Sections~\ref{section:formulation-MNL} and~\ref{section:conv-bilinear}. \revision{In addition, for the offline channel, we assume that the feasible assortment set $\X_0$ can be represented as $\bigl\{\x \in \{0,1\}^n \bigm| \BFA_0 \x \leq \BFb_0\bigr\}$, and, for each online consumer type $i \in M$, we assume that the convex hull of the operational constraints $\X_i$ admits an extended polyhedral representation, that is,}
\[
\revision{\conv(\X_i) = \bigl\{\x_i \bigm| \exists \boldsymbol{z}_i \text{ s.t. } \BFA_i\x_i + \BFB_i\z_i \leq \BFb_i \bigr\},}
\]
\revision{where $\BFA_i$ and $\BFB_i$ (resp. $\BFb_i$) is a matrix (resp. vector) with a proper dimension.} These developments lead  to the following equivalent formulation of~\eqref{eq:omni-mnl}, which we refer to as formulation~\our
\begin{subequations}:
\begin{align}
\our:\; \max \quad   &  \sum_{i \in M^+} \alpha_i\Bigl( \sum_{j \in N} r_{ij} u_{ij}y_{ij}  \Bigr) \notag \\
	\text{s.t.} \quad & \x \in \{0,1\}^n,\ \revision{\BFA_0\x \leq \BFb_0},\ u_{00}y_{00} + \sum_{j \in N}u_{0j}y_{0j} =1,\ \BFA_0\y_0 \leq \BFb_0y_{00} \label{eq:CH-0} \\
 &u_{i0} y_{i0} + \sum_{j \in  N} u_{ij} y_{ij} = 1 \text{ and }  \revision{\BFA_i\y_i+\BFB_i\boldsymbol{z}_i \leq \BFb_iy_{i0}} \qquad  \text{for  } i\in M \text{ and } j \in N \label{eq:CH-1}\\ 
&y_{0j}    \geq \alpha_0(S\cup j)  x_j  - \sum_{t \in  N \setminus (S \cup j)} u_{0t}   \alpha_0(S\cup j)  y_{0t}  
\qquad  \text{for }  j \in N  \text{ and } S \subseteq  N \setminus j \label{eq:CH-2} \\
&y_{ij} \leq    \alpha_i(S\cup j) x_{j} + \bigl(1-(u_{i0}+u_{ij}) \alpha_i(S\cup j) \bigr) y_{i0} - \sum_{t \in S } u_{it} \alpha_i(S\cup j)   y_{it} \notag \\
& \qquad \qquad \qquad \qquad \qquad\qquad \text{for } i \in M^+,\ j \in N \text{ and } S\subseteq N \setminus j \label{eq:CH-3} \\
&(u_{i0} + \sum_{j \in N}u_{ij}x_j ) y_{i0} \geq 1 \qquad \qquad  \qquad \qquad\qquad\qquad \text{ for } i \in M^+ \label{eq:CH-4},
\end{align} 
\end{subequations}
where $\alpha_i(S) := \frac{1}{\sum_{j\in \{0\}\cup S}u_{ij}}$ for each $i \in M^+$ and $S \subseteq N$ and implied constraints on offline choice probability  variables. \revision{The constraint~\eqref{eq:CH-0} models the offline feasible assortment set $\X_0$. It follows from Theorem~\ref{them:extended-formulation} that constraint~\eqref{eq:CH-1} describes the convex hull of the choice probability set $\Pi(\X_i)$.} The constraint~\eqref{eq:CH-2} is derived by applying~\eqref{eq:stair-under} in Theorem~\ref{them:staircase} to each product $j$ in the offline segment. The constraint~\eqref{eq:CH-3} is obtained by applying~\eqref{eq:stair-over} in Theorem~\ref{them:staircase} to each product $j$ cross all offline and online segments. We can interpret constraint~\eqref{eq:CH-2} (resp.~\eqref{eq:CH-3}) as polyhedral under-estimation (resp. over-estimation) of choice probability variable $y_{ij}$ of the $j^{\text{th}}$ product  in the $i^{\text{th}}$ segment. 
Note that we have not yet imposed additional constraints on the no-purchase probability variable $y_{i0}$ except the normalization requirement in~\eqref{eq:CH-1}. Next, we argue that it is natural to impose constraint~\eqref{eq:CH-4} as a convex under-estimation of the no-purchase probability variable. Recall that in formulation~\eqref{eq:omni-rlx}, ideally we wish to enforce the following  constraint
\[
y_{i0} = \frac{1}{u_{i0} + \sum_{j \in N}u_{ij} y_{ij}/y_{i0}} \quad \text{ for } i \in M^+,
\]
where $y_{ij}/y_{i0}$ models whether we offer product $j$ to consumer type $i$. This is nonconvex. However, we can relax this constraint as follows
\[
y_{i0} \geq \frac{1}{u_{i0} + \sum_{j \in N}u_{ij}y_{ij}/y_{i0}}  \geq \frac{1}{u_{i0} + \sum_{j \in N}u_{ij}x_{j}} \quad \text{ for } i \in M^+,
\]
where the second inequality holds since $y_{ij}/y_{i0} \leq x_j$ for all $i$ and $j$. If we now move the denominator of the right-hand side to the left, we arrive at constraint~\eqref{eq:CH-4}. 

Since the number of constraints in~\eqref{eq:CH-2} and~\eqref{eq:CH-3} is exponential in the number of products, even the continuous relaxation of $\textsf{CH}$ can not be efficiently solved by using a state-of-the-art solvers. To circumvent this difficulty, we propose a cutting plane algorithm, Algorithm~\ref{alg:cutting}, to \revision{implement formulation} $\textsf{CH}$. Instead of optimizing over all inequalities in~\eqref{eq:CH-2} and~\eqref{eq:CH-3}, Algorithm~\ref{alg:cutting} starts with processing a base formulation,  referred to as $\textsf{CH-0}$, which is constructed using inequalities in~\eqref{eq:improvedmccormick} instead of all inequalities in~\eqref{eq:CH-2} and~\eqref{eq:CH-3}
\begin{equation}
\begin{aligned}
\textsf{CH-0}: \max \quad &  \sum_{i \in M^+} \alpha_i\biggl( \sum_{j \in N} r_{ij} u_{ij}y_{ij}  \biggr) \notag\\
	\text{s.t.} \quad 	 
 	& \eqref{eq:CH-0},\ \eqref{eq:CH-1}  \\
               &(u_{i0} + \sum_{j \in N}u_{ij}x_j ) y_{i0} \geq 1 && \text{ for } i \in M^+ \\
    &y_{0j} \geq \alpha_0(N)  x_j \text{ and }  y_{0j} \geq \alpha_0(\emptyset)  x_j  + y_{00} - \alpha(\emptyset) && \text{ for } j \in N \\
           &y_{ij} \leq  \alpha_i(j) x_j \text{ and } y_{ij} \leq  \alpha_i(N \setminus j )  x_j + y_{i0}  -  \alpha(N\setminus j ) && \text{ for } i \in M^+, j \in N .
        \end{aligned}
\end{equation}
Next, we solve the continuous relaxation of $\textsf{CH-0}$ and obtain an optimal solution $(\bar{\x}, \bar{y}_0, \bar{\y})$. For each $i \in M^+$ and $j \in N$, we use the separation oracle, Algorithm~\ref{alg:sep}, to generate cuts, out of exponentially many ones in~\eqref{eq:CH-2} and~\eqref{eq:CH-3}, to cut off infeasible point $(\bar{x}_j, \bar{y}_{i0}, \bar{\y}_i)$ from the feasible region of $\textsf{CH}$. Adding newly generated cuts into the base formulation  $\textsf{CH-0}$ yields a tighter formulation $\textsf{CH-1}$. One can repeat this procedure $K$ times, and obtain a formulation, which is referred to as $\textsf{CH-K}$. \revision{Throughout the algorithm, the formulation is always exact as the base one is. Moreover, the cutting-plane algorithm generates a hierarchy of tighter formulations with the final one representing the formulation $\textsf{CH}$ in terms of the optimal objective value of the continuous relaxation~\citep{kelley1960cutting}}.
\begin{corollary}\label{cor:cutting-plane}
$\textsf{CH-0}$ \revision{is an exact formulation of~\eqref{eq:omni-mnl}. The optimal value of the continuous relaxation of $\textsf{CH-K}$ converges to that of $\textsf{CH}$ in a finite number of iterations }.
\end{corollary}

\begin{algorithm}\caption{A cutting-plane implementation of~\our }\label{alg:cutting}
\fontsize{9pt}{10pt}\selectfont 
\renewcommand{\baselinestretch}{1.3}\selectfont
\KwData{Formulation $\textsf{CH-0}$ and a positive integer $K$}
\KwResult{Formulation $\textsf{CH-K}$}
$k \gets 1$\;
\While{$k \le K$}
{
    $(\bar{\x}, \bar{y}_0, \bar{\y}) \gets $  an optimal solution to the continuous relaxation of $\textsf{CH-(k-1)}$\;
    $\textsf{Cuts} \gets \emptyset$ \;
    \For{$i$ in $M^+$}{
        \For{$j$ in $N$}{
        $\ell \gets$ call the separation oracle, Algorithm \ref{alg:sep}, to separate $(\bar{x}_j, \bar{y}_{i0}, \bar{\y}_i)$\;
        $\text{push}(\ell, \textsf{Cuts})$\;
        }
    }
    $\textsf{CH-k} \gets $ add $\textsf{Cuts}$ into formulation $\textsf{CH-(k-1)}$\;
    $k = k+1$\;
}
\end{algorithm}

\revision{
Next, we present formulation~\textsf{CH-Chain} for~\eqref{eq:omni-mnl} when consumers in each online segment make a choice decision according to the two-stage Luce model (2SLM), a choice model introduced in Section~\ref{section:2SLM}.
 In this setting, online consumers of type $i$ make the choice according to the MNL model based on the undominated products $\X_i$ defined as in~\eqref{eq:undominated}.
In other words, we need to derive a formulation for~\eqref{eq:omni-mnl} when $\X_i$ is~\eqref{eq:undominated} for each $i \in M$. Moreover, in Remark~\ref{rmk:2slm-formulation}, we show that $(y_{i0},\y_{i}) \in \conv(\Pi_i(\X_i))$ if and only if there exists $\z_i \in \R^n$ such that constraints in~\eqref{eq:2stageluce} are satisfied.
%
Recall that letting $y_{i0} = 1$, the first three constraints in~\eqref{eq:2stageluce} describe the \textit{chain polytope}, which is defined as the convex hull of~\eqref{eq:undominated}. Thus, the formulation of~(\ref{eq:omni-mnl}) obtained from replacing~\eqref{eq:CH-1} in formulation \textsf{CH} with~\eqref{eq:2stageluce} is referred to as formulation $\textsf{CH-Chain}$. Similarly, the formulation obtained from replacing \eqref{eq:CH-1} in formulation $\textsf{CH-0}$ with \eqref{eq:2stageluce} is referred to as $\textsf{CH-Chain-0}$, and the formulation generated by Algorithm \ref{alg:cutting} with input $\textsf{CH-Chain-0}$ is referred to as $\textsf{CH-Chain-K}$. }

\revision{To conclude this section, \revision{we remark that our formulations can also be solved  using other integer programming solution approaches~\citep{bertsimas1997introduction,conforti2014integer}. In particular, motivated by the success of Benders decomposition in solving large scale assortment optimization problems in~\cite{BertsimasMisic2019} and \cite{ChenMisic2021}, it is interesting to investigate the performance of such approach in our formulations.}}

\section{Computational experiments}\label{section:computation}

In this section, we present the results of our numerical experiments.  In particular, we report the performance of our formulations \textsf{CH} and \textsf{CH-Chain}, comparing them to other formulations using synthetic data.  \revision{In Appendix~\ref{sec:insight}, we offer managerial insights into the nature of~\eqref{eq:omni-mnl}.} All our numerical studies are performed in Python 3.10 on a Virtual Machine with 32 GB RAM and a 4-core Intel Core (Broadwell) @2.20 GHz processor.  All linear, second order conic and mixed-integer optimization problems were solved by using Gurobi 10.01 \citep{gurobi}. 



\subsection{\rvtwo{Alternative} formulations}
In our computational study, we consider two special cases of our model~\eqref{eq:omni-mnl}. For the unconstrained case, we consider three \rvtwo{alternative} formulations
\begin{itemize}
    \item[\textsf{CH-2}:] \revision{This is the formulation obtained by using Algorithm~\ref{alg:cutting} to implement \textsf{CH} with the number of iterations chosen to be $K = 2$. Our choice $K=2$ is based on the computational evidence that cutting-planes generated after two rounds almost close the remaining gap (see details in~\ref{EC:cutinground}).}
    \item[\conic:] This is a conic integer optimization formulation based on techniques in~\cite{sen2018conic}, and its derivation is detailed in Appendix \ref{section:conicformulation};
    \item[\milp:] This is a mixed-integer linear formulation based on the standard big-M linearization, and its derivation is detailed in Appendix \ref{section:milpformulation}.
\end{itemize}
For the case when customers in each online segment make a choice decision according to the two-stage Luce model, we consider
\begin{itemize}
    \item[\textsf{CH-Chain-2}:] \revision{This is the formulation obtained by using Algorithm~\ref{alg:cutting} to implement \textsf{CH-Chain} with the number of iterations chosen to be $K = 2$};
    \item[\conicC:] This is obtained by adding into \textsf{Conic} the linear constraints in \eqref{eq:undominated} that describe the set of undominated products of each online segment; 
    \item[\milpC:] This is obtained by adding into \textsf{MILP} the linear constraints in \eqref{eq:undominated}. 
\end{itemize}
\revision{Note that one can implement our inequalities~\eqref{eq:stair-under} and~\eqref{eq:stair-over} with \textsc{Gurobi}'s \textsc{Callback} function.  A computational result in~\ref{EC:callback} show that our implementation \textsf{CH-2} and \textsf{CH-Chain-2} outperforms the  \textsc{Callback} approach.}

\subsection{Problem setting}\label{subsec:num-setting}
We randomly generate data for the revenue and preference weight of each product across the offline and online segments. Details about data generation can be seen in Appendix \ref{subsection:generation_revenue_prefernce}. Generally, we set the revenue of each product for the offline consumer segment as a random variable of uniform distribution $U[10,20]$. We assume that the online segments are divided into two groups: the regular group and the VIP group. For each product, the revenue from the consumers in the regular group is the same as that from the offline consumer segment and the revenue from the consumers in the VIP group is lower than that from the regular group. We set the preference weight on each product across the offline and online consumer segments as a random variable from the uniform distribution $U[0,1]$. In addition, we fix the preference weight on the no-purchase option for the offline consumer segment to be $1$, while varying that for the online consumer segments in $\{2,5,10\}$. Furthermore, we assume that the probability of the offline consumer segment arriving in the system is $\alpha_0$ and that of each online segment $i\in M$ is $(1-\alpha_0)/m$ equally, i.e., $\sum_{i\in M}\alpha_i=1-\alpha_0$. 

We explore two choice models for the online segments: the standard MNL model and the 2SLM. The reason for considering the 2SLM exclusively for the online segments is that there is more data available in the online environment, which enables the firm to gain better insights into the consumer's choice behavior. Recall that to specify the 2SLM model, we need to use a partial order on products. In our experiments,  for each online segment, a partial order is randomly generated using a procedure in Appendix~\ref{subsection:generation_luce_graph}.

\subsection{Numerical results}\label{subsec:numerical_result}
\revision{We compare three different approaches in terms of the following performance measures: the average computation time (labeled as ``Time''), the minimal computation time (labeled as ``Min''), the maximal computation time (labeled as ``Max''), the standard deviation of the computation time (labeled as ``Std'')\footnote{The std is omitted (denoted by ``--'') if some instances are not solved successfully}, the average number of explored nodes by  Gurobi (labeled as ``Nds'') and the number of solved instances within a 3600 seconds time limit (labeled as ``\#''). Note that in our implementation of formulation \textsf{CH} and \textsf{CH-Chain}, we use the Algorithm \ref{alg:cutting} with $K=2$. Thus, the computation time of \textsf{CH-2} (resp. \textsf{CH-Chain-2}) includes the time \revision{spent by Algorithm \ref{alg:cutting} on cut generation} and the time for solving the formulation {\ourip} (resp. {\ouripC}).}



Table \ref{tab:Pref_NoLuce_Luce} reports the results for various configurations involving various problem sizes, that is, $(n,m) \in \{(100,50), (150,75), (200,100)\}$. \revision{For the preference weight values for the non-purchase option, we normalize that of the offline segment to $1$ and vary them of the online segments, that is, $u_{00}=1 $ and $u_{i0} \in \{2,5,10\},~ \forall i \in M$.} Each configuration includes 36 instances following the data generating process described in Section~\ref{subsec:num-setting} (see details in Appendix \ref{section:appendix_data_generation}), {and the offline consumer segment arrival probability for each instance is $\alpha_0 = 0.5$}.


\begin{table}[htbp]
\color{black}
        \fontsize{8pt}{10pt}\selectfont 
        \renewcommand\tabcolsep{1.8pt}
        \renewcommand{\arraystretch}{1.1}
        \centering
	\caption{Computational performance comparison }
	\label{tab:Pref_NoLuce_Luce}%
	\begin{threeparttable}   
		\begin{tabular}{p{4em}ccccccccccccccccccc}
			\toprule
   \multirow{2}{*}{$(n,m)$} & \multirow{2}{*}{$u_{i0}$} &
   \multicolumn{6}{c}{\textbf{CH-2}}&
   \multicolumn{6}{c}{\textbf{Conic}}&
   \multicolumn{6}{c}{\textbf{MILP}} \\
   \cmidrule(lr){3-8} \cmidrule(lr){9-14} \cmidrule(lr){15-20}
        & & Time & Min & Max & Std & Nds & \# 
        & Time & Min & Max & Std & Nds & \# 
        & Time & Min & Max & Std & Nds & \# \\
   \hline
   \multirow{3}{*}{$(100,50)$}
   & $2$ 
        & 8.2   & 6.0   & 13.7   & 1.7   & 1.0 & 36 
        & 75.8 & 5.7 & 240.2 & 48.7 & 14.4 & 36 
        & 46.8 & 3.4 & 175.9 & 42.2 & 411.2 & 36 
        \\
    & $5$
        & 15.2  & 6.6   & 55.2   & 11.1 & 1.0 & 36
        & 120.8 & 36.5 & 588.8 & 96.4 & 219.6 & 36
        & 131.8 & 17.1 & 611.4 & 131.9 & 3732.6 & 36
        \\
    & $10$
        & 28.7   & 10.0  & 80.7  & 14.5 & 3.0 & 36
        & 190.6 & 56.0 & 490.9 & 127.7 & 482.3 & 36
        & 663.6 & 50.6 & 3600 & -- & 19321.2 & 35
        \\
        
   \hline
   \multirow{3}{*}{$(150,75)$}
   & $2$ 
        & 21.5   & 15.1  & 38.6  & 4.2 & 1.0 & 36 
        & 347.5 & 121.0 & 845.9 & 146.0 & 44.1 & 36 
        & 251.2 & 39.4 & 544.5 & 132.3 & 577.0 & 36 
        \\
    & $5$
        & 48.3  & 23.1  & 173.8  & 31.8 & 3.0 & 36
        & 583.9 & 216.5 & 2034.9 & 430.9 & 506.3 & 36
        & 1190.9 & 141.2 & 3600 & -- & 10686.7 & 31
        \\
    & $10$
        & 191.7  & 38.4 & 378.8  & 89.4 & 15.8 & 36
        & 1040.3 & 346.5 & 2430.9 & 531.5 & 752.2 & 36
        & 3169.8 & 793.1 & 3600 & -- & 23726.2 & 10
        \\
        
   \hline
   \multirow{3}{*}{$(200,100)$}
   & $2$ 
        & 45.5   & 29.8  & 79.1  & 11.3 & 1.0 & 36 
        & 744.8 & 297.4 & 3600 & -- & 166.7 & 35
        & 329.2 & 85.9 & 3241.6 & 519.1 & 1493.4 & 36 
        \\
    & $5$
        & 196.8 & 59.7 & 3405.9  & 555.6 & 95.3 & 36
        & 2198.1 & 683.9 & 3600 & -- & 718.7 & 24
        & 3225.4 & 1658.2 & 3600 & -- & 12141.4 & 13
        \\
    & $10$
        & 544.5 & 111.4 & 2067.1 & 473.0 & 101.8 & 36
        & 3286.3 & 2024.8 & 3600 & -- & 892.4 & 17
        & 3600 & 3600 & 3600 & -- & 9841.4 & 0
        \\

\toprule

\multirow{3}{*}{$(n,m)$} & \multirow{3}{*}{$ u_{i0}$} &
   \multicolumn{6}{c}{\textbf{CH-Chain-2}}&
   \multicolumn{6}{c}{\textbf{Conic-Chain}}&
   \multicolumn{6}{c}{\textbf{MILP-Chain}} \\
   \cmidrule(lr){3-8} \cmidrule(lr){9-14} \cmidrule(lr){15-20}
        & & Time & Min & Max & Std & Nds & \# 
        & Time & Min & Max & Std & Nds & \# 
        & Time & Min & Max & Std & Nds & \# \\
   \hline
   \multirow{3}{*}{$(100,50)$}
   & $2$ 
        & 9.6   & 6.9   & 14.8   & 1.9 & 1.0 & 36 
        & 335.4 & 32.9 & 3600 & -- & 5234.4 & 35 
        & 3233.0 & 79.2 & 3600 & -- & 3.1e+5 & 4 
        \\
    & $5$
        & 17.8   & 9.2  & 45.0   & 8.0 & 1.4 & 36
        & 551.2 & 149.2 & 1250.6 & 227.7 & 2628.8 & 36 
        & 3600 & 3600 & 3600 & -- & 1.8e+5 & 0
        \\
    & $10$
        & 31.1   & 14.8  & 74.8  & 15.1 & 25.2 & 36
        & 709.4 & 290.7& 1874.0 & 345.1 &2201.3 & 36
        & 3600& 3600 & 3600 & -- & 1.1e+5 & 0
        \\
        
   \hline
   \multirow{3}{*}{$(150,75)$}
   & $2$ 
        & 31.8   & 22.5  & 43.0  & 5.0 & 1.0 & 36 
        & 2130.7 & 332.4 & 3600 & -- & 25308.7 & 28 
        & 3600 & 3600.0 & 3600 & -- & 1.1e+5 & 0
        \\
    & $5$
        & 69.3  & 33.0  & 172.8  & 34.1 & 10.2 & 36
        & 2948.4 & 853.8 & 3600 & -- & 12792.8 & 16
        & 3600 & 3600 & 3600 & -- & 45419.9 & 0
        \\
    & $10$
        & 144.8  & 60.9 & 307.2  & 80.7 & 38.2 & 36
        & 2837.1 & 1249.5 & 3600 & -- & 8029.9 & 25
        & 3600 & 3600 & 3600 & -- & 16591.3 & 0
        \\
        
   \hline
   \multirow{3}{*}{$(200,100)$}
   & $2$ 
        & 86.4  & 58.5  & 173.9  & 20.2 & 1.0 & 36 
        & 3510.1 & 2281.9 & 3600 & -- & 7338.2 & 5 
        & 3600 & 3600 & 3600 & -- & 25915.5 & 0 
        \\
    & $5$
        & 256.2 & 96.9 & 3116.7  & 494.4 & 109.7 & 36
        & 3600 & 3600 & 3600 & -- & 1546.8 & 0
        & 3600 & 3600 & 3600 & -- & 10538.1 & 0
        \\
    & $10$
        & 691.8 & 155.0 & 3557.7 & 638.9 & 162.9 & 36
        & 3600 & 3600 & 3600 & -- & 1507.6 & 0
        & 3600 & 3600 & 3600 & -- & 5246.1 & 0
        \\
\bottomrule
		\end{tabular}
	\end{threeparttable}

\end{table}%

\yun{Based on the results from Table \ref{tab:Pref_NoLuce_Luce}, it is evident that our formulation performs best in all the mentioned aspects for both choice models, MNL and 2SLM. Specifically, for the MNL, both {\ourip}  and {\conic} are able to solve all the $36$  instances with the default optimality gap $0.01\%$ and with a time limit of 3600s.} However, {\milp} even fails to solve an instance for $(n,m)=(100,50)$, not to mention the configuration with a larger problem size. For the 2SLM, {\ouripC} is able to solve all the instances with the average computation time being less than double that under MNL. This indicates the scalability of our approach as we explore the structure of  $\Pi(\X)$ effectively. However, the number of solved instances drops significantly for {\conic}, while {\milp} barely solves any instances for the 2SLM. The computational advantage of {\our} can be attributed to the fact that our formulation is tighter than the other two approaches, as indicated by the number of explored nodes under all three approaches. Especially when the preference weight of non-purchase for the online segment is small, it  takes only one node to solve {\ourip}.

\begin{figure}[htbp]
    \centering
 \includegraphics[width=0.4\textwidth,scale=0.8]{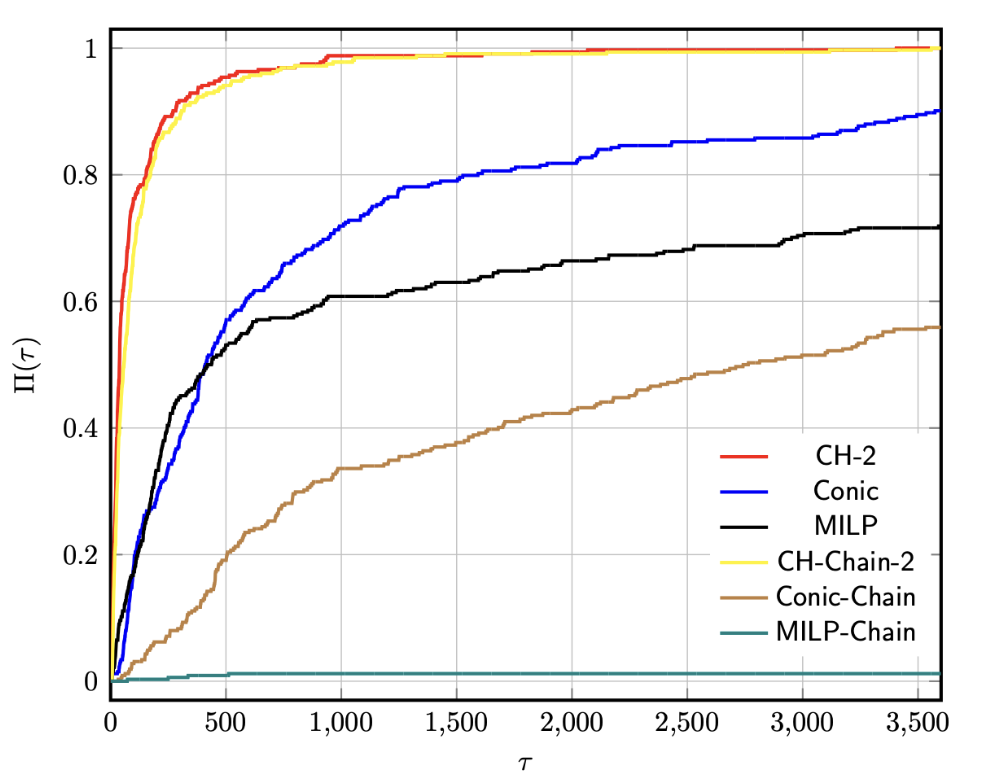}
\caption{Time performance profile in terms of CPU time.}
\label{fig:Perf_Profile_time}
\end{figure}

By examining the \revision{```Max'' and ``Std''  columns} in Table \ref{tab:Pref_NoLuce_Luce}, it becomes evident that our proposed approaches exhibit robustness 
across all the configurations. To further visualize the computational performance, we employ the performance profile, introduced by \citet{dolanBenchmarkingOptimizationSoftware2002}, which is commonly utilized to compare different solution approaches. Figure~\ref{fig:Perf_Profile_time} depicts the performance profile \footnote{Letting $ \tau_{p,c} $ denote the CPU time that approach $c$ consumed in solving problem $ p $ and $\Theta$ denote the set of $n_p$ test problems, the distribution function $\pi_c(\tau)$ for approach $ c $ is defined as 
$\pi_c(\tau) = {size\{p\in {\Theta}: \tau_{p,c} \le \tau\}}/{n_p}, \tau>0$,
which serves as a performance metric for the approach $c$, with respect to the CPU time. The performance profile of each approach $c$ is generated by plotting the corresponding distribution function $\pi_c(\tau)$.} in terms of CPU time on the instances tested in Table \ref{tab:Pref_NoLuce_Luce}.

The performance profile clearly demonstrates that our approach outperforms the other approaches in terms of CPU time. Moreover, the performance difference between \textsf{CH-2} and \textsf{CH-Chain-2} is very small, indicating the scalability of our approaches. 
Notably, both \textsf{CH-2} and \textsf{CH-Chain-2}  successfully solve $95\%$ of all instances within 500 seconds, while other approaches solve at most $60\%$ of all instances within the same time limit.
%

\section{Applications}\label{section:applications}
The goal of this section is to illustrate how the convex hull results developed in Section~\ref{section:convexhull} can be applied to solve various assortment optimization problems. 


\subsection{Constrained assortment planning under MNL}\label{section:MNL}
\revision{In this section, we use the convex hull results in Section~\ref{section:formulation-MNL} to derive compact formulations for constrained assortment optimization problems under MNL defined as follows:}
\begin{equation}\label{eq:cons-MNL}
 \max \biggl\{ \frac{\sum_{j \in N}r_{j}u_jx_{j}}{u_0 + \sum_{j \in N}u_jx_{j}} \biggm| \BFC \x \leq \BFd,\ \x \in \{0,1\}^n \biggr\} \tag{\textsc{Cons-MNL}},
\end{equation}
\revision{where $r_j$ is the revenue generated from selling product $j$, $u_j$ is the preference weight on product $j$, and $\BFC$ in an $k \times  n$ matrix and $\BFb$ is a $k$ dimensional vector. Our convex hull result implies that such formulation can be obtained as long as the convex hull of the feasible region can be. }

   \begin{corollary}\label{cor:consMNL}
   	\revision{Let $\BFA\x + \BFB\boldsymbol{z} \leq \BFb$ be a compact formulation of convex hull of $\bigl\{\x \in \{0,1\}^n \bigm|\BFC \x \leq \BFd\bigr\}$. Then, an LP formulation for~\eqref{eq:cons-MNL} is given as follows:}
   	\[
   	\max \biggl\{ \sum_{j \in N} r_ju_jy_j  \biggm| \BFA \y + \BFB\boldsymbol{z} \leq \BFb y_0,\ y_0 \geq 0,\ u_0 y_0 + \sum_{j \in N}u_jy_j = 1\biggr\}.
   	\] 
   \end{corollary}
\rvtwo{This result allows us to leverage convex hull results in combinatorial optimization~\citep{schrijver2003combinatorial,conforti2014integer}, generating compact LP formulations of~\eqref{eq:cons-MNL}. When the feasible assortment set is defined by a totally unimodular (TU) matrix, its convex hull coincides with its LP relaxation~\citep{conforti2014integer}, recovering the results in~\cite{avadhanula2016tightness} and \cite{sumida2021revenue}.}

\rvtwo{In practice, we often face constraints which are not TU matrix representable. For example,~\cite{housni2024assortment} consider~\eqref{eq:cons-MNL} with general covering constraints, which are common, as retailers often need to meet service-level agreements with suppliers or simply wish to impose diversity in the assortment to meet different customers’ needs. They propose an approximation algorithm for the general case which matches the problem's hardness up to a constant factor. Within the class of problems considered by~\cite{housni2024assortment}, our result helps to identify a subclass of instances that admit a polynomial-time solution method, as described below.
The products are categorized into $K$ categories $\mathcal{C}:=\{C_1, C_2, \ldots, C_K\}$ where $ C_k \subseteq N$ for each $k$, and the offered assortment should include at least $\ell_k \in \Z$ products that belongs to the category $C_k$ for each $k$. In other words, the feasible assortments defined by covering constraints are given as follows:}
\begin{equation}\label{eq:covering}
    \biggl\{\x \in \{0,1\}^n \biggm| \sum_{j \in C_k} x_j \geq \ell_k  \text{ for } k = 1, \ldots, K \biggr\}, \tag{\textsc{Covering}}
\end{equation}
\rvtwo{which is compactly denoted as $\BFC^\mathcal{C}\x \geq \ell$. Assume that $\ell_k = 1$ for each $k$. If the coefficient matrix $\BFC^\mathcal{C}$ belongs to the class of \textit{balanced} matrices, which strictly includes TU matrices,  then the convex hull of~\eqref{eq:covering} coincides with its LP relaxation~\cite[Section 3]{berge1972balanced}. 
More importantly,~\cite{conforti1999decomposition} provide a polynomial time algorithm for recognizing balanced matrices. Consequently,  Corollary~\ref{cor:consMNL}, together with this celebrated recognition algorithm,  enables us to identify and solve instances of~\eqref{eq:cons-MNL} with a class of covering constraints, including, but not limited to, the TU case, in polynomial time. 
} 

\rvtwo{
Besides going beyond the TU case, our approach leads to computationally efficient formulations. To illustrate this, we consider the assortment optimization under the two-stage Luce model~\citep{EcheniqueSaito2019}. It follows readily from Section~\ref{section:2SLM} that this is a special case of~\eqref{eq:cons-MNL} with the following type of constraints describing all undominated products:
\[
\bigl\{\x \in \{0,1\}^n \bigm| 0 \leq x_{a_1} + \cdots + x_{a_k} \leq 1 \text{ for chains } a_1 \prec \cdots \prec a_k \text{ in a partial order } \mathcal{P} \bigr\}.
\]
The convex hull of this set coincides with its continuous relaxation~\citep{stanley1986two}, but has an exponential number of inequalities. Corollary~\ref{cor:consMNL} allows us to leverage the state-of-the-art compact formulation for the chain polytope, proposed by~\cite{fawzi2022lifting} and discussed in Remark~\ref{rmk:2slm-formulation}, to solve the assortment optimization under two-stage Luce model using an LP with $\mathcal{O}(n^2)$ inequalities and $n $ additional variables.
}

\subsection{QAP under the independent demand model}
\label{section:poly}

Oftentimes, due to the low access cost in the online channel as opposed to the offline channel (where a click replaces a physical visit), the rate of non-purchases is much higher in the online environment. Consequently, the independent demand assumption can serve as a reliable choice model. According to \citet{vanRyzin2005}, vast revenue management systems are built upon such a independent demand model  (IDM) in practice. In this section, we aim to capture the distinction in choice models between offline and online channels. We assume that a consumer in the online segment makes purchase decisions in accordance with the independent choice model, while the offline customers continue to adhere to the MNL model.

For each online type $i \in M$ and for each product $j \in N$, we use $\theta_{ij}>0$ to denote the probability that a consumer in online segment $i$ purchases product $j$. Letting $r_{ij}$ be the revenue obtained from selling a product $j \in N$ to a consumer in online segment $i \in M$, if we decide to offer assortment $\x_i \in \{0,1\}^n$ then the expected revenue obtained from the online segment $i$ is $R^{\text{IDM}}_{i}(\x_i) := \sum_{j \in N}r_{ij}\theta_{ij}x_{ij}$.
For offline consumers, we use $u_j$ to denote their preference weight on the product $j$, and $u_0$ to denote their preference weight on the no-purchase option.  Given an offline  assortment $\x_0\in \{0,1\}^n$, the expected revenue is given as 
\[
R_0^\text{MNL}(\x_0) := \frac{\sum_{j \in N}r_{0j}u_jx_{0j}}{u_0 + \sum_{j \in N}u_jx_{0j}}.
\]
Moreover, we assume that the set of feasible assortments that we can offer to the offline store needs to satisfy \textit{precedence constraints}, see~\cite{sumida2021revenue} for applications of precedence constraints in assortment planning. More specifically, letting $G = (N,E)$ be a directed graph, where $N$ is a set of nodes and $E$ is a set of arcs, we assume that the set of feasible assortments is given by $x_{0j}\geq x_{0k}$ for $(j,k) \in E$. Our goal is to solve the following variant of~(\ref{eq:omni-mnl}):
\begin{equation}\label{eq:mnl-ind}
	\begin{aligned}
 \max_{}\quad & \alpha_0R_0^\text{MNL}(\x_0) + \sum_{i \in M} \alpha_{i}R^{\text{IDM}}_{i}(\x_i)\\
\text{s.t.}\quad &\x_i \in \{0,1\}^n  \qquad \text{for } i \in M^+ \\
&\x_0 \geq \x_i  \qquad  \text{for } i \in M \\
& x_{0j} \geq x_{0k} \qquad  \text{for } (j,k) \in E. \\
\end{aligned}\tag{\textsc{QAP-IDM}}
\end{equation}

The main result of this section is to show that inequalities~(\ref{eq:stair-under}) yield a polynomial time solvable linear programming (LP) formulation for~(\ref{eq:mnl-ind}),  that is 
\begin{equation}\label{eq:variant-LP}
\begin{aligned}
   \max_{}\quad & \alpha_0 \sum_{j \in N} r_{0j} u_{j}y_{j}   +  \sum_{i \in M} \alpha_i \biggl( \sum_{j \in N} r_{ij} \theta_{ij}x_{j}  \biggr) \\
   \text{s.t.} \quad & u_0y_0 + \sum_{j \in N}u_jy_j=1 \\
   &0\leq y_j\leq y_0 \qquad \text{ for } j \in N \\
   & y_j \geq y_k \qquad \text{ for } (j,k) \in E \\
   &(x_j,y_0,\y) \text{ satisfies } (\ref{eq:stair-under}) \quad \text{ for } j \in N \text{ and } S \subseteq N \setminus \{j\}. 
\end{aligned}\tag{\textsc{QAP-IDM-LP}}
\end{equation}
In particular, we use a randomized rounding approach to prove the correctness of the LP formulation. This approach has been used in~\cite{teo1998geometry} to provide an elegant proof of the integrality of the stable marriage polytope, and more applications can be found in Section 3.3 of~\cite{bertsimas2005optimization}. The general idea of the rounding approach is as follows. We solve~(\ref{eq:variant-LP}) and obtain an optimal solution $(x, y_0, \y)$ and optimal revenue $r^*$. Then, from $\y$ we create a new random assortment, that is feasible to~(\ref{eq:mnl-ind}), and show that the expected revenue generated from the random assortment is $r^*$. This shows the equivalence between~(\ref{eq:mnl-ind})  and~(\ref{eq:variant-LP}). The key in the proof is to design an appropriate randomization method, which is inspired by the exactness proof of~\eqref{eq:stair-under} in Theorem~\ref{them:convexhull}. 

\begin{theorem}\label{them:polyLP}
\revision{If $r_{ij} \geq 0$ for $i \in M^+$ and $j \in N$ then} $(\ref{eq:variant-LP})$ is a polynomial-time solvable LP formulation of $(\ref{eq:mnl-ind})$.
\end{theorem}
\rvtwo{Note that, in reality, due to the precedence constraints,  it is possible that some of products are unprofitable. For example, a retailer needs to select an unprofitable product so that it can stock other products with high profit margins due to precedence constraints imposed by its suppliers. In this case, our inequalities~\eqref{eq:stair-under} and~\eqref{eq:stair-over} fail to provide an LP formulation.}

Last, we discuss a consequence of Theorem~\ref{them:polyLP}. In~\cite{cao2023revenue}, they consider assortment optimization problems when customers choose under a mixture of MNL and IDM, that is 
  \begin{equation}
\begin{aligned}\label{eq:mixture}
 \max_{} \Biggl\{ \sum_{i \in N}r_i\biggl( \alpha \frac{v_ix_i}{v_0 + \sum_{i \in N} v_ix_i}  +  (1-\alpha) \theta_ix_i \biggr) \Biggm|  \x \in \X \subseteq \{0,1\}^n   \Biggr\},
\end{aligned}\tag{\textsc{Mixture}}
\end{equation}
where $r_i$ is the revenue obtained from selling product $i \in N$, $v_i$ denote consumers' preference on product $i \in N$,  $\theta_i$ is the probability that a consumer select product $i \in N$, and $\alpha$ is the probability that an arriving customer is in the first segment. It is shown in~\cite{cao2023revenue} that~\eqref{eq:mixture} is polynomial time tractable by solving a linear program when $\X = \{0,1\}^n$, and, unfortunately, is NP-hard when the feasible region $\X$ is defined by a cardinality constraint. A natural follow-up question is to characterize a class of constraints under which~\eqref{eq:mixture} is polynomial time solvable. Next, we show that Theorem~\ref{them:polyLP} implies that~\eqref{eq:mixture} with precedence constraints is polynomial time solvable.
\begin{corollary}\label{cor:mixture}
  Let $G = (N,E)$ be a directed graph. If $\X= \bigl\{\x \in \{0,1\}^n \bigm| x_j \geq x_k \text{ for } (j,k) \in E \bigr\}$ \revision{ and $r_{ij} \geq 0$ for $i \in M^+$ and $j \in N$} then inequalities~\eqref{eq:stair-under} yield a polynomial time solvable LP formulation for~\eqref{eq:mixture}.
\end{corollary}

\subsection{Other configurations for QAP}\label{section:frontwarehouse}

\revision{Thus far, we have focused on one particular configuration for QAP, where the online assortment is a subset of the products available in offline stores. However, two alternative configurations merit consideration. The first involves the firm establishing a network of front warehouses in urban areas, enabling rapid delivery of fresh food to nearby consumers. These front warehouses serve only online customers and do not cater to offline shoppers. A prominent example is Dingdong (Cayman), a Chinese company listed on Nasdaq, which offers a similar rapid delivery promise as Alibaba's Hema Fresh. The second configuration allows firms to use the backroom space of offline stores for storing products \textit{exclusively} for online customers. We will demonstrate how to adjust the parameter settings to enable our QAP model to accommodate both configurations.}

\revision{For the first configuration involving front warehouses, we let the index $i= 0$ denote the front warehouse and set $\alpha_0 = 0$ to indicate the absence of offline consumer traffic at the front warehouse. Furthermore, we define $\X_0 = \{\x_0\mid \sum_{j\in N} x_{0j} \leq k\}$ in \eqref{eq:omni-mnl} to represent the limited storage capacity of the front warehouse in urban areas. For the second configuration, the index $i = 0$ represents a dummy consumer type with $\alpha_0 = 0$, while $i = m+1$ denotes the offline consumer type, which does not have access to the exclusive online products stored in the backroom space of offline stores. Then, we redefine $M = \{1, \dots, m, m+1\}$.  The set of available products for the dummy consumer type, denoted as $A = \{j \in N \mid x_{0j} = 1 \}$, represents all products offered in the offline store, including those stored in the backroom exclusively for online customers. In addition, the set of available products for type ${m+1}$, denoted as $B = \{j \in N \mid x_{m+1,j} = 1 \}$, comprises products available to offline consumers. By applying the \eqref{eq:omni-mnl} formulation under this setting, together with $\X_0 = \{\x_0\mid \sum_{j\in N} x_{0j} \leq k\}$, we observe that the optimal solution ensures $B \subseteq A$. This models that the assortment for offline consumers is a subset of the total assortment. Additionally, $A \cap \bar{B}$ represents the products stored in the backroom space exclusively for online customers. In other words, by setting appropriate parameters,  \eqref{eq:omni-mnl} can address the second configuration.}

\revision{With such configuration of parameters and assuming uniform pricing across consumer types (i.e., $r_{ij} = r_{i'j}$ for $i \neq i'$), our QAP formulation reduces to the problem studied in \cite{HousniTopaloglu2022}, which focuses on evaluating the value of personalization.  In Section~\ref{EC:personalization}, we present a numerical study to demonstrate the computational efficacy of our formulation under such setting.}

\subsection{Inventory Decision under QAP}
\label{subsec:inventory_QAP}

\rvtwo{To evaluate the performance of the assortment generated by \eqref{eq:omni-mnl}, it is essential to account for the effect of potential stockouts, as they directly influence realized demand. Therefore, inventory planning plays a critical role in realizing the full value of the assortment solution. In this section, we adopt the rounding procedure proposed by \cite{zhang2024leveraging} to determine implementable inventory levels. We assume that at the beginning of each day, the retailer replenishes inventory to predetermined target levels for all products in the assortment. Performance is measured by total profit over a fixed time horizon, where one consumer arrives in each epoch. Total profit is defined as the revenue from product sales minus the corresponding inventory ordering costs.}

\rvtwo{The inventory targets are derived from the solution to \eqref{eq:omni-mnl}. Specifically, solving \eqref{eq:omni-mnl} yields the optimal consumer-product matching probabilities, from which we compute the purchase probability for each product. Multiplying these probabilities by the total number of consumer arrivals gives the expected demand assuming no stockouts. Since these demand values are typically fractional, we apply a rounding procedure to convert them into integer order quantities. Each product’s expected demand is rounded either up or down, with the total inventory constrained to match the rounded total demand. To decide which products are rounded up, we rank them by a weighted price metric that captures their expected revenue contribution across consumer segments. Products with higher weighted prices are given priority, resulting in the final integer inventory levels. This procedure allows us to compare the fluid approximation revenue with the simulated revenue generated under the integer inventory decisions.}

\rvtwo{Our rounding procedure differs from that of \cite{zhang2024leveraging} in two key aspects: the assortment is personalized for different consumer segments, and product prices may vary across segments. The detailed procedure and numerical results are provided in Online Appendix~\ref{ec:inventory_simulation}. Despite these differences, similar to the findings in \cite{zhang2024leveraging}, our rounding approach demonstrates strong performance, especially as the time horizon increases. }

\rvtwo{However, the current setting assumes a fixed number of consumers $T$ per day. In practice, the number of arriving consumers is random, introducing additional complexity to inventory decisions. We leave the study of such scenarios to future research.}

\section{Conclusion}\label{section:conclusion}
In this paper, we study the assortment problem within the quick-commerce context. Recognizing the face-inducing property of the physical linkage constraint, we adopt a convex hull representation to depict the customer choice behavior for each online segment. Exploiting the algebraic and geometric structures of the MNL choice model, we unveil a series of convex hull results that shed light on the set of choice probabilities under general operational constraints and on the interaction between offline assortment decisions and online/offline choice probabilities. These methodologies effectively address the challenges posed by the quick-commerce assortment problem.

We hope that our results will inspire further investigation of integer programming methods in addressing assortment optimization challenges. To elaborate, we envision three specific avenues for future research. First, it is natural to inquire whether our geometric insights into Luce's choice axioms can be applied to derive compact polyhedral formulations for other Luce-type choice models proposed in~\cite{kovach2022focal,tserenjigmid2021order,EcheniqueSaito2019,echenique2018perception}. Second, it is interesting to investigate the performance of inequalities~\eqref{eq:stair-over} and~\eqref{eq:stair-under} in solving other assortment optimization problems. \revision{Last, it is worth to study how to adapt the Benders decomposition based solution approach proposed in~\cite{BertsimasMisic2019} to solve our formulation, and investigate its computational benefit.}


\bibliographystyle{informs2014} 
\bibliography{ref}
\newpage  

\newpage

\ECSwitch

\ECDisclaimer

\ECHead{Additional discussions, missing proofs, and formulations}

\section{Proofs}
\subsection{Proof of Proposition~\ref{prop:more-offering}}

Before proving this result, we introduce few notations. Let $\x^{\textsc{QAP}}$ be an \textit{maximal} element, with respect to the ordering associated with the nonnegative orthant, of the optimal solution set of~\eqref{eq:omni-mnl}, and let $\x^{\textsc{RO}}$ be an solution satisfying  $(\ref{eq:2SNR})$. Consider the expected revenue from the offline channel of $(\ref{eq:2SNR})$ and $(\ref{eq:omni-mnl})$, that is, $\pi^{\textsc{RO}}_0 = R^{\text{MNL}}_0(\x^{\textsc{RO}}_0)$ and $\pi^{\textsc{QAP}}_0 = R^{\text{MNL}}_0(\x^{\textsc{QAP}}_0)$, respectively. It follows readily that  $\pi^{\textsc{RO}}_0 \geq \pi^{\textsc{QAP}}_0$ since  $\x^{\textsc{RO}}_0$ is obtained by solely maximizing the offline channel's revenue. Next, given that the offline channel assortment is $x^{\textsc{QAP}}_0$, consider the online assortment optimization problem 
\begin{equation}\label{eq:QAP-online}
     \max_{\x_i}\Bigl\{R_i^{\textsc{MNL}} (\x_i) \Bigm|  \x_i \leq \x^{\textsc{QAP}}_0,\ \x_i \in \{0,1\}^n \Bigr\} \qquad \text{ for } i \in M,
\end{equation}
where the optimal value is denoted as $\pi_i^{\textsc{QAP}}$. Since $\x^{\textsc{QAP}}$ is optimal, it follows readily that $\pi_i^{\textsc{QAP}} = R_i^{\textsc{MNL}}(\x^{\textsc{QAP}}_i)$, that is, $\x^{\textsc{QAP}}_i$ is an optimal solution of~\eqref{eq:QAP-online}. In other words, $\x_i^{\textsc{QAP}}$ is obtained via the revenue-ordered policy. Last, for each $i \in M$, let $\pi_i^{\textsc{RO}} = R_i^{\textsc{MNL}}(\x_i^{\textsc{RO}})$.


First, we prove $\x^{\textsc{RO}}_0 \leq \x^{\textsc{QAP}}_0$ by contradiction. Suppose that there exists product $k$ such that $\x^{\textsc{RO}}_{0k}=1$ and  $\x^{\textsc{QAP}}_{0k}=0$. We will show that adding product $k$ into the offline assortment $\x^{\textsc{QAP}}_0$ will not reduce the total revenue. This contradicts that $\x^{\textsc{QAP}}$ is an maximal element of the optimal solution set of $(\ref{eq:omni-mnl})$. Let $S^{\textsc{QAP}}_0$ be the assortment for offline channel corresponding to the solution $\x^{\textsc{QAP}}$, and let $\x^{\textsc{QAP} \cup \{k\}}_0$ be the binary vector corresponding to the $S^{\textsc{QAP}}_0 \cup \{k\}$. Then, 
\begin{eqnarray*}
R^{\text{MNL}}_0\bigl(\x^{\textsc{QAP} \cup \{k\}}_0\bigr) & = & \frac{\sum_{j \in S^{\textsc{QAP}}_0 \cup \{k\}}  u_{0j} r_{0j}}{u_{00}+\sum_{j \in S^{\textsc{QAP}}_0 \cup \{k\}} u_{0j}} \\
& = & \frac{u_{00}+\sum_{j \in S^{\textsc{QAP}}_0} u_{0j}}{u_{00}+\sum_{j \in S^{\textsc{QAP}}_0 \cup \{k\}} u_{0j} } \frac{\sum_{j \in S^{\textsc{QAP}}_0} u_{0j} r_{0j}}{u_{00}+\sum_{j \in S^{\textsc{QAP}}_0} u_{0j}} + \frac{u_{0k}}{u_{00}+\sum_{j \in S^{\textsc{QAP}}_0 \cup \{k\}} u_{0j}} r_{0k}\\
& = & \frac{u_{00}+\sum_{j \in S^{\textsc{QAP}}_0} u_{0j}}{u_{00}+\sum_{j \in S^{\textsc{QAP}}_0 \cup \{k\}} u_{0j}} \pi^{\textsc{QAP}}_0 + \frac{u_{0k}}{u_{00}+\sum_{j \in S^{\textsc{QAP}}_0 \cup \{k\}} u_{0j}} r_{0k} \\
& \geq & \pi^{\textsc{QAP}}_0 ,
\end{eqnarray*} 
where the inequality follows from $r_{0k} \geq \pi^{\textsc{RO}}_0 \geq \pi^{\textsc{QAP}}_0$. Thus, we obtain a feasible solution $(\x^{\textsc{QAP} \cup \{k\}}_0, \x^{\textsc{QAP}}_1, \ldots,  \x^{\textsc{QAP}}_m)$ of $(\ref{eq:omni-mnl})$, whose objective value is at least the optimal value of $(\ref{eq:omni-mnl})$.

Next, we argue that  $\x^{\textsc{RO}}_i \leq \x^{\textsc{QAP}}_i$ for each $i \in M$. Since $\x^{\textsc{RO}}_0 \leq \x^{\textsc{QAP}}_0$, it follows that for each $i \in M$, $\pi_i^{\textsc{QAP}} \geq \pi_i^{\textsc{RO}}$. Moreover, $\x_i^{\textsc{QAP}}$ and $\x_i^{\textsc{RO}}$ can be obtained using revenue-ordered argument on $\x_0^{\textsc{QAP}}$ and $\x_0^{\textsc{RO}}$, respectively. Therefore,  by the same price rank assumption and the maximality of $\x^{\textsc{QAP}}$, we conclude that $x_i^{\textsc{RO}} \leq \x_i^{\textsc{QAP}}$. \hfill \Halmos

\subsection{Proof of Proposition~\ref{prop:complexity}}
The proof of this proposition draws inspiration from the proof technique used in Theorem 1 of \citet{rusmevichientong2014assortment}, where a reduction from the Partition problem is employed. In this proof, we will use $[n]$ to denote $\{1, \ldots, n\}$, and for a given vector $c \in \R^n$ and $S \subseteq [n]$ we will use $c(S)$ to denote $\sum_{j\in S}c_j$. We will prove that a special case of~(\ref{eq:omni-mnl}) is NP-hard, which is defined as follows. Given a set of $n$ products, a revenue vector $r \in \R^n$, two utility vectors $u \in \Z^n_+$ and $v \in \Z^n_+ $, and define
\[
R_0(S)=\frac{\sum_{j \in S} r_{j} u_{j}}{1+\sum_{j\in S}u_j} \quad \text{and} \quad R_1(T) = \frac{\sum_{j \in T}r_{j} v_{j} }{1+\sum_{j\in T}v_j},
\] 
we consider the following maximization problem 
\begin{equation}\label{eq:QAP-Special}
\max\bigl\{R_0(S)+ R_1(T) \bigm| T \subseteq S \subseteq [n] \bigr\}. \tag{\textsc{QAP-Sepcial}}	
\end{equation}
To prove the NP-hardness of~(\ref{eq:QAP-Special}), we map an arbitrary instance of \textsc{Partition} problem to an equivalent~(\ref{eq:QAP-Special}) problem. The \textsc{Partition} problem is defined as follows.  Consider a set of $n$ items and each item $i$ has a size of $c_i \in \Z_+$. Is there a subset $S \subseteq [n]$ such that $c(S) = c([n] \setminus S)$? Let $C = c([n])/2$. Note that $c(S)=c([n] \setminus S)$ if and only if $c(S)=C$. Therefore, we may assume without loss of generality that $C \in \mathbb{Z}_{+}$. 
Given an instance of the \textsc{Partition} problem, $n$ items, a vector $c \in  \mathbb{Z}^n_{+}$ and a $C \in  \mathbb{Z}_{+}$, we construct an instance of~(\ref{eq:QAP-Special}) by assuming that there are $n+1$ products, and letting 
\begin{equation}\label{eq:hardness-parameters}
	r_{j}=\begin{cases}1 & \text { if } j \in [n]\\ 1+\frac{3}{14C} & \text { if } j=n+1\end{cases} 
    \qquad
    u_{j}= \begin{cases}2 c_{j} & \text { if } j \in [n] \\ 6C & \text { if } j=n+1\end{cases}
    \qquad
    v_{j}= \begin{cases} 4 c_{j}  & \text { if } j \in [n] \\4C & \text { if } j=n+1\end{cases}.
\end{equation}

First, we show that the product $n+1$ must be selected in an optimal solution of the instance of~(\ref{eq:QAP-Special}) specified as in~\eqref{eq:hardness-parameters}. Assume that $(S^*, T^*)$ is an optimal assortment but $n+1 \notin S^*$. Then, 
  \[
  \begin{aligned}
  R_0(S^* \cup \{n+1\}) &= \frac{1+u(S^*)+r_{n+1}u_{n+1}}{1+u(S^*) + u_{n+1}} \\
  &= \frac{1 + u(S^*)}{1+u(S^*) + u_{n+1}}\cdot  \frac{u(S^*)}{1 + u(S^*)} + \frac{u_{n+1}}{1+u(S^*) + u_{n+1}}r_{n+1}  \\
  &=	\frac{1 + u(S^*)}{1+u(S^*) + u_{n+1}}\cdot  R_0(S^*) + \frac{u_{n+1}}{1+u(S^*) + u_{n+1}}r_{n+1} \\
  & > R_0(S^*),
  \end{aligned}
  \]
  where the strict inequality holds since for each $j \in [n]$, $r_{n+1} > r_j \geq R_0(S^*)$. 
Therefore, we obtain a feasible solution $\bigl(S^* \cup \{n+1\}, T^*\bigr)$ with an objective value $R_0(S^*\cup\{n+1\}) + R_1(T^*) > R_0(S^*) + R_1(T^*) $, contradicting to the optimality of $(S^*,T^*)$. Similarly, assume that $(S^*,T^*)$ is an optimal solution with $n+1 \in S^*$ but $n+1 \notin T^*$. Then, 
\[
R_1(T^*\cup\{n+1\}) = \frac{1+v(T^*)}{1+v(T^*)+v_{n+1}}R_1(T^*) + \frac{v_{n+1}}{1+v(T^*)+v_{n+1}}r_{n+1} > R_1(T^*). 
\]
Again, we obtain a feasible solution $(S^*,T^*\cup \{n+1\})$ which has a larger objective value than $R_0(S^*) + R_1(T^*)$, a contradiction.

Next, we argue that  $S^* = T^*$ for any optimal solution $(S^*,T^*)$ to the instance of~(\ref{eq:QAP-Special}) specified as in~\eqref{eq:hardness-parameters}. Suppose that $T^* \subset S^*$, that is there exists $k$ such that $k \in S^*$ but $k \notin T^*$. Then, 
\[
\begin{aligned}
R_1(T^* \cup \{k\}) - R_1(T^*) &= \frac{v_k+v(T^* \setminus \{n+1\})+r_{n+1}v_{n+1}}{1+v(T^*)+v_k} - \frac{v(T^* \setminus \{n+1\})+r_{n+1}v_{n+1}}{1+v(T^*)} \\
&=\frac{v_k(1+ v(T^*)) - v_k\cdot (v(T^* \setminus \{n+1\})+r_{n+1}v_{n+1})  }{(1+v(T^*))(1+v(T^*) +v_k)} \\
&=\frac{v_k +v_kv_{n+1} - v_kr_{n+1}v_{n+1} }{(1+v(T^*))(1+v(T^*) +v_k)}\\
& = \frac{2}{14}\frac{v_k}{(1+v(T^*))(1+v(T^*) +v_k)} \\
&>0	.
\end{aligned}
\]
Thus, $(S^*,T^*\cup\{k\})$ is a feasible solution whose objective value is larger than $R_0(S^*) + R_1(T^*)$, a contradiction.

Last, let $K$ denote the target revenue $K=(16C+15/7)/(1+8C)$, we show the \textsc{Partition} problem has a solution if and only if the optimal value of the instance of~(\ref{eq:QAP-Special}) specified as in~\eqref{eq:hardness-parameters} is $K$. Let 
\[
F(z) :=  \frac{6C(1+\frac{3}{14C})+2z }{1+6C+2z} + \frac{4C(1+\frac{3}{14C})+4 z}{1+4C+4z}.
\]
Then, we obtain that
\[
\begin{aligned}
\max \Bigl\{R_0(S)+ R_1(T) \Bigm| T \subseteq S \subseteq [n+1]\Bigr\}  &=\max \Bigl\{F(z) \Bigm| z = c(S),\ S\subseteq [n]\Bigr\} .
\end{aligned}
\]
The derivative of $F$ is given as follows,
\[
F'(z)= \frac{4}{7} \frac{(2C-2z)(2+10C+8z)}{(1+6C+2z)^2 (1+4C+4z)^2},
\]
which is strictly positive over $[0,C)$ and strictly negative over $(C,\infty)$. Therefore, $F$ has a unique maximum at $C$, that is,
\[
F(z) \leq F(C)=\frac{16C+\frac{15}{7}}{1+8C} \quad \text{ for } z \in  [0,C) \cup (C, \infty).
\]
 Hence, we obtain
\[
\max \Bigl\{R_0(S)+ R_1(T) \Bigm| T \subseteq S \subseteq [n+1]\Bigr\} = \max \Bigl\{ F\bigl(c(S)\bigr) \Bigm| S \subseteq [n] \Bigr\} \leq F(C) = K. 
\]
In other words, there exist a pair of assortments $(S,T)$ with $T \subseteq S \subseteq [n+1]$ whose objective value is $K$ if and only if the inequality holds as equality, where the latter is equivalent to that there exists a subset $S$ of $[n]$ such that  $c(S) = C$. \hfill \Halmos

\subsection{Proof of Theorem~\ref{theorem:partial-convexification}}
To simplify the following presentation, we will use $\hat{\r}_i$ to denote the vector $(r_{i0}\cdot u_{i0}, \cdots, r_{in}\cdot u_{in})$. Let $v_{\text{QAP}}$ and $v_{\text{CONV}}$ denote the optimal objective value of~(\ref{eq:omni-mnl}) and~(\ref{eq:omni-rlx}), respectively. Let 
\[
v_{\mibp}:=\max \Biggl\{ \sum_{i \in M^+}\alpha_i \bigl(\hat{\r}_i^\top \y_i\bigr)\Biggm|(\ref{eq:on-mnl}) ,\ (\ref{eq:off-bi})  \text{ and } 	 	 (\ref{eq:link-bi}) \Biggr\}.
\]
Clearly, it follows readily from the discussion before Theorem~\ref{theorem:partial-convexification} that  $v_{\text{QAP}} = v_{\text{MIBP}}$. On the other side, letting $\Pi:=\prod_{i=1}^m\Pi_i(\X_i)$, the value $v_{\text{CONV}}$ is equal to $v'_{\text{CONV}}$, where
\begin{equation*}
v'_{\text{CONV}}:=\max \Biggl\{ \sum_{i \in M^+}\alpha_i \bigl(\hat{\r}_i^\top \y_i\bigr)\Biggm|  (\y_{\cdot 0},\y) \in \conv(\Pi),\ (\ref{eq:off-bi})  \text{ and } 	 	 (\ref{eq:link-bi}) \Biggr\},
\end{equation*}
and $\y_{\cdot 0}:=(\y_{10}, \ldots, \y_{m0})$ and $\y :=(\y_1, \ldots, \y_m)$. This holds since $\prod_{i=1}^m\conv\bigl(\Pi_i(\X_i)\bigr) = \conv\bigl(\prod_{i=1}^m\Pi_i(\X_i)\bigr)=\conv(\Pi)$.

Next, we will prove that $v_{\text{MIBP}} = v'_{\text{CONV}}$. Clearly, $v_{\text{MIBP}} \leq v'_{\text{CONV}}$ since two problems maximize the same objective function but the feasible region of the left-hand side one is contained in that of the right-hand side one.  To show the opposite direction $v_{\text{MIBP}} \geq v'_{\text{CONV}}$, we consider an optimal solution $(\bar{\x}, \bar{y}_{00}, \bar{\y}_{0},\bar{\y}_{\cdot 0},\bar{\y})$ of the right-hand side problem. The next step is to perform a decomposition on the online choice probability vector $\bar{\y}$. Let $J:=\{j \in N \mid \bar{x}_{j} = 0\}$ and define 
\[
F:=\conv(\Pi) \cap \bigl\{(\y_{\cdot 0},\y ) \bigm| y_{ij} = 0  \text{ for } i \in M \text{ and } j \in J \bigr\},
\]
and
\[
\Pi':=\Pi \cap \bigl\{( \y_{\cdot 0},\y ) \bigm| y_{ij} = 0 \text{ for } i \in M \text{ and } j \in J \bigr\}.
\]
Now, Lemma~\ref{lemma:facial decomposition} can be recursively invoked to obtain $\conv(\Pi') = F$ since, for every $i \in M $ and $j \in J$, the inequality $y_{ij} \geq 0$ is valid for $\Pi$.  This implies that the point $(\bar{\y}_{\cdot 0}, \bar{\y})$, which lies on $F$, is expressible as a convex combination of points $(\y^t_{\cdot 0}, \y^t)_{t \in T}$ of $\Pi'$, that is $(\y_{\cdot 0}, \bar{\y}) = \sum_{t \in T}\lambda^t (\y^t_{\cdot 0}, \y^t)$, where $\sum_{t \in T}\lambda^t = 1$ and $\lambda \geq 0$. The proof is complete by observing
\begin{equation*}
	\begin{aligned}
		v'_{\text{CONV}} &= \sum_{i \in M^+}\alpha_i \bigl(\hat{\r}_i^\top \bar{\y}_i\bigr) = \sum_{t \in T} \lambda^t \Biggl( \sum_{i \in M^+} \alpha_i \bigl(\hat{\r}_i^\top \bar{\y}_i^t\bigr) \Biggr)  \leq \sum_{t \in T}\lambda^tv_{\text{MIBP}} = v_{\text{MIBP}},
	\end{aligned}
\end{equation*}
where the inequality holds because $\lambda^t$ is nonnegative and for each $t$ the point $(\bar{\x},\bar{y}_{00},\bar{\y}_0,\y^t_{\cdot 0},\y^t)$ is a feasible solution to the underlying maximization problem of $v_{\text{MIBP}}$, the last equality holds as $\sum_{t \in T} \lambda^t = 1$. \hfill  \Halmos

\subsection{Proof of Lemma~\ref{lemma:fractional-LP}}
Clearly, the convex hull of $\Pi(\X)$ is contained in $H \cap \conv(K)$  since the latter set is a convex relaxation of $\Pi(\X)$. To show the reverse containment, we consider a point $(y_0,\y)$ in $H \cap \conv(K)$. Clearly, $(y_0,\y) \in H$, where $H$ is a hyperplane not passing $0$, implies that $(y_0,\y) \neq 0$. Besides, $(y_0,\y) \in \conv(K) \setminus \{0\}$ implies that there exist a set of points $(y^t_0,\y^t)_{t \in T} \subseteq K \setminus \{0\}$, $\lambda \geq 0 $ and $\sum_{t \in T} \lambda^t=1$ such that
\[
(y_0,\y) = \sum_{t \in T} \lambda^t \cdot (y^t_0,\y^t), 
\]
For $t \in T$, let $\beta^t := u_0y^t_0 + \sum_{j\in N} u_j y^t_j$. Then, $(y_0^t,\y^t) \geq 0$ and $(y_0^t,\y^t) \neq 0$  imply  $\beta^t > 0$ for every $t \in T$. Moreover, we obtain that 
\[
\sum_{t \in T}\lambda^t \beta^t =  \sum_{t \in T} \lambda^t \cdot \biggl(u_0 y^t_0 + \sum_{j\in N} u_j y^t_j\biggr) =   u_0 y_0 + \sum_{j\in N} u_j y_j = 1,
\]
where the last equality holds due to $(y_0,\y) \in H$. It turns out that
\[
(y_0,\y) = \sum_{t \in T}(\lambda^t \beta^t)\cdot \biggl(\frac{1}{\beta^t}(y^t_0,\y^t)\biggr)\quad\text{and} \quad \frac{1}{\beta^t}(y^t_0,\y^t) \in K \cap H \text{ for } t \in T.
\]
In other words, $(y_0,\y)$ in $\conv(K\cap H)$, showing $H \cap \conv(K) \subseteq \conv(K\cap H)$. 

Next, letting $R:=\bigl\{(y_0,\y)\bigm| \y \in y_0 \cdot \conv(\X),\ y_0 > 0 \bigr\}$, we prove that $R = \conv(K)$. Clearly, $\conv(K)\subseteq R$ since $K \subseteq R$ and $R$ is convex. To show the reverse containment, we consider a point $(y_0,\y) \in R$. Then, $(y_1/y_0, \ldots, y_n/y_0)$ belongs to $\conv(\X)$, and thus is expressible as a convex combinations of points $\{\x^t\}_{t \in T}$ in $\X$. Therefore, $(y_0,\y)$ is a convex combination of points $\{(y_0,y_0\cdot\x^t)\}_{t \in T}$, which is a subset of the cone $K$. This shows that $(y_0,\y) \in \conv(K)$. \hfill \Halmos

\subsection{Proof of Theorem~\ref{them:extended-formulation}}
The proof is complete by observing that  
\[
\begin{aligned}
\conv\bigl(\Pi(\X)\bigr) &= \bigl\{(y_0,\y) \bigm| (y_0,\y) \in H,\ \BFA\y + \BFB \boldsymbol{z} \leq by_0,\ y_0 >0 \bigr\} \\
 &= \bigl\{(y_0,\y) \bigm| (y_0,\y) \in H,\ \BFA\y + \BFB \boldsymbol{z} \leq by_0,\ y_0 \geq 0 \bigr\},
\end{aligned}
\]
where the first equality holds due to the extended formulation of $\X$ and Lemma~\ref{lemma:fractional-LP}. To see the second equality, we assume that there exists a point $(\bar{y}_0,\bar{\y})$ of the last set with $\bar{y}_0 = 0$. Since $\conv(\X)$ is a subset of $[0,1]^n$ it follows readily that for each $j$ $y_j\leq y_0$ is implied by the system of inequalities $\BFA\y + \BFB \boldsymbol{z} \leq by_0$. Thus, we obtain $\bar{y}_j \leq \bar{y}_0 = 0$. Therefore, $(\bar{y}_0,\bar{\y}) = 0$. This contradicts with $(\bar{y}_0,\bar{y}) \in H$ as $H$ does not contain $0$. \hfill \Halmos

 \subsection{Proof of Proposition~\ref{prop:hullcomplexity}} 
In this proof, we only argue that the separation problem for $\B^\leq$ is NP-hard, and remark that the hardness result for $\B$ follows from a similar argument by dropping the non-positivity requirement on $c$.  Consider an assortment optimization with product costs 
\begin{equation}\label{eq:mnl-cost}
\max \biggl\{ \frac{\sum_{j \in N}r_ju_jx_j}{u_0+\sum_{j \in N}u_jx_j} + c^\top \x \biggm| \x \in \{0,1\}^n \biggr\},\tag{\textsc{MNL-Cost}}	
\end{equation}
where $u_0>0$,  $u$ and $r$ are positive vectors in $\R^n$, and $c$ is non-positive in $\R^n$.~\cite{kunnumkal2019tractable} has shown that~(\ref{eq:mnl-cost}) is NP-hard. On ther other hand, by introducing variables $y_j= x_j/(u_0 + \sum_{j \in N} u_jx_j)$ and $y_0 = 1/(u_0 + \sum_{j \in N} u_jx_j)$, solving an instance of~(\ref{eq:mnl-cost}) reduces to optimizing a linear function over $\B^\leq$, that is,
\[
\max \biggl\{\sum_{j \in N} r_ju_jy_j + c^\top\x \biggm| (\x,y_0,\y) \in \B^\leq \biggr\},
\]
which is equivalent to optimizing a linear function over $\conv(\B^\leq)$ \citep[see Lemma 1.3 in][]{conforti2014integer}. In other words, optimization a linear function over the convex hull of $\B^\leq$ is NP-hard. Thus, by the equivalence of separation and optimization, \citep[see Theorem 6.49 in][]{grotschel2012geometric}, the separation problem for $\B^\leq$ is NP-hard. \hfill \Halmos

\subsection{Proof of Theorem~\ref{them:staircase}}\label{appendix:staircase}
\begin{lemma}\label{lemma:fractional-LP-2}
\revision{For any set $\X \subseteq \{0,1\}^n$, the vertex set of $\conv(\Pi(\X))$ is  $H \cap \bigl\{(y_0,\y) \bigm| \y \in y_0 \cdot \X,\ y_0 > 0\bigr\}$.}
\end{lemma}
\noindent \textbf{Proof of Lemma~\ref{lemma:fractional-LP-2}: }  \revision{Due to $\vertex\bigl(\conv(\Pi(\X))\bigr) = \X$, it suffices to show that  the vertices of $\conv(\Pi(\X))$ is given by}
\[
 V := H \cap \Bigl\{(y_0,\y) \Bigm| \y \in y_0 \cdot \vertex\bigl(\conv(\X)\bigr),\ y_0 > 0\Bigr\}.
\]
The vertices of the convex hull of $\Pi(\X)$ are contained in $V$ by observing that 
\[
 \conv\bigl(\Pi(\X)\bigr) =  H \cap \Bigl\{(y_0,\y) \Bigm| \y \in y_0 \cdot \conv(\X),\ y_0 > 0\Bigr\} = \conv(V) ,
\]
\revision{where the first equality holds Lemma~\ref{lemma:fractional-LP}, and the second equality follows from Lemma~\ref{lemma:fractional-LP} and $\conv(\X) = \conv\bigl(\vertex\bigl(\conv(\X)\bigr)\bigr)$.  To prove the reverse containment, we consider a point $(\bar{y}_0,\bar{\y})$ of $V$, and assume that it is not a vertex of the convex hull of $\Pi(\X)$. Then, it can be expressed as a convex combination of distinct points in $\Pi(\X)$. Therefore, $(\bar{y}_1/\bar{y}_0, \ldots, \bar{y}_n/\bar{y}_0)$ can be expressed as a convex combination of distinct points in $\X$, contradicting with the fact that $(\bar{y}_1/\bar{y}_0, \ldots, \bar{y}_n/\bar{y}_0)$ is a vertex of $\conv(\X)$.} \hfill  \Halmos

\noindent \textbf{Proof of Theorem~\ref{them:staircase}: } We start with introducing the vertices of $Y$. It follows readily from Lemma~\ref{lemma:fractional-LP-2} that vertices of $Y$ can be characterized using subsets of $N$. Namely, for each subset $S \subseteq N$, we denote as  $\ext(S)$ the vertex of $Y$ associated with $S$, where each coordinate is given as follows: 
\[
\ext_i(S) = \begin{cases}
	\alpha(S) & \text{for } i \in 0 \cup S \\
	0 & \text{otherwise} 
\end{cases}.
\]

Next, we establish the validity of inequality~(\ref{eq:stair-under}). Let $j \in N$ and let $S$ be a subset of $N \setminus j$. For every point $(x_j,y_0,\y) \in  \B_j$, 
\[
\begin{aligned}
y_j &= x_j \cdot y_0 \\
&=   x_j \cdot \max \bigl\{ y_0, \alpha\bigl(S \cup j\bigr) \bigr\} + x_j \cdot \bigl(y_0 -  \max\bigl\{ y_0, \alpha\bigl(S \cup j \bigr) \bigr\} \bigr) \\
& \geq x_j\cdot  \alpha\bigl(S \cup j\bigr)  +  y_0 -    \max\bigl\{ y_0, \alpha\bigl(S \cup j \bigr) \bigr\} \\
\end{aligned}
\]
where the second equality holds since the right hand side telescopes to the left hand side. The relaxation step is done by replacing $ \max \bigl\{ y_0, \alpha\bigl(S \cup j\bigr) \bigr\}$ with its lower bound $\alpha\bigl(S \cup j\bigr)$, and $x_j$ with its upper bound $1$, whose validity is due to the non-negativity of $x_j$ and the non-positivity of $\max\{ y_0, \alpha(N) \} - \max\{ y_0, \alpha(S \cup j )\}$. Now, we interpret the difference term as a function on the polytope $Y$, that is,  
\[
D(y_0,\y \mid S):=y_0 - \max\bigl\{ y_0, \alpha(S \cup j ) \bigr\} \quad \text{ for } (y_0,\y ) \in Y.
\]
The function $D(\cdot \mid S)$ is concave, and is further relaxed to a linear function by affinely interpolating $D(\cdot\mid S)$ over the following vertices of $Y$ 
\[
\ext\bigl(S \cup j \bigr) \qquad \text{and} \qquad \ext\bigl(S \cup j \cup t \bigr) \quad \text{ for } t \notin S\cup j.
\]
In particular, we obtain a linear function 
\[
L(y_0,\y \mid S ):= - \sum_{t \notin S \cup j} \bigl(u_{t}\cdot	\alpha (S \cup j ) \bigr) \cdot y_t \quad \text{ for } (y_0,\y ) \in Y.
\]
To show $L(\cdot \mid S )$ is is a valid underestimator of $D(\cdot \mid S)$ over the polytope $Y$, it suffices to verify that $D(y_0,\y \mid S) \geq L(y_0,\y \mid S)$ for every vertex $(y_0,\y)$ of $Y$ due to the concavity of $D(\cdot \mid S)$ and the linearity of $L(\cdot \mid S)$. To complete this, we discuss three cases.
\begin{itemize}
	\item For $T \subseteq S \cup j$, $D\bigl(\ext(T) \bigm| S \bigr) =0= L\bigl(\ext(T) \bigm| S \bigr)$. 
	\item For $T$ with $S \cup j \subseteq T$, 
 \[
 D\bigl(\ext(T) \bigm| S \bigr) = \alpha(T) - \alpha(S\cup j) = \biggl(-\sum_{t \in T \setminus (S \cup j) } u_t \biggr) \alpha(S\cup j)  \alpha(T) =  L\bigl(\ext(T) \bigm| S\bigr). 
 \]
\item For $T$ containing elements in both $S\cup j$ and its complement, 
\[
D\bigl(\ext(T) \bigm| S \bigr) \geq D\bigl(\ext( S\cup j \cup T) \bigm| S \bigr) = L\bigl(\ext( S\cup j \cup T) \bigm| S \bigr) \geq L\bigl(\ext(T) \bigm| S \bigr),
\]
where the first inequality holds since $D\bigl(\ext(A)\bigm| S \bigr) \geq D\bigl(\ext(B)\bigm| S \bigr)$ for $A\subseteq B \subseteq N$, the equality follows from the second case, and the last inequality holds by the non-negativity of $u_t$ and of $\ext_t(T) - \ext_t(S\cup j \cup T)$ for every $t \notin S \cup j$.
\end{itemize}

Now, we show the validity of inequality~(\ref{eq:stair-over}) by using a similar two step procedure. For a given $j \in N$ and a given subset $S$ of $N \setminus j$. For every point $(x_j,y_0,\y) \in  \B_j$,
\[
\begin{aligned}
y_j = x_j \cdot y_0 &=  x_j \cdot \min\bigl\{ y_0, \alpha( S \cup j ) \bigr\} + x_j \cdot y_0 -   x_j \cdot \min \bigl\{ y_0, \alpha(S \cup j ) \bigr\} \\
& \leq \alpha(S \cup j ) \cdot x_j + y_0 - \min\bigl\{ y_0, \alpha(S \cup j ) \bigr\}.
\end{aligned}
\]
Here, the difference term $y_0 - \min\bigl\{ y_0, \alpha(S \cup j ) \bigr\}$ is convex, and the second relaxation step  affinely interpolates it over the following vertices of $Y$
\[
\ext(S \cup j ) \qquad \text{and} \qquad \ext\bigl(S \cup j \setminus t\bigr) \quad \text{ for } t \in S.
\]
This yields a linear function $\bigl(1-(u_0+u_j)\cdot\alpha(S\cup j)\bigr) \cdot y_0 + \sum_{t \in S } \bigl(-u_{t}\cdot	 \alpha(S \cup j)\bigr) \cdot y_t$, completing the proof of validity. \hfill   \Halmos


\subsection{Proof of Corollary~\ref{cor:mccormick}}
    For $S = N \setminus j$,~(\ref{eq:stair-under}) yields $y_j \geq  \alpha(N) \cdot x_j$. For $S = \emptyset$, we obtain
    \[
    \begin{aligned}
        y_j &\geq \alpha(j) \cdot x_j  - \sum_{t\in N\setminus j} \alpha(j) u_t\cdot y_t \\
        & = \alpha(j) \cdot x_j +y_0 - \frac{U(j)}{U(j)}y_0 - \sum_{t\in N \setminus j} \frac{u_t}{U(j)}\cdot y_t \\ 
         & = \alpha(j) \cdot x_j + y_0 - \frac{u_jy_0-u_jy_j +1}{U(j)},
    \end{aligned}
    \]
    where the inequality follows from~(\ref{eq:stair-under}). Therefore, $(u_0+u_j)y_j \geq x_j + (u_0+u_j)y_0- (u_jy_0-u_jy_j+1)$. This yields $y_j \geq \alpha(\emptyset)x_j+y_0 - \alpha(\emptyset)$.
    
  For $S = \emptyset$,~(\ref{eq:stair-over}) yields $y_j \leq  \alpha(j) \cdot x_j$.   For $S = N \setminus j$, we obtain 
  \[
    \begin{aligned}
        y_j &\leq  \alpha(N) \cdot x_j + y_0  -  \alpha(N) (u_0+u_j)\cdot y_0  - \alpha(N) \sum_{t \in N\setminus j}u_t y_t\\
        & =  \alpha(N) \cdot x_j + y_0 - \alpha(N)u_jy_0 - \alpha(N) (1-u_jy_j) 
    \end{aligned}
    \]
    where the first inequality follows from~(\ref{eq:stair-over}).  This shows that $U(N)y_j \leq x_j +U(N)y_0 - u_jy_0-1+u_jy_j$, thus yielding $U(N\setminus j)y_j\leq x_j +U(N\setminus j)y_0-1$.

    Last, we argue that inequalities~\eqref{eq:mccormick} are implied. Clearly, $y_j\leq \alpha(j)x_j \leq \alpha(\emptyset)x_j$ as $ x_j \geq 0$ and $\alpha(\emptyset) \geq \alpha(j)$.   In addition, 
    \[
    \begin{aligned}
    y_j &\leq \alpha(N) \cdot x_j + y_0 - \alpha(N)u_jy_0 - \alpha(N) (1-u_jy_j) \\    & \leq \alpha(N) \cdot x_j + y_0 - \alpha(N),	
    \end{aligned}
    \]
where the second inequality holds since $\alpha(N)u_j(y_j - y_0) \leq 0$.    
\hfill \Halmos



\subsection{Proof of Theorem~\ref{them:poly-sep-stair}}
We start with proving the correctness of Algorithm~\ref{alg:sep} on solving the separation problem for~$(\ref{eq:stair-under})$. Let $\mathcal{R}^*$ be the optimal value of~(\ref{eq:sep-cover}). Clearly, $\bigl(u_0+U(N) -  \sum_{t \notin S \cup\{j\}} u_t \bigr) \mathcal{R}^* = \bigl(u_0+u_j + U(S)\bigr)\mathcal{R}^* \geq \bar{x}_j - \sum_{t \notin S \cup j} u_t \cdot \bar{y}_t$ for all $S \subseteq N \setminus \{j\}$ with equality holding at optimal $S$. Therefore, $\bigl(u_0+U(N)\bigr)\cdot\mathcal{R}^* \geq \bar{x}_j  +  \sum_{t \notin S\cup j}(\mathcal{R}^* - \bar{y}_t) \cdot u_t $  for all $S \subseteq N \setminus \{j\}$. It follows readily that an optimal solution can be recovered by solving 
\[
\max \Biggl\{\sum_{s \notin S\cup j}(\mathcal{R}^* - \bar{y}_s)u_s \Biggm|  S \subseteq N \setminus \{j\} \Biggr\}.
\]
It is optimal to offer each product whose revenue $\bar{y}_s$ exceeds $\mathcal{R}^*$. Therefore, an optimal solution is $\bigl\{s \in N \setminus \{j\} \bigm| \bar{y}_s  \geq \mathcal{R}^* \bigr\}$. Let $\sigma$ be a permutation of $N \setminus \{j\}$ such that $\bar{y}_{\sigma(1)} \geq \bar{y}_{\sigma(2)} \geq \cdots \geq \bar{y}_{\sigma(n-1)}  $. It follows that an optimal solution is one of the nested sets 
\[
E^\sigma_0:=\emptyset \quad E^\sigma_1:=\bigl\{\sigma(1)\bigr\} \quad E^\sigma_2:=\bigl\{\sigma(1),\sigma(2)\bigr\}\quad  \ldots \text{ and } E^\sigma_{n-1}:= \bigl\{\sigma(1), \sigma(2), \ldots, \sigma(n-1)\bigr\}.
\]
Thus, we can conclude that $\bar{y}_j \geq \mathcal{R}^* $ if and only if $\bar{y}_j \geq \alpha(S\cup j)\cdot \bar{x}_j-  \sum_{t \notin S \cup j} u_{t} \cdot	\alpha(S \cup j)\cdot \bar{y}_t$ for every $S \in \{E^\sigma_k\mid k=0,\ldots,n-1\}$.  The latter condition is equivalent to 
\[
\bar{x}_j\leq \min \biggl\{ U(S\cup j)\cdot\bar{y}_j + \sum_{t \notin S\cup j}u_t\bar{y}_t \biggm| S \in \{E^\sigma_0, \ldots, E^\sigma_{n-1}\} \biggr\}.
\]
Now, this minimization problem achieve its optimal at $\bigl\{\sigma(1), \ldots,\sigma(k) \bigm| \bar{y}_{\sigma(k)} \geq \bar{y}_j \bigr\}$, which coincides with $S^*$ used in Algorithm~\ref{alg:sep}, completing the proof.

Next, we prove the correctness of Algorithm~\ref{alg:sep} on solving the separation problem for~$(\ref{eq:stair-over})$. Let $\mathcal{C}^*$ be the optimal value of~(\ref{eq:sep-pack}). 	Clearly, $(u_0+ u_j + U(S)) \cdot \mathcal{C}^* \leq \bar{x}_j +\bar{y}_0 \cdot (u_0+ u_j + U(S))  - (u_0+u_j)\cdot \bar{y}_0 - \sum_{t \in S}u_t \cdot \bar{y}_t $ for all $S \subseteq N\setminus \{j\}$. Therefore, $(u_0+u_j)\mathcal{C}^* \leq \bar{x}_j + \sum_{t \in S}  (\bar{y}_0 - \mathcal{C}^*-\bar{y}_t) \cdot u_t$ for all $S \subseteq N \setminus \{j\}$. It follows readily that an optimal solution can be recovered by solving
\[
\min \Biggl\{\sum_{s \in S} (\bar{y}_0  -\bar{y}_s - \mathcal{C}^*) \cdot u_s \Biggm|  S \subseteq N \setminus \{j\} \Biggr\}.
\]
It is optimal to offer each product whose revenue $\bar{y}_s$ exceeds $\bar{y}_0-\mathcal{C}^*$. Therefore, an optimal solution is $\bigl\{s \in N \setminus \{j\} \bigm| \bar{y}_s  \geq \bar{y}_0 -\mathcal{C}^* \bigr\}$. It follows that it is one of the nested sets $E^\sigma_0, \ldots, E^\sigma_{n-1}$. Thus, we can conclude that $\bar{y}_j\leq \mathcal{C}^*$ if and only if $\bar{y}_j \leq \alpha(S \cup j ) \cdot \bar{x}_j  + \bigl(1-(u_0+u_i)\cdot\alpha(S\cup j)\bigr) \cdot \bar{y}_0 + \sum_{t \in S } \bigl(-u_{t}\cdot	 \alpha(S \cup j)\bigr) \cdot \bar{y}_t$ for every $S \in \{E^\sigma_0, \ldots, E^\sigma_{n-1}\}$. The latter condition is equivalent to 
\[
\bar{x}_j\geq \max\biggl\{ u_j\bar{y}_j + \sum_{t\in S}\bigl(\bar{y}_t - (\bar{y}_0-\bar{y}_j)\bigr)u_t\biggm| S \in \{E^\sigma_0, \ldots, E^\sigma_{n-1}\} \biggr\}.
\]
Now, this maximization problem achieves its optimal at $\bigl\{\sigma(1), \ldots,\sigma(k) \bigm| \bar{y}_{\sigma(k)} \geq \bar{y}_0-\bar{y}_j\bigl\}$, which coincides with the set $T^*$ used in Algorithm~\ref{alg:sep}, completing the proof.

Last, the algorithm involves sorting $\bar{y}_t$ in descending order once.  Sorting has a complexity of $\mathcal{O}(n\log n)$, so the overall time complexity of Algorithm~\ref{alg:sep} is $\mathcal{O}(n\log n)$.
\Halmos

\subsection{Proof of Theorem~\ref{them:convexhull}}
\revision{We start with a decomposition lemma.}
\begin{lemma}\label{lemma:decom}
	\revision{$\conv(\B_j) = \conv(\B_j^\leq) \cap \conv(\B_j^\geq)$.}
\end{lemma}
\noindent \textbf{Proof of Lemma~\ref{lemma:decom}:} \revision{Let $L$ (resp. $R$) denote the set in the left hand side (resp. the right hand side). Clearly, $L \subseteq R$. To show $R\subseteq L$, we consider a point $(\bar{x}_j, \bar{y}_0, \bar{\y})$ in $R$. Then, there exists $w \in \B^\leq_j \cap\{ x_j =0\}$ and $r \in \B^\leq_j \cap\{ x_j =1\}$ such that  $(\bar{x}_j, \bar{y}_0, \bar{\y}) = (1-\bar{x}_j) \cdot w + \bar{x}_j \cdot r$.  Moreover, it can be verified that $w$ and $r$ are in $\B_j^\geq$, showing that $w \in \B_j \cap \{x_j=0\}$ and $r \in\B_j \cap \{x_j=1\}$. This establishes that $(\bar{x}_j, \bar{y}_0, \bar{\y})$ is a convex combination of points in $\B_j$. In other words, it belongs to $L$, completing the proof.} \hfill \Halmos 

\noindent \textbf{Proof of Theorem~\ref{them:convexhull}: } \revision{Due to Lemma~\ref{lemma:decom}, it suffices to characterize $\conv(\B_j^\leq)$ and $\conv(\B_j^\geq)$ individually}. First, we show that $\conv(\B_j^\geq) = R^\geq_j$, where $R^\geq_j:=\bigl\{(x_j,y_0,\y) \in [0,1] \times Y \bigm|(\ref{eq:stair-under})  \bigr\}$. It follows readily that $\conv(\B_j^\geq) \subseteq R^\geq_j$ since  by Theorem~\ref{them:staircase} the latter set is a convex relaxation of $\B_j^\geq$. To prove the reverse containment $R^\geq_j \subseteq \conv(\B_j^\geq)$, we consider an extreme point $( \bar{x}_j,\bar{y}_0, \bar{\y})$ of the polytope $R^\geq_j$ and will express it as a convex combination of points in $\B_j^\geq$. It follows readily that $( \bar{x}_j,\bar{y}_0, \bar{\y})$ satisfies one of the inequalities in~(\ref{eq:stair-under}) with equality due to its extremality.  Let $\sigma$ be a permutation of $N$ such that $\bar{y}_{\sigma(1)} \geq \bar{y}_{\sigma(2)} \geq \cdots \geq \bar{y}_{\sigma(n)}$, and let $\sigma^{-1}$ be its inverse. For $k \in N$, define $S_k := \{\sigma(1), \ldots, \sigma(k)\}$ and $S_0 := \emptyset$. By Theorem~\ref{them:poly-sep-stair}, the inequality (\ref{eq:stair-under}) corresponding with $S = S_{\sigma^{-1}(j)-1}$ is tight at $(\bar{x}_j,\bar{y}_0,\bar{\y})$. Next, we use the nested sets to construct $n+1$ points in  $\B^\geq_j$ given as follows,
\[
\bigl(0,\ext(S_0)\bigr), \cdots, \bigl(0,\ext(S_{\sigma^{-1}(j)-1})\bigr), \bigl(1,\ext(S_{\sigma^{-1}(j)})\bigr), \ldots, \bigl(1,\ext(S_{n})\bigr),
\]
and a multiplier $(\lambda_k)_{k=0}^n$ given as follows,
\begin{equation*}
\lambda_0 = \frac{\bar{y}_0-\bar{y}_{\sigma(1)}}{\alpha(\emptyset)} \qquad  \lambda_k = \frac{\bar{y}_{\sigma(k)} -\bar{y}_{\sigma(k+1)} }{\alpha(S_k)} \text{ for } k \in N \setminus \{n\} \quad \text{and} \quad \lambda_n = \frac{\bar{y}_{\sigma(n)}}{\alpha(N)} .
\end{equation*}
This $\lambda$ is a convex multiplier since $\lambda \geq 0$ and 
\begin{equation}\label{eq:convexhull-lambda}
\sum_{k = 0}^n\lambda_k = \frac{\bar{y}_0}{\alpha(\emptyset)} + \sum_{k=1}^n\bar{y}_{\sigma(k)}\biggl(\frac{1}{\alpha(S_k)} - \frac{1}{\alpha(S_{k-1})}\biggr) =u_0\bar{y}_0 + \sum_{k = 1}^nu_{\sigma(k)}\bar{y}_{\sigma(k)} = 1,
\end{equation}
where the last equality holds since $y \in Y$. Therefore, we can conclude that $(\bar{x}_j,\bar{y}_0,\bar{\y}) \in \conv(\B_j^\geq)$ since
\[
(\bar{x}_j,\bar{y}_0,\bar{\y}) = \sum_{k\leq \sigma^{-1}(j)-1}\lambda_k \cdot \bigl(0,\ext(S_k)\bigr) + \sum_{k \geq \sigma^{-1}(j)}\lambda_k \cdot \bigl(1,\ext(S_k)\bigr).
\]
To see the convex combination decomposition, we observe that
\begin{equation}\label{eq:convhull-cc1}
\begin{aligned}
\sum_{k=0}^n\lambda_k\ext(S_k) &= (\bar{y}_0 - \bar{y}_{\sigma(1)}) \frac{\ext(\emptyset)}{\alpha(\emptyset)} + \sum_{k=1}^{n-1}(\bar{y}_{\sigma(k)}-\bar{y}_{\sigma(k+1)})\frac{\ext(S_k)}{\alpha(S_k)} + \bar{y}_{\sigma(n)} \frac{\ext(N)}{\alpha(N)}\\
&= (\bar{y}_0 - \bar{y}_{\sigma(1)})\cdot \bigl(1, \chi(\emptyset) \bigr)+ \sum_{k=1}^{n-1}(\bar{y}_{\sigma(k)}-\bar{y}_{\sigma(k+1)})\cdot \bigl(1, \chi(S_k)\bigr) + \bar{y}_{\sigma(n)}\cdot \bigl(1, \chi(S_n)\bigr)\\
&=(\bar{y}_0, \bar{\y}),
\end{aligned}
\end{equation}
where for $S \subseteq N$, $\chi(S)$ is the characteristic vector of $S$, that is $\chi_k(S) = 1$ if $k \in S$ and $\chi_k(S) = 0$ otherwise. Moreover, we obtain
\begin{equation}\label{eq:convexhull-cc2}
\begin{aligned}
\sum_{k\leq \sigma^{-1}(j)-1}\lambda_k\cdot 0 + \sum_{k\geq\sigma^{-1}(j)} \lambda_k\cdot 1 &= \frac{1}{\alpha(S_{\sigma^{-1}(j)})}\bar{y}_{j} + \sum_{k \geq \sigma^{-1}(j)+1}\biggl(\frac{1}{\alpha(S_{k})} - \frac{1}{\alpha(S_{k-1})}\biggr)\bar{y}_{\sigma(k)} \\
& =  \frac{1}{\alpha(S_{\sigma^{-1}(j)})}\bar{y}_{j}  + \sum_{k \geq \sigma^{-1}(j)+1} u_{\sigma(k)}\bar{y}_{\sigma(k)} \\
& =  \frac{1}{\alpha(S_{\sigma^{-1}(j)})}\bar{y}_{j}  + \sum_{k \notin  S_{\sigma^{-1}(j)} } u_{k}\bar{y}_{k} \\
& = \bar{x}_j,
\end{aligned}
\end{equation}
where the last equality holds since, by Theorem~\ref{them:poly-sep-stair}, the inequality~(\ref{eq:stair-under}) with $S = S_{\sigma^{-1}(j)-1}$ is tight at $( \bar{x}_j,\bar{y}_0, \bar{\y})$.   

Now, we show that $\conv(\B_j^\leq) = R^\leq_j$, where $R^\leq_j:=\bigl\{(x_j,y_0,\y) \in [0,1] \times Y \bigm|(\ref{eq:stair-over})  \bigr\}$. Clearly, $\conv(\B_j^\leq) \subseteq R^\leq_j$. 	To prove the reverse containment, we consider a point $( \bar{x}_j,\bar{y}_0, \bar{\y})$ of $R^\leq_j$ and express it as a convex combination of points in $\B_j^\leq$.  Clearly, $( \bar{x}_j,\bar{y}_0, \bar{\y})$ satisfies one of the inequalities~(\ref{eq:stair-over}) with equality. Let $\tilde{\y}$ be a vector such that $\tilde{y}_t = \bar{y}_t$ for $t \neq j$ and $\tilde{y}_t = \bar{y}_0-\bar{y}_t$ for $t =j$. Let $\sigma$ be a permutation of $N$ such that $\tilde{y}_{\sigma(1)} \geq \tilde{y}_{\sigma(2)} \geq \cdots \geq \tilde{y}_{\sigma(n)}$, and let $\sigma^{-1}$ be its inverse. For $k \in N$, define $S_k := \{\sigma(1), \ldots, \sigma(k)\}$, and $S_0:=\emptyset$. Now, we consider the following $n+1$ points in $\B_j^\leq$
\[
\bigl(1,\ext(S_0 \cup j )\bigr), \ldots, \bigl(1,\ext(S_{ \sigma^{-1}(j) -1} \cup j )\bigr), \bigl(0,\ext(S_{\sigma^{-1}(j)} \setminus j )\bigr),\ldots, \bigl(0,\ext(S_{n} \setminus j )\bigr),
\]
and a vector $(\lambda)_{k=0}^n$ defined as follows: 
\[
\begin{aligned}
	&\lambda_k= \frac{\tilde{y}_{\sigma(k)}-\tilde{y}_{\sigma(k+1)}}{\alpha(S_k \cup j  )} \quad \text{for } k \in \{ 0,\ldots, \sigma^{-1}(j)-1 \}\qquad \text{and} \quad \lambda_k= \frac{\tilde{y}_{\sigma(k)}-\tilde{y}_{\sigma(k+1)}}{\alpha(S_k \setminus j )} \quad \text{for } k \in \{\sigma^{-1}(j),\ldots, n\},
	\end{aligned}
	\]
where $\tilde{y}_{\sigma(0)}$ denote $\bar{y}_0$ and $\tilde{y}_{\sigma(n+1)}$ denote $0$. Here, $\lambda$ is a convex multiplier since $\lambda \geq 0$ and 
	\[
	\sum_{k = 0}^n \lambda_k   = (u_0+u_j)\bar{y}_0-u_j(\bar{y}_0-\bar{y}_j) +\sum_{k \in N\setminus j}u_k\bar{y}_k =1.
	\]
We can conclude that $(\bar{x}_j,\bar{y}_0,\bar{\y}) \in \conv(\B_j^\leq)$ since $( \bar{x}_j,\bar{y}_0, \bar{\y})$ can be expressed as a convex combination of points in $\B^\leq_j$, that is,
	\[
	(\bar{x}_j,\bar{y}_0,\bar{\y}) = \sum_{k \leq \sigma^{-1}(j)-1} \lambda_k \bigl(1,\ext(S_k \cup j )\bigr) + \sum_{k \geq \sigma^{-1}(j)} \lambda_k \bigl(0,\ext(S_k \setminus j)\bigr).
	\]
To see this, we observe 
\[
\begin{aligned}
\sum_{k=0}^{\sigma^{-1}(j)-1}\lambda_k\ext(S_k\cup j) &+ \sum_{k=\sigma^{-1}(j)}	^n\lambda_k\ext(S_k\setminus j) =  \sum_{k=0}^{\sigma^{-1}(j)-1}(\tilde{y}_{\sigma(k)}-\tilde{y}_{\sigma(k+1)})\frac{\ext(S_k\cup j)}{\alpha(S_k \cup j)}  \\
&+ \sum_{k=\sigma^{-1}(j)}^{n}(\tilde{y}_{\sigma(k)}-\tilde{y}_{\sigma(k+1)})\frac{\ext(S_k\setminus j)}{\alpha(S_k \setminus j)}  \\
& = \sum_{k=0}^{\sigma^{-1}(j)-1}(\tilde{y}_{\sigma(k)}-\tilde{y}_{\sigma(k+1)})\bigl(1,\chi(S_k\cup j) \bigr) + \sum_{k=\sigma^{-1}(j)}^{n}(\tilde{y}_{\sigma(k)}-\tilde{y}_{\sigma(k+1)})\bigl(1,\chi(S_k\setminus j) \bigr)\\
& = (\bar{y}_0,\bar{\y}),
\end{aligned}
\]
and
\[
\sum_{k=0}^{\sigma^{-1}(j)-1}1 \cdot \lambda_k + \sum_{k=\sigma^{-1}(j)}^{n} 0 \cdot \lambda_k= (u_0+u_j)\cdot\bar{y}_0+\Bigl(u_0+u_j+ \frac{1}{\alpha(S_{\sigma^{-1}(j)-1})} \Bigr) \cdot(\bar{y}_j-\bar{y}_0) + \sum_{t \in S_{\sigma^{-1}(j)-1}}u_t\bar{y}_t = \bar{x}_j,
\]
where the last equality holds since by Theorem~\ref{them:poly-sep-stair}, the inequality~(\ref{eq:stair-over}) with $S = S_{\sigma^{-1}(j)-1}$ is tight at $( \bar{x}_j,\bar{y}_0, \bar{\y})$.  Therefore, $( \bar{x}_j,\bar{y}_0, \bar{\y}) \in \conv(R_j^\leq)$.
\hfill \Halmos

\subsection{Proof of Corollary~\ref{cor:cutting-plane}} 
\revision{The second statement follows from~\cite{kelley1960cutting}. Next, we prove the first statement.  Due to Theorem~\ref{theorem:partial-convexification},~\eqref{eq:omni-mnl} is equivalent to~\eqref{eq:omni-rlx}. Now, we argue that \textsf{CH-0} is an exact reformulation of~\eqref{eq:omni-mnl}.  Clearly, it follows from Theorem~\ref{them:extended-formulation} that the constraint~\eqref{eq:CH-1} describes the convex hull of $\Pi_i(\X_i)$. Moreover, since $\x$ is binary,~\eqref{eq:mccormick} exactly models the relation $y_{ij} = x_jy_{i0}$. This, together with Corollary~\ref{cor:mccormick}, implies that the last two constraints in \textsf{CH-0} exactly models~\eqref{eq:off-bi} and~\eqref{eq:link-bi}. \hfill \Halmos }

\subsection{Proof of Corollary~\ref{cor:consMNL}} 
\revision{By C-C transformation, the constrained assortment optimization~\eqref{eq:cons-MNL} is equal to}
\[
\max \biggl\{ \sum_{j \in N}r_ju_jy_j \biggm| u_0 y_0 + \sum_{j \in N}u_jx_j = 1,\  \BFC \BFy \leq \BFd y_0 ,\ y_j \in \{0,y_0\} \text{ for } j \in N  \biggr\}.
\]
\revision{Since the objective function is linear in $(y_0,\BFy)$, it suffices to derive the convex hull of the feasible region~\cite[Lemma 1.3]{conforti2014integer}. Now, the result follows from Theorem~\ref{them:extended-formulation}.} \hfill \Halmos

\subsection{Proof of Theorem~\ref{them:polyLP}}
We start with showing that~(\ref{eq:variant-LP}) is polynomial-time solvable. Notice that~(\ref{eq:variant-LP}) has a linear objective function, and the separation problem of its feasible region can be solved using Algorithm~\ref{alg:sep}.   Therefore, due to the equivalence of separation and optimization~\citep[see Theorem 6.4.9 in][]{grotschel2012geometric},~(\ref{eq:variant-LP}) is polynomial time solvable.

Next, we present a proof that~(\ref{eq:variant-LP}) is a linear programming formulation of~(\ref{eq:mnl-ind}). The main idea is to use the optimal solution obtained from~\eqref{eq:variant-LP} to construct a feasible random assortment that satisfies the constraints of~\eqref{eq:mnl-ind}. Because \eqref{eq:variant-LP} is a relaxation of \eqref{eq:mnl-ind}, the average revenue of such randomized assortments should be no larger than the optimal objective value of \eqref{eq:variant-LP}. If the average revenue of the constructed randomized assortments is exactly equal to the optimal objective value of \eqref{eq:variant-LP}, then any possible realization of our constructed randomized assortments is the optimal assortment for~\eqref{eq:mnl-ind}. Therefore, it follows that \eqref{eq:variant-LP} can indeed serve as a linear programming formulation of \eqref{eq:variant-LP} because their optimal objectives are the same.

Let $(\bar{\x}, \bar{y}_0, \bar{\y})$ be an optimal solution of~(\ref{eq:variant-LP}) and let $r^*$ be the optimal value. First, it is evident that for each $j \in N$, at least one of the inequalities~\eqref{eq:stair-under} are tight. 
Suppose not, that is, there exists $j \in N$ such that none of inequalities~\eqref{eq:stair-under} is tight at $(\bar{x}_j,\bar{y}_0,\bar{\y})$.   Then, we can increase $\bar{x}_j$ by some $\epsilon >0$ without violating any constraint in \eqref{eq:stair-under}. Thus, the objective value of \eqref{eq:variant-LP} increases by $\sum_{i \in M} \alpha_i r_{ij} \theta_{ij} \epsilon$, contradicting with the optimality of $(\bar{\x}, \bar{y}_0, \bar{\y})$.

Second, we construct a randomized assortment based on $(\bar{\x}, \bar{y}_0, \bar{\y})$ that satisfies the constraints of~\eqref{eq:mnl-ind}.    In particular, let $\sigma$ be a permutation of $N$ such that $\bar{y}_{\sigma(1)} \geq \bar{y}_{\sigma(2)} \geq \cdots \geq \bar{y}_{\sigma(n)}$.  Let $\sigma^{-1}$ be the inverse of $\sigma$, that is, $\sigma^{-1}(j)$ denote the index $\tau$ such that $\sigma(\tau)=j$. Define 
\[
C_0 := \emptyset\quad \text{and} \quad C_j := \bigl\{\sigma(1), \ldots, \sigma(j) \bigr\} \text{ for } j  \in N. 
\]
Theorem~\ref{them:poly-sep-stair} implies that for each $j \in N$, the inequality~(\ref{eq:stair-under}) with $S = C_{\sigma^{-1}(j)-1}$ is tight at $(\bar{x}_j,\bar{y}_0,\bar{\y})$. Then, we construct a random assortment $\tilde{S}$ as follows:
\begin{equation*}
  \tilde{S} =
  \begin{dcases}
    C_0 & \textsf{ w.p. }  \lambda_0 := \frac{\bar{y}_0-\bar{y}_{\sigma(1)}}{\alpha(\emptyset)} \\
    C_j & \textsf{ w.p. } \lambda_k := \frac{\bar{y}_{\sigma(k)} -\bar{y}_{\sigma(k+1)} }{\alpha(S_k)} \text{ for } j \in N \setminus \{n\}\\
    C_n & \textsf{ w.p. } \lambda_n := \frac{\bar{y}_{\sigma(n)}}{\alpha(N)}.
  \end{dcases}
\end{equation*}
It is clear that $\lambda \geq 0$, and by~(\ref{eq:convexhull-lambda}) in the proof of Theorem~\ref{them:convexhull}, $\sum_{j}\lambda_j=1$. Moreover, for each $j\in \{0, \ldots, n\}$, $C_j$ satisfies the pre-order given by the graph $G$ because of the definition of $\sigma$ and $y_j \geq y_k$ for $(j,k) \in E$. Let $X(\tilde{S})$ be a random vector with $x_j(\tilde{S}) = 1$ if $j \in \tilde{S}$. It is clear that $X(\tilde{S})$ is feasible to \eqref{eq:mnl-ind}. Moreover, based on the definition of $C_k$, we can obtain 
\[
\mathbb{E}\bigl[x_j(\tilde{S})\bigr] = \pr\bigl(x_j(\tilde{S}) = 1 \bigr) = \sum_{k\leq\sigma^{-1}(j)-1} \lambda_k \cdot 0 + \sum_{k\geq \sigma^{-1}(j)}\lambda_k\cdot 1 = \bar{x}_j,
\]
where the second equality holds since product $j$ does not belong to $C_k$ if and only if $k \leq\sigma^{-1}(j)-1$, and the last equality holds due to~\eqref{eq:convexhull-cc2} in the proof of Theorem~\ref{them:convexhull}. In addition, we define another random variable  $Y(\tilde{S}) = \ext(C_k)$ with probability $\lambda_k$, and obtain that 
\[
\mathbb{E}\bigl[Y(\tilde{S})\bigr] = \sum_{k=0}^n\lambda_k \ext(C_k) =  (\bar{y}_0,\bar{\y}),  
\]
where the second equality follows from~\eqref{eq:convhull-cc1} in the proof of Theorem~\ref{them:convexhull}.   

Last, we compute the expected revenue of the random assortment. Let $\hat{\r}_0 = (0,u_{01}r_{01}, \ldots, u_{0n}r_{0n})$ and, for $i \in M$, let $\hat{\r}_i:=(\theta_{i1}r_{i1}, \ldots, \theta_{in}r_{in} )$. Now, the proof is complete by observing that 
\[
\begin{aligned}
\mathbb{E}[\alpha_0R_0^\text{MNL}(X(\tilde{S})) + \sum_{i \in M} \alpha_{i}R^{\text{IDM}}_{i}(X(\tilde{S}))] &=  \alpha_0 \hat{\r}_0^ \intercal \mathbb{E}[Y(\tilde{S})] + \sum_{i \in M}\alpha_i \hat{\r}_i \mathbb{E}[X(\tilde{S})] \\
& = \alpha_0\hat{\r}_0 (\bar{y}_0,\bar{\y}) + \sum_{i \in M}\alpha_i \hat{\r}_i \bar{\x} \\
& = r^*,
\end{aligned}
\]
where the first equality holds by the linearity of expectation, the second equality follows from the expectation of $X(\tilde{S})$ and $Y(\tilde{S})$ derived above, and the last equality holds since $(\bar{\x},\bar{y}_0,\bar{\y})$ is an optimal solution. 
\hfill \Halmos

\subsection{Proof of Corollary~\ref{cor:mixture}}
Consider a relaxation of~\eqref{eq:mixture} given as follows,
\begin{equation}\label{eq:mixture-rlx}
 \max_{} \Biggl\{ \sum_{i \in N}r_i\biggl( \alpha \frac{v_ix_{i}}{v_0 + \sum_{i \in N} v_ix_{i}}  +  (1-\alpha) \theta_iz_{i} \biggr) \Biggm| \x \geq \z ,\  \x \in \X \subseteq \{0,1\}^n ,\ \z \in \{0,1\}^n  \Biggr\} \tag{\textsc{Mixture-Rlx}}.
\end{equation}
If $\X$ is defined by precedence constraints then~\eqref{eq:mixture-rlx} is a special case of~\eqref{eq:mnl-ind}. Thus, by Theorem~\ref{them:polyLP},  inequalities~\eqref{eq:stair-under} yield a polynomial time solvable LP formulation for~\eqref{eq:mixture-rlx}. The proof is complete if we show that the objective value of~\eqref{eq:mixture}, denoted as $R$,  equals to the objective value of~\eqref{eq:mixture-rlx}, denoted as $R'$.  Clearly, $R'\geq R$. To argue $R'\leq R$, we consider an optimal solution $(\x^*,\z^*)$ of~\eqref{eq:mixture-rlx}, and observe that
\[
R' = \sum_{i \in N}r_i\biggl( \alpha \frac{v_ix^*_{i}}{v_0 + \sum_{i \in N} v_ix^*_{i}}  +  (1-\alpha) \theta_iz^*_{i} \biggr) \leq \sum_{i \in N}r_i\biggl( \alpha \frac{v_ix^*_{i}}{v_0 + \sum_{i \in N} v_ix^*_{i}}  +  (1-\alpha) \theta_ix^*_{i} \biggr) \leq R,
\]
where the first inequality holds due to $\x^* \geq \z^*$ and the non-negativity of $r_i$, $1-\alpha$ and $\theta_i$, and the second inequality holds since $\x^*$ is a feasible solution to~\eqref{eq:mixture}. \hfill \Halmos 



\section{Integer programming formulations for~(\ref{eq:omni-mnl})}

\subsection{Formulation~\textsf{Conic}}\label{section:conicformulation}
Here, we consider a special case of~(\ref{eq:omni-mnl}) where $\X_i = \{0,1\}^n$ for each $i \in M^+$ and use ideas in~\cite{sen2018conic} to derive formulation~\textsf{Conic} for~(\ref{eq:omni-mnl}). First, using the algebraic property of C-C transformation, we obtain an mixed-integer bilinear reformulation of~(\ref{eq:omni-mnl}), that is
\begin{equation*}~\label{eq:CC-QAP}
\begin{aligned}
\max \quad &  \sum_{i \in M^+} \alpha_i\biggl( \sum_{j \in N} r_{ij} u_{ij}y_{ij}  \biggr)\\
	\text{s.t.} \quad 	 &  (y_{i0},\y_i) \in \Pi_i(\X_i) \quad\text{for } i \in M^+ \\
 &\y_i =  y_{i0}\x_i \quad \text{and} \quad \x_i \in \X_i  \quad \text{for } i \in M^+  \\	
	 	 &\x_0 \geq \x_i  \quad \text{for } i \in M. \\
	\end{aligned}
\end{equation*}
Second, we use the McCormick inequalities~\citep{mccormick1976computability} to linearize the product of a continuous variable and a binary variable $y_{i0}x_j$ for each $i\in M^+$, with estimations about the lower bound $y^L_{i0}$ and the upper bound $y^U_{i0}$ on variable $y_{i0}$. More specifically, we define $y^L_{i0} := 1/\sum_{j=0}^nu_{ij} $ and $y^U_{i0} := 1/u_{i0}$. Third, following the approach in \cite{sen2018conic}, conic constraints corresponding to the bilinear terms $y_{ij} = y_{i0}x_j,\, \forall i\in M^+$ are also imposed, while the objective function needs to be rewritten. Thus, letting $ \tilde{r}_i = \max_{j \in N} r_{ij},\, \forall i \in M^+ $, we give the conic integer formulation \eqref{eq:conic-conic} for \eqref{eq:omni-mnl} as follows,

\begin{equation}\label{eq:conic-conic}
\begin{aligned}
\max \quad 
&  \sum_{i \in M^+} \alpha_i\biggl( \sum_{j \in N} \left(r_{ij}-\tilde{r}_{i}\right) u_{ij}y_{ij} + &{}&\tilde{r}_{i}\left(1-u_{i0}y_{i0}\right) \biggr)  \notag\\
\text{s.t.} \quad 	 
&x_{ij} \in \{0,1\} &&  \text{for  } i \in M^+ \text{ and } j \in N \\
&u_{i0} y_{i0} + \sum_{j \in  N} u_{ij} y_{ij} \ge 1  && \text{for  } i \in M^+\\
&0\leq y_{ij} \leq y_{i0}  && \text{for  } i\in M^+ \text{ and } j \in N \\
&w_{i} = u_{i0} + \sum_{j \in N} u_{ij}x_{ij}  && \text{for  } i \in M^+ \\
&w_{i} y_{i0} \ge 1  && \text{for  } i \in M^+\\
&w_{i}y_{ij} \ge x_{ij}^2  && \text{for  } i \in M^+ \text{ and } j \in N \\
& x_{0j} \ge x_{ij}   && \text{for  } i \in M \text{ and } j \in N  \\
& y_{ij} \geq  y^L_{i0} \cdot x_{ij}  && \text{for } i\in M^+ \text{ and }   j \in N \\
& y_{ij} \geq y_{i0} + y^U_{i0} \cdot x_{ij} -  y^U_{i0}  && \text{for } i\in M^+ \text{ and }   j \in N \\
& y_{ij} \leq  y_{i0} + y^L_{i0} \cdot x_{ij} -  y^L_{i0}  && \text{for  } i\in M^+ \text{ and } j \in N  \\
& y_{ij} \leq y^U_{i0} \cdot x_{ij}  && \text{for } i\in M^+ \text{ and } j \in N ,
\end{aligned}\tag{\textsc{Conic}}
\end{equation}
where the second constraint is obtained
by relaxing the constraint $u_{i0}y_{i0} + \sum_{j\in N} u_{ij}y_{ij} = 1 $ in $\Pi_i(\mathcal{F}_i), \, \forall i \in M^+$, which combined with the third constraint describes the convex hull of the choice probability set $\Pi_i(\mathcal{F}_i)$ for the case $\mathcal{F}_i=\{0,1\}^n$; it is sufficient to only use lower bounds on $ y_{ij} $ and $y_{i0}$ since the objective coefficients are nonpositive; the fourth through the sixth constraints are obtained from the conic quadratic relaxation of $y_{ij} = y_{i0}x_j,\, \forall i\in M^+$ as the approach introduced in \cite{sen2018conic}, and the last four constraints are derived using the McCormick inequalities. 


\subsection{Formulation~\textsf{MILP}}\label{section:milpformulation}
Here, we continue with the special case mentioned in Section \ref{section:conicformulation} and give a mixed-integer linear formulation \textsf{MILP} following the approach in \citep{bront2009column, MendezMiranda-BrontVulcanoZabala2014}.
The bilinear terms $y_{ij} = y_{i0}x_j,\, \forall i\in M^+$ can be linearized using the standard big-M approach. That is, for any bilinear term $y_{ij}=y_{i0}x_{ij}$, where $y_{i0}$ is continuous and nonnegative and $x_{ij}$ is binary, the following linear inequalities model a relaxation for the bilinear term: $0\le y_{ij} \le y_{i0}, ~y_{i0} - y_{ij} \le (1-x_{ij})U$ and $y_{ij} \le Ux_{ij}$, where $U$ is a sufficiently large upper bound on $y_{i0}$ and can be selected as $y^U_{i0} = 1/u_{i0}$. The mixed-integer linear formulation is given as follows,

\begin{equation}\label{eq:milp}
\begin{aligned}
\max \quad 
&  \sum_{i \in M^+} \alpha_i\biggl( \sum_{j \in N} r_{ij} u_{ij}y_{ij}  \biggr)  \notag\\
\text{s.t.} \quad 	 
&x_{ij} \in \{0,1\} \quad &&  \text{ for  } i \in M^+ \text{ and } j \in N \\
&u_{i0} y_{i0} + \sum_{j \in  N} u_{ij} y_{ij} = 1 \quad && \text{ for  } i \in M^+\\
&0 \le y_{ij} \le y_{i0} \quad &&  \text{ for  } i \in M^+ \text{ and } j \in N \\
&u_{i0} y_{ij} \le x_{ij} \quad &&  \text{ for  } i \in M^+ \text{ and } j \in N \\
&u_{i0}(y_{i0}-y_{ij})\le (1-x_{ij}) \quad &&  \text{ for  } i \in M^+ \text{ and } j \in N .
\end{aligned}\tag{\textsc{MILP}}
\end{equation}
\section{Generation of the data}\label{section:appendix_data_generation}

\subsection{Generation of the revenue and preference weight}\label{subsection:generation_revenue_prefernce}
Here, we introduce the generation of the data about revenue and preference weight of each product across the offline and online consumer segments. 

Firstly, the price of each product for the offline consumer segment is randomly drawn from a uniform distribution $U(10,20)$. The online segments are divided into two groups: the regular group and the VIP group. In the regular group, the price of all products is the same as that for the offline consumer segment. However, for the VIP group, we offer consumers a price discount. We assume that the discounted price is equal to the original price multiplied by a random variable drawn from the uniform distribution $U(0.8,1)$. 

Secondly, the preference weight of each product for the offline customer segment is randomly drawn from a uniform distribution $U(0,1)$. For the online segments, we assume that each segment has its own most preferred product, and its preference weight is set to $1$. Subsequently, all the remaining products for each online segment are assigned a preference weight randomly drawn from a uniform distribution $U(0,1)$. Throughout this section, we assume that the number of segments does not exceed the number of products. This allows us to set the most preferred product of each online segment to be unique. Furthermore, the preference weight for the no-purchase option are set to be $1$ for the offline consumer segment. However, for the online consumer segments, the preference weight for the no-purchase option varies, specifically set to $2$, $5$, and $10$. This difference is based on the observation that people from the offline channel tend to have a higher probability of making a purchase compared to those from the online channel.


\subsection{Generation of the data about the dominance relationship }\label{subsection:generation_luce_graph}
Here we give a detailed description about the generation of the directed acyclic graph for the 2SLM. Algorithm \ref{alg:luce_data_generation} generates $m$ random directed acyclic graphs and each graph indicates the dominance relationships among the products in each online consumer segment. For each online consumer segment, we assume that total $s\coloneqq \lfloor 0.25\cdot n \rfloor$ products are involved in the 2SLM. The remaining products are not influenced even if the consumer's behavior follows the 2SLM. Specifically, a subset of $s$ products is firstly sampled from $N$, then a random permutation is given for the selected products. Following the permutation, arcs are constructed by joining a product with at most three products from its $w$ successors. Each arc indicates the dominating relationship between its two vertices. In Algorithm \ref{alg:luce_data_generation}, we also set some conditions that help to construct the graph structure we want. That is, for each of the first $w=6$ products, at least one arc starts from it, while for each of the last $w$ products, only the arcs arriving at it are allowed. This guarantees that at least $w$ dominance relationships and at most $w$ dominated products exist for each online consumer segment. Algorithm \ref{alg:luce_data_generation} generates for each online consumer segment the directed graph without cycles. As a result, each pair of products has at most one arc joining them while may have multi distinct paths through them. Moreover, the more products involved, the more complex the directed acyclic graph.

\begin{algorithm}[h]
\fontsize{9pt}{10pt}\selectfont 
\renewcommand{\baselinestretch}{1.2}\selectfont
\caption{Generation of the directed graph for the 2SLM}\label{alg:luce_data_generation}
\KwData{products set $N$ and products number $n=|N|$; involved products numer $s\coloneqq \lfloor 0.25\cdot n \rfloor$; dominating width $w=6$.}
\KwResult{Directed acyclic graph $G(N, E^i),~ \forall i\in M$\;}
\For{$i$ in $M$}{
$E^i \gets \emptyset$ and $l \gets \lfloor 0.5\cdot w \rfloor $\;
$\{\mathcal{V}_1,\mathcal{V}_2, \dots, \mathcal{V}_s\}  \gets $  sorting of a random subset $\textsf{sample}(N, s)$\;
    \For{$v$ in $\{1, \dots, s-w\}$}{
    \eIf{$v<w$}
    {$k \gets $  $U[1,l]$\;}
    {$k \gets $  $U[0,l]$\;}
    $\mathcal{S} \gets $ random subset $\textsf{sample}(\{\mathcal{V}_{v+1}, \dots, \mathcal{V}_{v+w}\}, k)$\;
    $E^i \gets E^i \cup \{(\mathcal{V}_v, \mathcal{V}_{v^\prime})\}_{v^\prime \in \mathcal{S}}$ \;
    }
}
\end{algorithm}
For a given directed acyclic graph, by the graph theory, the reachability matrix defining the partial orders and the matrix defining the cover relation among the nodes can be obtained by graph manipulation on the adjacent matrix. The set of minimal elements is just composed of nodes with positive in-degree and zero out-degree. Then, we can formulate the convex hull of $\Pi(\mathcal{F})$ for {\ourC} (in Section \ref{section:CH}). We use the Python package \textsf{NetworkX} \footnote{https://networkx.org/} to help us tackle the necessary graph manipulations.

\section{Additional Numerical Studies}\label{sec:insight}
Thanks to the efficiency of {\our} and \textsf{CH-Chain}, in this subsection, we are able to provide some managerial insights regarding the quick-commerce assortment planning based on our solution approach \revision{through additional numerical studies}.

\subsection{Assortment Map}

\begin{figure}[htbp]
	\centering
	\begin{subfigure}{0.49\textwidth}
		\centering
            \captionsetup{font=normalsize}
            \includegraphics[width=0.8\linewidth]{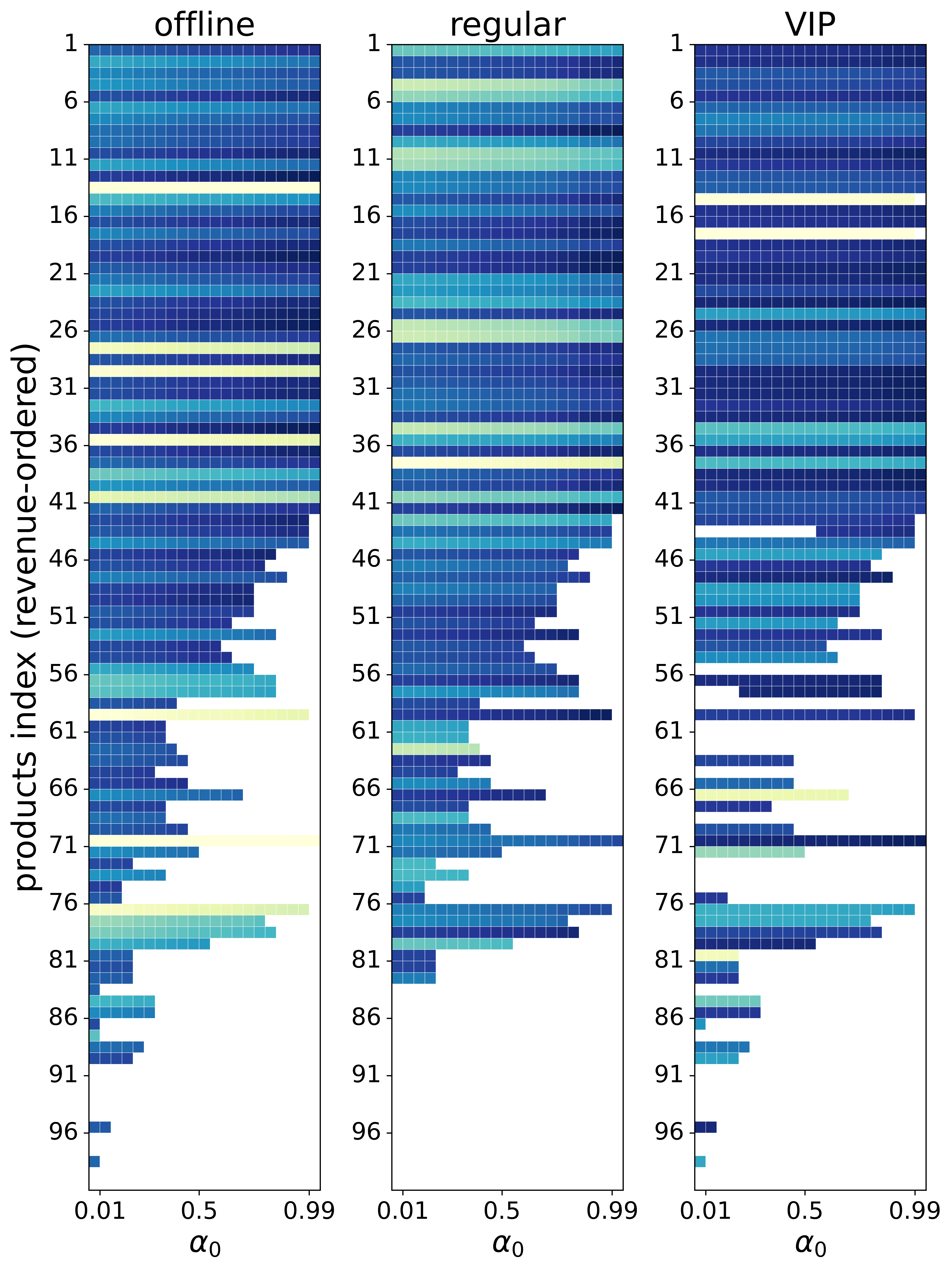}
		\caption{under MNL}
		\label{fig:map_assort_alpha0}
	\end{subfigure}
	\begin{subfigure}{0.49\textwidth}
		\centering
            \captionsetup{font=normalsize}
            \includegraphics[width=0.8\linewidth]{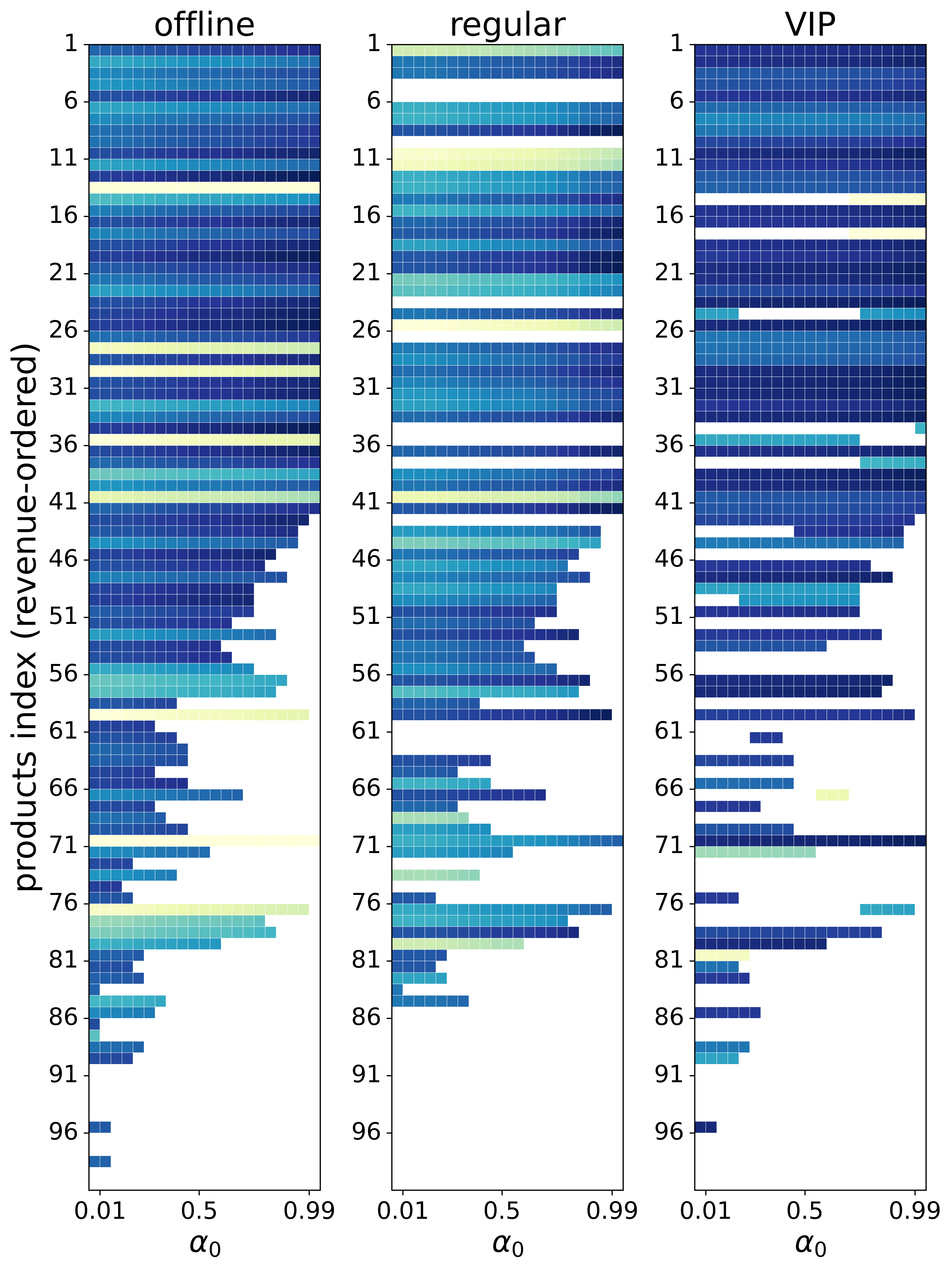}
		\caption{under 2SLM}
		\label{fig:map_assort_alpha0_Luce}
	\end{subfigure}
        \caption{\revision{Assortment maps by varying $\alpha_0$ with $ u^{\text{on}}_0=5 $}}
	\label{fig:MAP_ASSORT_vary_alpha0}
\end{figure}

\begin{figure}[htbp]
	\centering
	\begin{subfigure}{0.49\textwidth}
		\centering
            \captionsetup{font=normalsize}
            \includegraphics[width=0.8\linewidth]{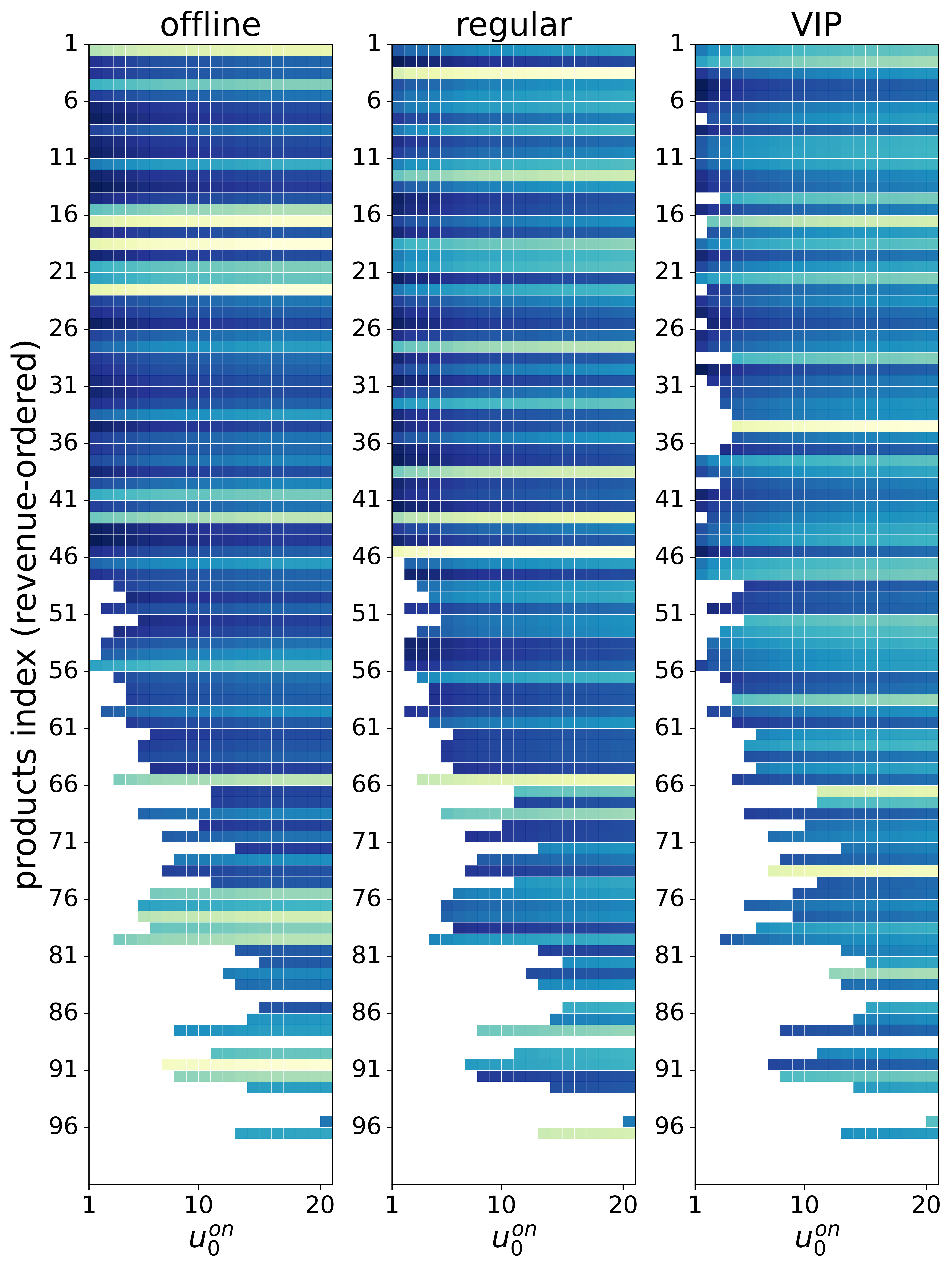} 
		\caption{under MNL}
		\label{fig:map_assort_v0}
	\end{subfigure}
	\begin{subfigure}{0.49\textwidth}
		\centering
            \captionsetup{font=normalsize}
            \includegraphics[width=0.8\linewidth]{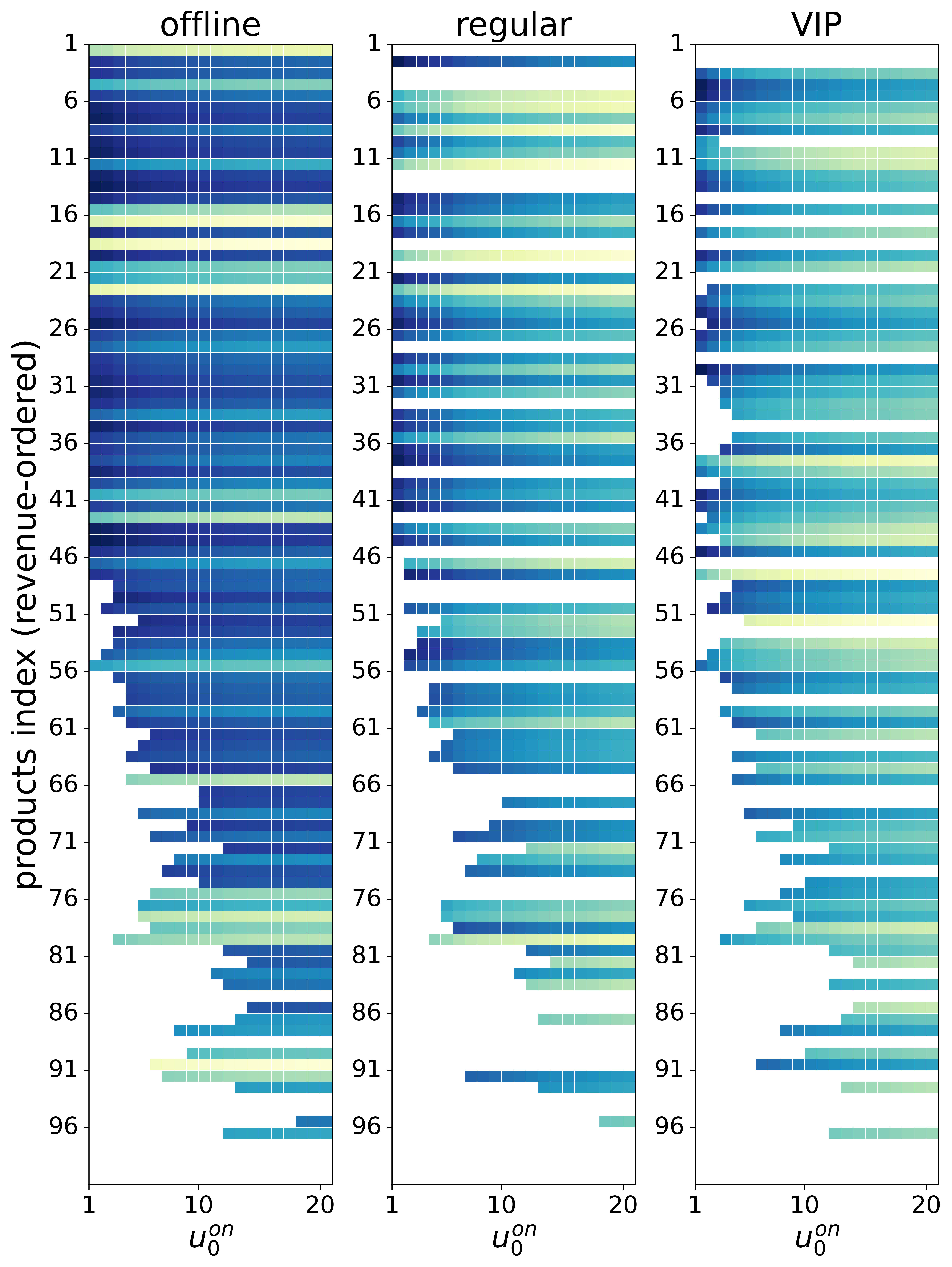} 
		\caption{under 2SLM}
		\label{fig:map_assort_v0_Luce}
	\end{subfigure}\\
         \caption{\revision{Assortment maps by varying $u^{\text{on}}_0$ with $\alpha_0=0.5$.}}
	\label{fig:MAP_ASSORT_vary_v0}
\end{figure}

We employ an \textsf{Assortment Map} to visually capture the impact of different parameters on the optimal assortment. Figure \ref{fig:MAP_ASSORT_vary_alpha0} presents the assortment map by varying the proportion of the offline channel $\alpha_0$, across a range from $0.01$ to $0.99$, using increments of $0.05$.  Additionally, Figure \ref{fig:MAP_ASSORT_vary_v0} displays the assortment map by varying the online consumers' preference weight on the no-purchase option, i.e., $u^{\text{on}}_0$, ranging from $1$ to $20$ with a step of $1$. 

In both figures, we visually depict the purchase probabilities for three distinct consumer segments: the offline consumer segment, a regular online consumer segment, and a VIP online consumer segment. These depictions are presented under two scenarios: one employing the MNL choice model and the other utilizing the two-stage Luce choice model. The y-axis in all figures corresponds to the product index, which is arranged according to the descending order of product price for the offline consumer segment, despite potential variations in price orders between the offline and VIP segments. A \yun{smaller} product index indicates a higher price for the corresponding product in the offline channel regardless of its price for the VIP segment. \yun{Only the results for the first $100$ products are disclosed because few of the remaining products are offered.} Each square block in the subfigures represents the purchase rate corresponding to a specific product under different values of a parameter.  Darker shades of blue indicate higher purchase probabilities. These assortment maps provide valuable insights into the influence of different parameters on the optimal assortment for the quick commerce assortment planning problem.

In Figure \ref{fig:MAP_ASSORT_vary_alpha0}, we observe that when $\alpha_0$ is close to 1, indicating a predominant share of traffic from the physical store, the offline assortment aligns with the ranking based on product price.  As  $\alpha_0$ diminishes and more traffic shifts to online channels, more products are activated in the optimal offline assortment. This change occurs because the decreased importance of the offline segment makes the cannibalization effect in the offline channel less relevant from the perspective of the entire market. As a result, it becomes more beneficial to provide a broader assortment of products in the offline channel to cater to the needs of different online segments. However, under smaller $\alpha_0$,  the optimal offline assortment departs from adhering solely to the product ranking based on price. This deviation underscores the necessity to account for the intricate influence of online segments more comprehensively on the assortment planning process.

In Figure~\ref{fig:MAP_ASSORT_vary_v0}, as the preference weight of the nonpurchase option for online segments increases, the cannibalization effect in the online channel becomes less significant. Consequently, it becomes more advantageous to offer a broader assortment of products for online segments. As a result,  more products are offered in the optimal assortment in the offline channel.

Furthermore, when comparing the case with the MNL model, it becomes evident that the pattern of the optimal assortment for the online segments becomes more intricate under the 2SLM, primarily due to the influence of dominant relationships among products. However, it is important to note that these dominant relationships among different online segments are generated independently, which subsequently limits their impact on the offline segment. As a consequence, assortment planning in the offline segment maintains a relatively stable nature.

\subsection{Benefit of the Joint Optimization}

Next, we will assess the benefits of jointly optimizing assortments from both the offline and online channels. In Section~\ref{subsec:2S-RO}, we introduce a two-step revenue-ordered policy (RO) procedure. This strategy involves initially determining the assortment for the offline channel and then utilizing this assortment as a candidate set for each online segment to select their respective optimal assortments. We denote the resulting revenue of the firm using the RO procedure as $rev^\text{RO}$. Furthermore, we use $rev^\text{opt}$ to represent the optimal revenue achieved by solving the joint optimization of quick-commerce assortment planning. To gauge the benefits of this joint optimization approach, we quantify the benefit of joint optimization as follows:
\begin{equation*}
Gain^\text{Joint} \coloneqq \frac{rev^\text{opt} - rev^\text{RO}}{rev^\text{opt}},
\end{equation*}

\begin{figure}[htbp]
\centering
\begin{subfigure}{0.49\textwidth}
  \centering
  \captionsetup{font=normalsize}
  \includegraphics[width=0.99\linewidth]{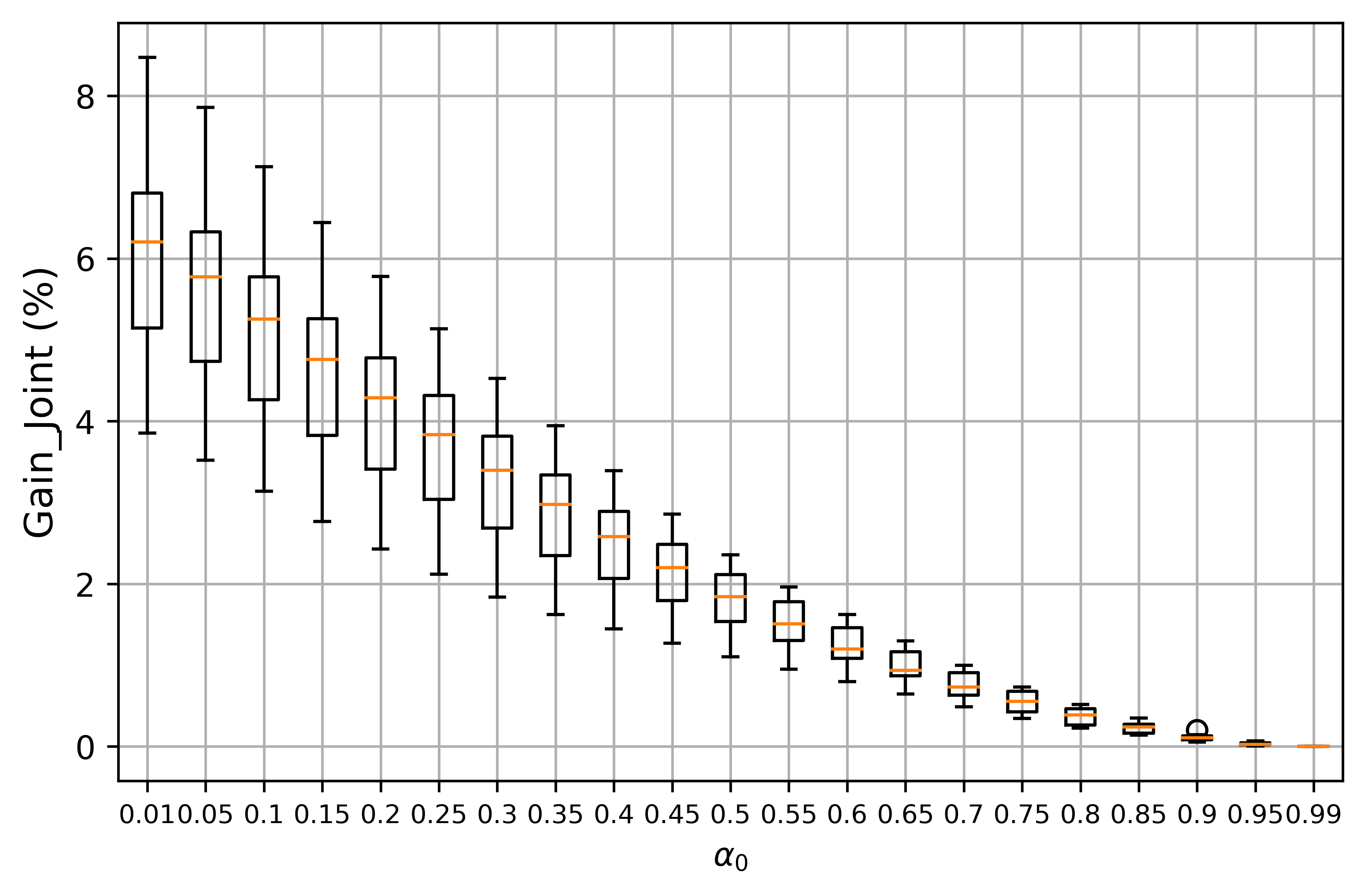}
  \caption{varying $\alpha_0$ with $u^{\text{on}}_{0}=5$}
    \label{fig:revenue_gain_SO_a0}
\end{subfigure}
\begin{subfigure}{0.49\textwidth}
    \centering
    \captionsetup{font=normalsize}
    \includegraphics[width=0.99\linewidth]{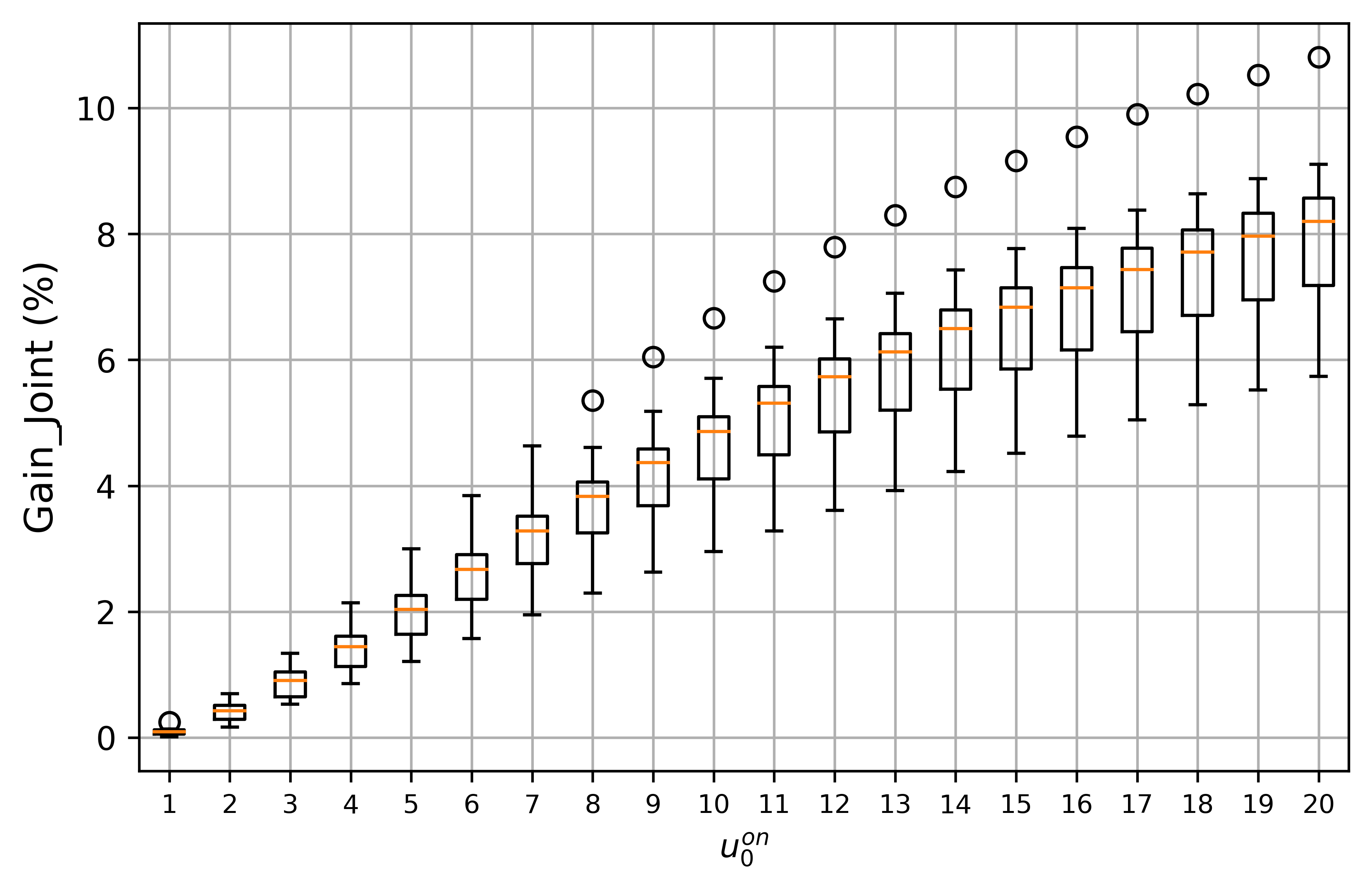}
    \caption{varying $u^{\text{on}}_{0}$ with $\alpha_0=0.5$}
    \label{fig:revenue_gain_SO_v0}
\end{subfigure}
\caption{The Benefit of Joint Optimization under $(n,m)=(200,100)$.}
\label{fig:REVENUE_SO}
\end{figure}

In Figure~\ref{fig:REVENUE_SO}, we present box plots of $Gain^\text{Joint}$ for each scenario using 12 random instances. Figure~\ref{fig:revenue_gain_SO_a0} delves into the impact of $\alpha_0$ on $Gain^\text{Joint}$. Notably, as $\alpha_0$ decreases, $Gain^\text{Joint}$ increases. This trend arises from the growing importance of accounting for the online channel's influence on the offline assortment, which becomes particularly pronounced as the market share of the online channel expands. This is further evidenced by Figure~\ref{fig:map_assort_alpha0},  where the optimal offline assortment is the smallest with $\alpha_0$ close to 1. As $\alpha_0$ decreases, the optimal offline assortment expands. This reflects how the joint optimization process incorporates the requirements of the online channel when devising the optimal offline assortment.

In Figure~\ref{fig:revenue_gain_SO_v0}, we explore the influence of $u^{\text{on}}_0$ on $Gain^\text{Joint}$. Clearly, as $u^{\text{on}}_0$ rises, so does $Gain^\text{Joint}$. This outcome can be attributed to the heightened need for more extensive assortments in the online segments due to a larger preference weight assigned to the non-purchase option. However, the offline assortment determined by the RO procedure neglects these attributes, and the resulting significant disparity between desired offline and online assortments leads to a higher revenue loss attributed to the RO policy. This is also evident from Figure~\ref{fig:map_assort_v0} that the optimal offline assortment expands when $u^{\text{on}}_0$ increases under the joint optimization approach.

%
%
%
%
%
%
%
%
%
%
\subsection{Revenue Loss of the Model Mis-Specification}
Leveraging the data accessibility in the online channel, the firm can pinpoint the dominant product relationships for each online segment. Unfortunately, these influential relationships are often disregarded during assortment planning. Let us delve into a scenario where this misspecification occurs: the firm develops an assortment strategy assuming that consumers, who actually follow the 2SLM, behave in the same manner as consumers in the conventional MNL model. Our primary objective is to quantitatively assess the resulting loss from this misspecification.

In cases of model misspecification, the firm tackles the joint optimization for the quick-commerce assortment planning problem under the presumption that consumer behavior adheres to the MNL model. Nevertheless, with the implemented assortment, consumers make selections in accordance with their inherent two-stage procedure following the 2SLM. We denote the resulting revenue arising from this model misspecification as $rev^\text{mis}$. Hence, the metric used to gauge the loss attributed to model misspecification can be expressed as follows:
\begin{equation*}
	Loss^{\text{mis}} \coloneqq \frac{rev^\text{opt} - rev^\text{mis}}{rev^\text{opt}},
\end{equation*}

\begin{figure}[htbp]
\centering
\begin{subfigure}{0.49\textwidth}
  \centering
  \captionsetup{font=normalsize}
  \includegraphics[width=0.99\linewidth]{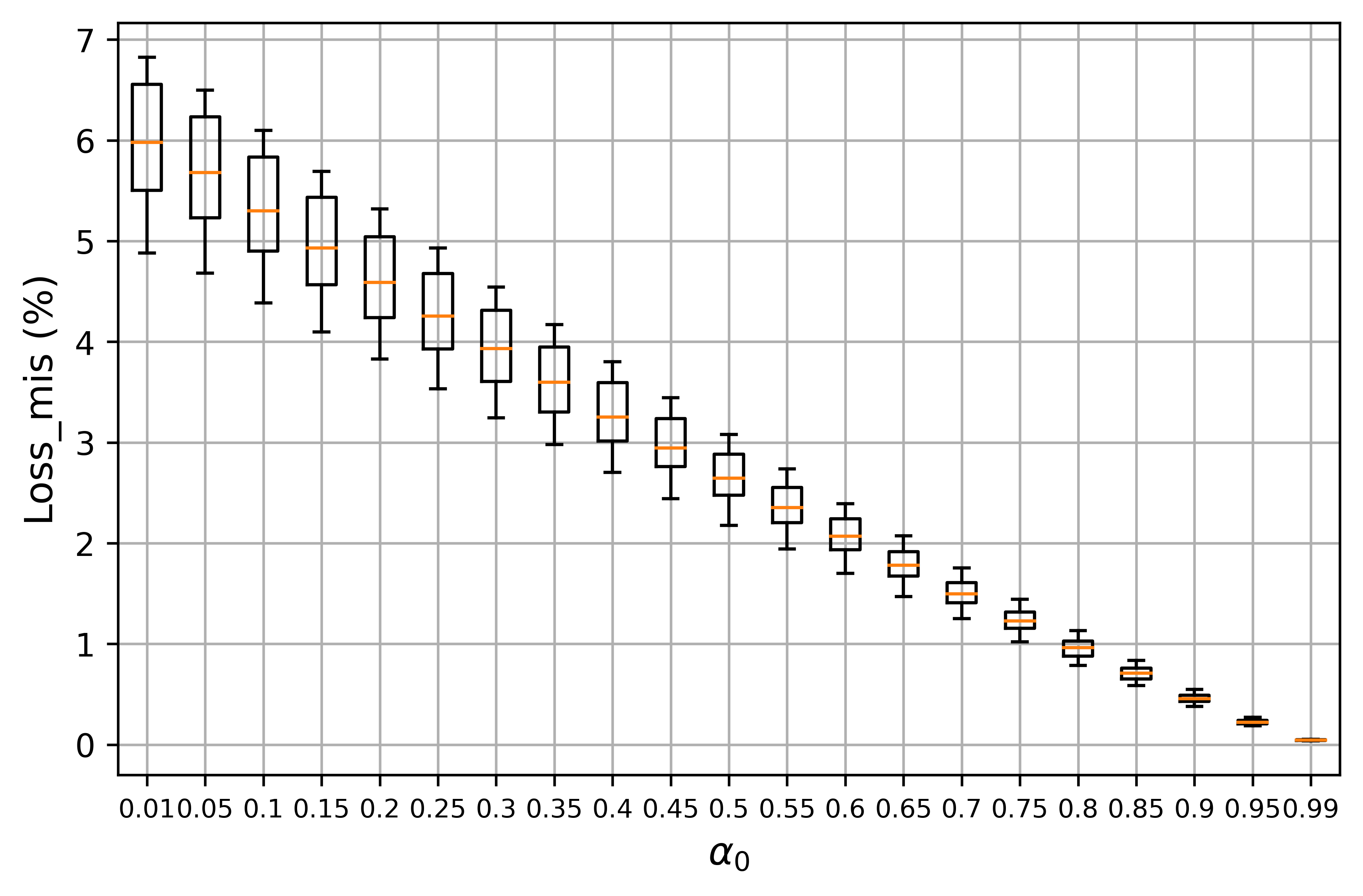}
  \caption{\normalsize varying $\alpha_0$ with $u^{\text{on}}_{0}=5$}
    \label{fig:revenue_loss_Luce_a0}
\end{subfigure}
\begin{subfigure}{0.49\textwidth}
    \centering
    \captionsetup{font=normalsize}
    \includegraphics[width=0.99\linewidth]{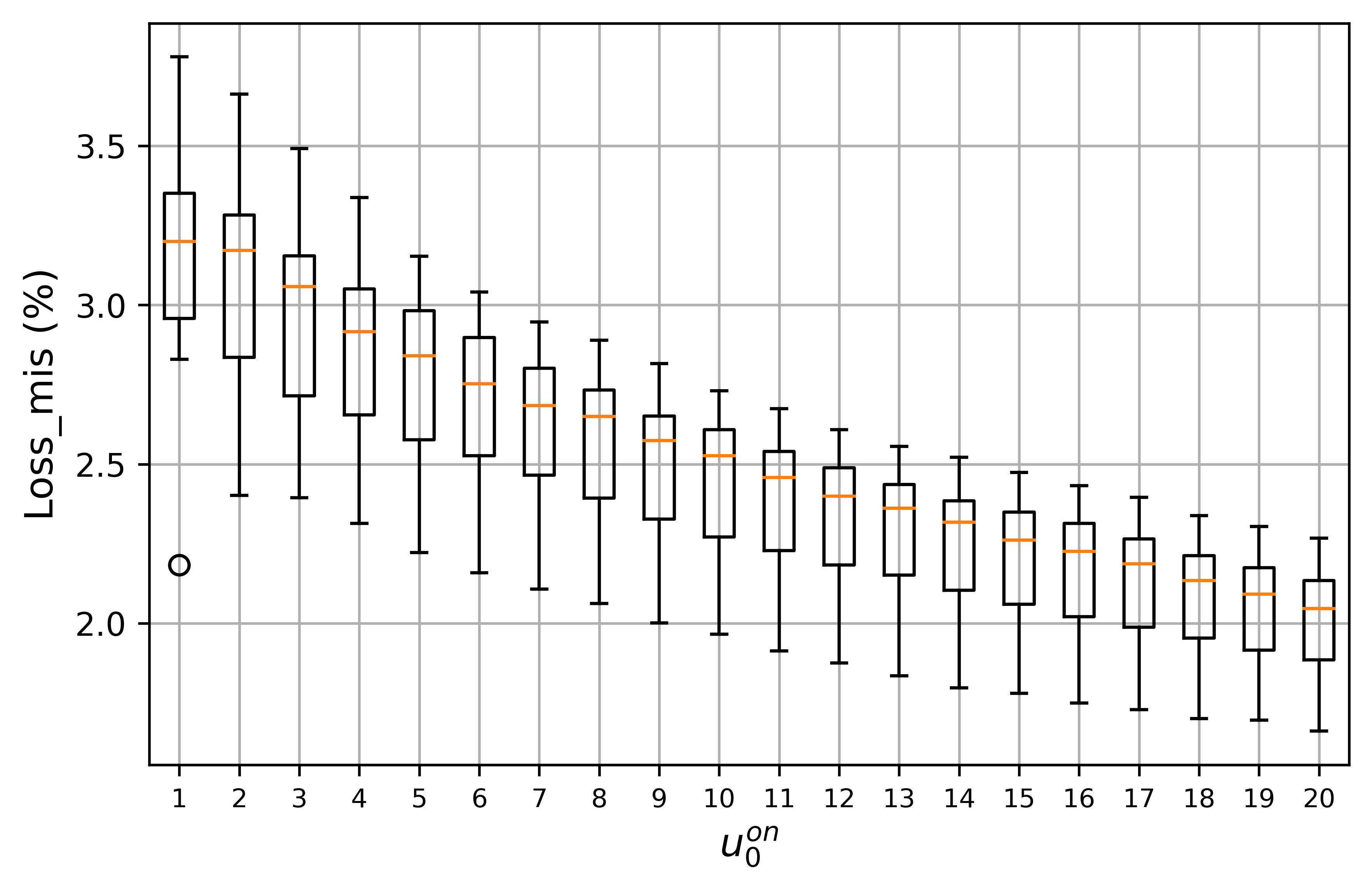}
    \caption{varying $u^{\text{on}}_{0}$ with $\alpha_0=0.5$}
    \label{fig:revenue_loss_Luce_v0}
\end{subfigure}
\caption{The Loss of Model Mis-Specification under $(n,m)=(200,100)$.}
\label{fig:REVENUE_Luce}
\end{figure}

The trends observed in Figure~\ref{fig:REVENUE_Luce} highlight that the loss amplifies with a decrease in $\alpha_0$ or $u^{\text{on}}_0$. This pattern can be elucidated as follows. As $\alpha_0$ decreases, the significance of online segments in the overall revenue generation becomes more pronounced. On the other hand, reducing $u^{\text{on}}_0$ leads to an increase in revenue from online segments due to the decreased likelihood of nonpurchase. Consequently, the contribution of online segments to the total revenue becomes more prominent. In both scenarios, neglecting the inherent behavior of online consumers results in elevated losses.

\section{Additional computational results}\label{section:app-computation}

\revision{The numerical experiments below were performed in the same environment as that used in Section \ref{section:computation}. That is, using Python 3.10 as the programming language on a Virtual Machine with 32 RAM and a 4-core Intel Core (Broadwell) @2.20 GHz processor, and using Gurobi 10.01 as the optimization solver. We accessed this environment for the additional numerical studies below from Aug. 2024 to Oct. 2024. (We accessed the environment for the numerical experiment in Section \ref{section:computation} from June. 2023 to Jul. 2023.)}


\subsection{\revision{The number of cutting-plane generation rounds}}\label{EC:cutinground} 

\revision{
Recall that Algorithm \ref{alg:cutting} adds  valid inequalities to the exact model $\textsf{CH-0}$, producing a tighter model $\textsf{CH-K}$. We address two questions regarding the impact of different values of $K$: (1) How much is the gap reduced with each additional round of adding new constraints? and (2) How does the computational performance vary with different values of $K$?}

\subsubsection{How much is the gap reduced with each additional round for adding new constraints?}
\begin{figure}[htbp]
    \centering
    \includegraphics[width=0.6\linewidth]{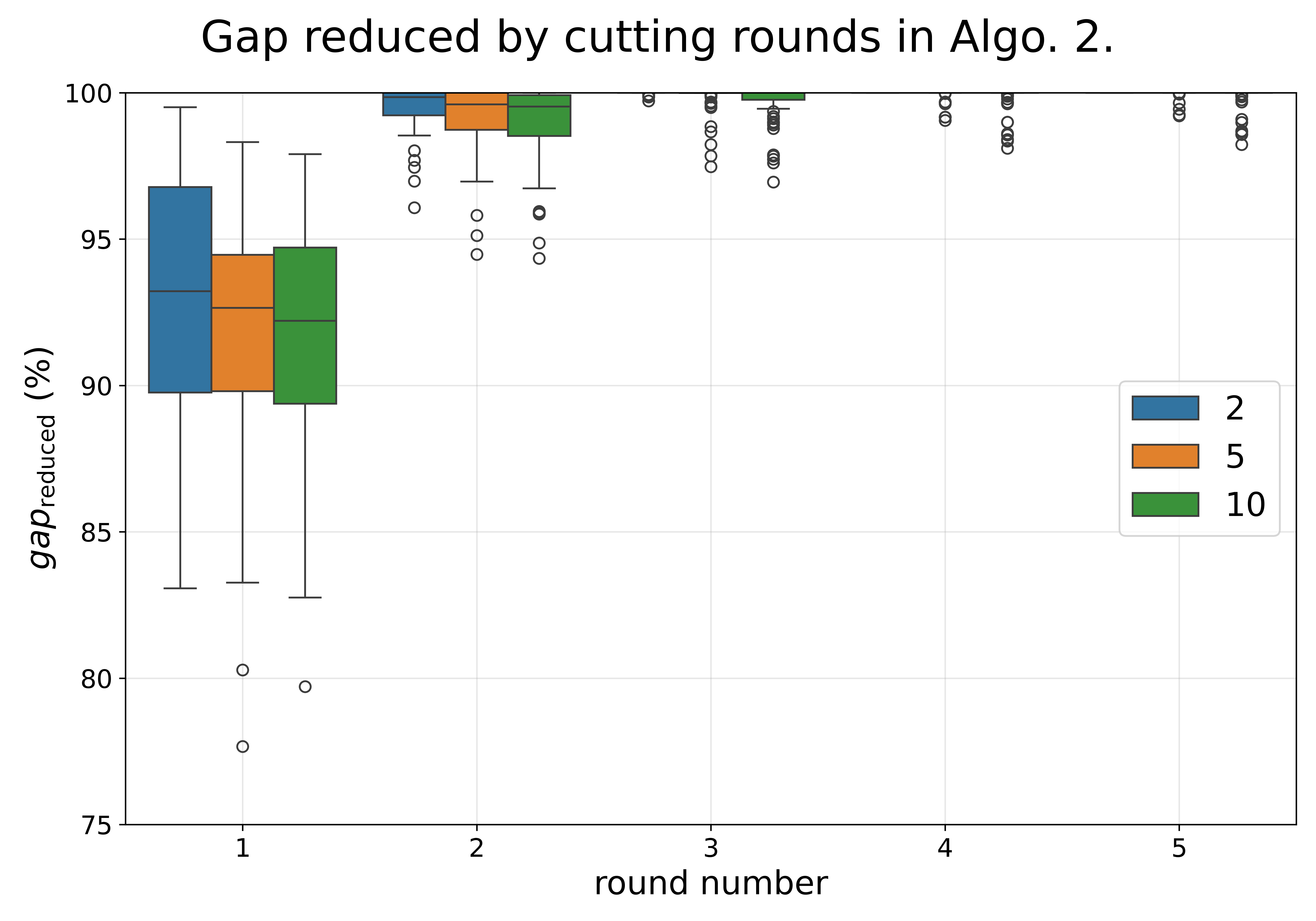}
    
    \caption{Gap reduced by the rounds of adding cutting planes.}
    \label{fig:gapClosedChangeWithK_compareK}
    \end{figure}

\revision{
    To address this, we set $K=5$ in  Algorithm \ref{alg:cutting} and calculated the reduced gap with each round, which is defined as 
    \begin{equation*}
        \text{gap}_{reduced}^{k} = \frac{\text{Obj}^0-\text{Obj}^k} 
        {\text{Obj}^0 - \text{Obj}^*}
    \end{equation*}
    where $\text{Obj}^k$ is the objective value of the relaxation of \texttt{CH-k} and $\text{Obj}^*$ is the optimal objective value. We test on $36$ random instances with size $(n,m)=(150,75)$.
    In Figure \ref{fig:gapClosedChangeWithK_compareK}, we plot the distribution of the reduced gap for parameter $u_{i0} \in \{2, 5, 10\}$ separately. 
    We observe that the average reduced gap exceeds $90\%$ after the first round and then increases to the proximity of $100\%$ as the number of cutting rounds grows. Additionally, the reduced gap stabilizes after the fourth round of adding the cutting planes on the tested instances.
    }

\subsubsection{\revision{How does the computational performance vary with different values of $K$?}}
\begin{figure}[htbp]
\centering
\includegraphics[width=0.6\linewidth]{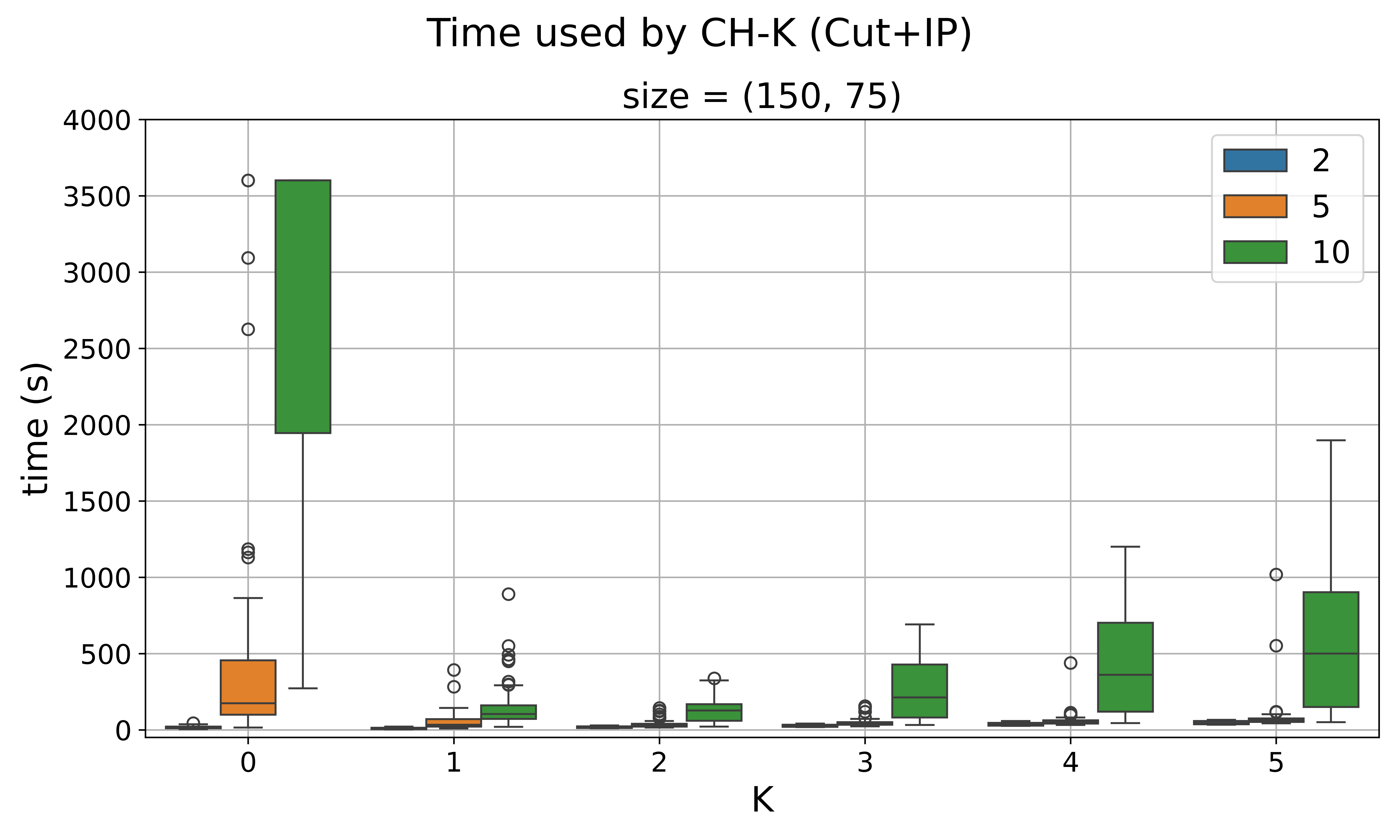}
\caption{Time used by \textsf{CH-K}.}
\label{fig:timeChangewithK}
\end{figure}

\revision{
    To address this, we set \revision{$K\in \{0,1,2,3,4,5\}$} and test each of them on $36$ instances with size $(n,m)=(150,75)$. For $K=0$, we directly solve the integer programming problem \textsf{CH-0}. For every \revision{$K\in \{1,2,3,4,5\}$}, we firstly add the corresponding \eqref{eq:stair-over} and \eqref{eq:stair-under} constraints by $K$ round(s) before solving the final integer programming problem \textsf{CH-K}. Figure \ref{fig:timeChangewithK} presents the distribution of the total computation time, including the time for adding cutting planes and solving the final integer programming problem \textsf{CH-K}. 
    First, Figure \ref{fig:timeChangewithK} shows that directly solving \textsf{CH-0} takes significant more time than our formulations, which introduce the cutting planes \eqref{eq:stair-over} and \eqref{eq:stair-under} into \textsf{CH-0}.
    Second, it appears that the total consumed time decreases with increasing $K\, (>0)$ when the preference weight on the no-purchase option is small (i.e., $u_{i0} \in \{2, 5\}$) but increases with increasing $K$ when the preference weight on the no-purchase option is large (i.e., $u_{i0} = 10$). However, this pattern is not consistent on all instances. While the total time required for the constraint generation (via Algorithm \ref{alg:cutting}) inevitably rises as $K$ increases, the continuous relaxation also becomes tighter due to the accumulation of additional valid constraints, which, in turn, increases the model's complexity. 
    \revision{To effectively handle varying values of the preference weight on the no-purchase option,} we select $K=2$ for use in Algorithm \ref{alg:cutting} across all the instances. 
    This choice represents a balance between the additional time required for constraint generation, the complexity introduced by these constraints, and the overall tightness of the formulation \textsf{CH-K}. 
    }

\subsection{Implementation via \textsc{Gurobi}'s \textsc{Callback} function}\label{EC:callback} 
\revision{
In this subsection, we use \textsc{Gurobi}'s \textsc{Callback} function to implement our inequalities~\eqref{eq:stair-over} and \eqref{eq:stair-under}, which we refer to as the \textsf{Callback} approach. \revision{Specifically, these inequalities are added at each current exploration node (\textsc{MIPNODE}\footnote{\url{www.gurobi.com/documentation/current/refman/cb_codes.html\#sec:CallbackCodes}}) with \textsc{OPTIMAL} status. In other words, whenever Gurobi solves a node to optimality, which is not always required during the branching process,  \textsf{Callback} is triggered.}}

\revision{
We evaluate the \textsf{Callback} approach on the same instances used in Section \ref{subsec:numerical_result} (specifically, Table \ref{tab:Pref_NoLuce_Luce}). 
In Table \ref{table:callback}, \revision{we  report the average, minimum, maximum, and standard deviation of the solving time across the 36 instances for each configuration, denoted as ``Ave'', ``Min'', ``Max'', ``Std'', respectively. Additionally, we include the number of instances solved within a 3600-second time limit, denoted by ``\#''.} For the case that some instances are not solved successfully within the time limit 3600 seconds, we omit the standard deviation (denoted by ``--'').}

\revision{
Compared with Table \ref{tab:Pref_NoLuce_Luce} in Section \ref{subsec:numerical_result}, Table \ref{table:callback} reveals two keys findings. First, the \textsf{Callback} approach outperforms both the {\conic} (or {\conicC}) and {\milp} (or {\milpC}) approaches in terms of average solving time and the number of solved instances. \revision{Second, it shows that the \textsf{Callback} approach exhibits significantly longer average solving time and greater variability compared to our formulation {\our} (or {\ourC}). As demonstrated in Appendix~\ref{EC:cutinground}, setting $K=2$ in Algorithm~\ref{alg:cutting} significantly reduces the gap. Thus,  the additional time required to generate cuts during the  \textsf{Callback} procedure at non-root nodes outweighs the benefit of further improving the gap at those nodes.}
}

\begin{table}[htbp]
\fontsize{9pt}{10pt}\selectfont 
\renewcommand\tabcolsep{2pt}
\renewcommand{\arraystretch}{1.1}
\caption{Computational performance of the CALLBACK-approach}
\centering
\revision{
\begin{tabular}{cccccccccccc}
\toprule
\multirow{2}{*}{(n,m)} & \multirow{2}{*}{$u_{i0}$} & \multicolumn{5}{c}{\textbf{without 2SLM}} & \multicolumn{5}{c}{\textbf{with 2SLM}}\\
\cmidrule(lr){3-7} \cmidrule(lr){8-12}
& & Ave & Min & Max  & Std & \# & Ave & Min & Max & Std & \# \\
\hline 
\multirow{3}{*}{(100,50)}  & 2 & 12.5 & 4.1  & 29.6 & 6.2 & 36  & 11.4 & 4.4  & 22.6 & 5.7 & 36 \\
  & 5 & 38.2 & 13.1 & 109.8 & 23.7 & 36 & 40.2 & 9.8  & 158.8 & 29.5 & 36 \\
  & 10 & 70.5 & 20.8 & 136.3 & 32.7 & 36 & 92.4 & 25.8 & 427.7 & 73.8 & 36 \\
\hline
\multirow{3}{*}{(150,75)}  & 2 & 46.5 & 31.2 & 91.9 & 14.4 & 36 & 41.8 & 22.9 & 65.3 & 10.8 & 36 \\
  & 5 & 172.4 & 60.0 & 633.5 & 115.0 & 36 & 173.9 & 60.0 & 693.2  & 134.7 & 36 \\
  & 10 & 495.2 & 181.7 & 1859.3 & 400.0 & 36 & 541.2 & 141.1 & 1354.6 & 340.4 & 36 \\
\hline
\multirow{3}{*}{(200,100)} & 2 & 115.3 & 34.9 & 501.9 & 82.3 & 36 & 177.9 & 37.6 & 400.1 & 95.7 & 36 \\
 & 5 & 754.0 & 206.2 & 3600 & -- & 34 & 776.6 & 160.5 & 3600 & -- & 35 \\
 & 10 & 2,783.6 & 435.6 & 3600 & -- & 16 & 2799.2   & 951.0 & 3600 & -- & 15 \\
\bottomrule
\end{tabular}\label{table:callback}
}
\end{table}

\subsection{On the joint assortment and personalization problem}\label{EC:personalization} 
\revision{
In this subsection, we present numerical results on the joint assortment and personalization problem proposed in~\cite{HousniTopaloglu2022}. The preference weights and the prices are generated following the approach in \cite[Sec. 7.1.2]{HousniTopaloglu2022}. Additionally, we consider a class of larger instances with $(n,m)=(200,500)$. Following the same setup as in \cite[Sec. 7.1.2]{HousniTopaloglu2022}, we impose the cardinality constraint and vary the upper bound on the number of products offered to the offline customer segment over $K \in \{\frac{1}{10}n, \frac{1}{2}n\}$. Each configuration is repeated $36$ times randomly. In Table \ref{tab:personalization}, we report the average, minimum, maximum and standard deviation of the solving time across the $36$ instances, denoted as ``Time'', ``Min'', ``Max'', and ``Std'', respectively. The result demonstrates that our method efficiently solves the joint assortment and personalization problem for the instances with the same size as those in \cite{HousniTopaloglu2022}, as well as the instances with much larger size. Furthermore, the small standard deviations in Table \ref{tab:personalization} highlight the robustness of our method.
 }
\begin{table}[htbp]
    \fontsize{9pt}{10pt}\selectfont 
    \renewcommand\tabcolsep{2pt}
    \renewcommand{\arraystretch}{1.1}
    \caption{Testing on joint assortment and personalization problem}
	\label{tab:personalization}
	\centering
 \revision{
    \begin{tabular}{p{5em}ccccccc}
    \toprule
    Instances & m & n  & K & Time & Min & Max & Std\\
    \hline
    \multirow{12}{*}{\makecell{Same size \\as in \\literature}} 
    & \multirow{4}{*}{10} & \multirow{2}{*}{50} & 5 & 0.59   & 0.49  & 0.68   & 0.05 \\
    & & & 25 & 0.63   & 0.57   & 0.75   & 0.04 \\
    \cmidrule(l){3-8}
    & & \multirow{2}{*}{100} & 10 & 1.21   & 1.08   & 1.65   & 0.11 \\
    & & & 50 & 0.96   & 0.84   & 1.12   & 0.07 \\
    \cmidrule(l){2-8}
    & \multirow{4}{*}{50} & \multirow{2}{*}{50} & 5 & 2.29   & 2.02   & 2.78   & 0.22 \\
    & & & 25 & 1.81   & 1.65   & 2.10   & 0.10 \\
    \cmidrule(l){3-8}
    & & \multirow{2}{*}{100} & 10 & 4.69   & 4.21   & 5.26   & 0.26 \\
    & & & 50 & 3.83   & 3.48   & 4.22   & 0.18 \\
    \cmidrule(l){2-8}
    & \multirow{4}{*}{100} & \multirow{2}{*}{50} & 5 & 4.75   & 4.23   & 7.42   & 0.69 \\
    & & & 25 & 3.40   & 3.13   & 3.77   & 0.19 \\
    \cmidrule(l){3-8}
    & & \multirow{2}{*}{100} & 10 & 9.52  & 8.39   & 11.11   & 0.69 \\
    & & & 50 & 7.17   & 6.68   & 7.79   & 0.30 \\
    \hline
    \multirow{2}{*}{Larger size} & \multirow{2}{*}{500} & \multirow{2}{*}{200} & 20 & 121.11 & 113.62 & 135.87 & 4.36 \\
    &  &  & 100 & 104.31 & 99.18 & 110.54  & 2.59\\
    \bottomrule
    \end{tabular}
}
\end{table}

\subsection{Inventory Rounding Procedure for \eqref{eq:omni-mnl}: Method and Numerical Results}
\label{ec:inventory_simulation}

{\color{black}
\rvtwo{We consider a setting where $T$ consumers arrive each day. The retailer must determine both the assortments and the inventory target levels for the products included in the assortments. Let $V_{\text{fluid}}$ denote the retailer’s maximum expected profit from these $T$ consumers under a hypothetical make-to-order setting. Let $V_{\text{simulation}}$ represent the average realized profit under a make-to-stock setting, where inventory constraints are active. In Algorithm~\ref{alg:inventory-simulation}, $V_{\text{fluid}}$ is computed based on the formulation in \eqref{eq:omni-mnl}. We also propose an inventory rounding procedure to determine the inventory target levels.}

\rvtwo{To incorporate inventory costs, we denote $c_j$ as the unit ordering cost for product $j$, and set the effective price of product $j$ for consumer segment $i$ as $r_{ij} + c_j$. When solving the assortment optimization problem, we continue to use $r_{ij}$ to represent the unit ``revenue'' for product $j$ and segment $i$, even though it effectively corresponds to the unit profit. The main components of the inventory rounding procedure and the simulation logic are summarized below:}
\begin{enumerate}
    \item \textbf{\rvtwo{Determining purchase probability for each product}.}  
    Solving the complete problem \eqref{eq:omni-mnl} provides the optimal decision variables \( y^*_{ij} \), from which we derive the purchase probability for product \( j \) as  
    \(
    \beta_j = \sum_{i \in M^+} \alpha_i u_{ij} y^*_{ij},
    \) 
    where \( \alpha_i \) is the arrival probability of consumer segment \( i \), and \( u_{ij} \) is the utility of product \( j \) to segment \( i \).  The corresponding optimal objective value of \eqref{eq:omni-mnl} is denoted by $V^*$.

    \item \textbf{\rvtwo{Calculating expected demand}.}  
    Assume the total time horizon consists of \( T \) periods. Then, the expected demand (assuming there is no stockout) for product \( j \) is  
    \(Q_j = T \beta_j\),
    and the fluid approximation revenue over the $T$ periods is defined as $V_{\text{fluid}} = T\times V^*$ \citep{zhang2024leveraging}.
    
    \item \textbf{Rounding procedure.} 
    To obtain the implantable inventory decisions, fractional order quantities are rounded to integers using the approach proposed by \cite{zhang2024leveraging}. Let \( \lfloor Q_j \rfloor \) and \( \lceil Q_j \rceil \) denote the floor {(rounded down value)} and ceil (rounded {up} value) of \( Q_j \), respectively. \rvtwo{Define $\delta = \left\lceil \sum_j Q_j \right\rceil - \sum_j \lfloor Q_j \rfloor$, which represents the number of products that must be rounded up to ensure that the total inventory aligns with the total expected demands.}
    
    Since product prices may vary across consumer segments, we define the \emph{weighted price} of each product \( j \) as  
    \(
    \tilde{r}_j = \sum_i \alpha_i r_{ij}
    \). Then rank the products in the optimal assortment based on their weighted prices. Let \( \bar{S}^* \) be the set of top $\delta$
    products (by weighted price) within the optimal assortment.
    Then, the final integer order quantity \( \tilde{Q}_j \) for product \( j \) is defined as  
    \[
    \tilde{Q}_j = 
    \begin{cases}
        \lceil Q_j \rceil & \text{if } j \in \bar{S}^*, \\
        \lfloor Q_j \rfloor & \text{otherwise}.
    \end{cases}
    \]
    \item \textbf{Simulation.} According to the solution to the complete problem, we set the assortments as $S = \{S_0, S_1, \dots, S_M\}$ where $S_i = \{j\in N | y^*_{ij} > 0\}$ represents the assortment for the customer segment $i\in M^+$. \rvtwo{We then conduct a Monte Carlo simulation. For each sample path, the process proceeds as follows}:
    \begin{enumerate}
        \item \label{simulate_period_t} At the beginning of each period $t$, denote the inventory vector by $I^t = [I^t_1, \dots, I^t_N]$ where $I^t_j$ is the inventory of the product $j$ and let $A^t(I^t) = \{j \in N | I^t_j >0 , \forall j \in N\}$ be the available product set for period $t$.
        \item Randomly draw a customer from $M^+$ with the distribution $\{\alpha_0, \alpha_1, \dots, \alpha_M\}$. Suppose the customer is from segment $i$. The customer chooses a product $j \in \tilde{S}^t_i = S_i \cap A^t(I^t)$ or opts not to purchase (denoted $j=0$), following the standard multinomial logit (MNL) model. If $j \ne 0$, the firm earns a revenue of $p_{ij} = r_{ij} + c_j$.
        \item Update the inventory for the next period $t+1$ as 
        \[
        I^{t+1}_j = \begin{cases}
            I^{t}_j-1, & \text{ if } j=d^{t},\\
            I^{t}_j, & \text{ otherwise},
        \end{cases}
        \]
        where $d^{t}$ is the choice of the customer in period $t$. Then go to \eqref{simulate_period_t} until the end of period $T$.
        \item Finally, the company’s net revenue over the selling horizon is calculated as the total accumulated revenue across the $T$ periods minus the total cost of the initial inventory.
    \end{enumerate}
\end{enumerate}
}

\begin{algorithm}[hb]
\fontsize{9pt}{12pt}\selectfont 
\caption{Inventory Simulation}\label{alg:inventory-simulation}
\KwData{$T>0$; $u_{ij}$ for $i\in M^+, j\in N\cup\{0\}$; $r_{ij}$, for $i\in M^+, j\in N$; $c_{j}$ for $j\in N$\;}
\KwResult{$V_{\text{fluid}}$, $V_{\text{simulation}}$\;}
$(\{y^*_{ij}\}_{i\in M^+,\, j\in N}, V^*) \gets$ solve the QAP\;
$V_{\text{fluid}} \gets T \times V^*$\; 
$\beta_j \gets \sum_{i\in M^+} \alpha_i u_{ij} y^*_{ij}$, for each $j\in N$\;
$Q_j \gets T \times \beta_j$, for each $j\in N$\; 
$\delta \gets \lceil \sum_{j\in N} Q_j\rceil - \sum_{j\in N}\lfloor Q_j\rfloor$\;
$\tilde{r}_j \gets \sum_{i\in M^+} \alpha_i r_{ij}$, for each $j\in N$\;
$\bar{S}^* \gets $ the index of top $\delta$ products within $\{j\in N | y_{0j}^* > 0\}$ by weighted price $\tilde{r}_j$\;
$\tilde{Q}_j \gets \begin{cases}
        \lceil Q_j \rceil & \text{if } j \in \bar{S}^*, \\
        \lfloor Q_j \rfloor & \text{otherwise}.
    \end{cases}$\;
\For{$l = 1:1000$}
{
$t \gets 0$; $I^0 \gets \tilde{Q}$; $d^0 \gets 0$; $V^l_{\text{simulation}} \gets -\sum_{j\in N} \tilde{Q}_j c_j$\;
$S_i \gets \{j\in N | y^*_{ij} > 0\}$, for each $i\in M^+$\;
\While{$t\le T$}{
    $I^{t+1} \gets  \begin{cases}
            I^{t}_j-1, & \text{ if } j=d^{t},\\
            I^{t}_j, & \text{ otherwise},
        \end{cases}$\;
    $t \gets t+1$\;
    $A^{t} \gets \{j \in N | I^{t}_j > 0\}$\;
    $i \gets $ randomly drawn from $M^+$ with distribution $\{\alpha_0, \alpha_1, \dots, \alpha_M\}$\;
    $\tilde{S}^t_i \gets S_i \cap A^t$\;
    $d^t \gets $ randomly drawn from $\tilde{S}^t\cup\{0\}$ with distribution $\{\frac{u_{ij}}{u_{i0} + \sum_{j^\prime\in \tilde{S}^t}u_{ij^\prime}}\}_{j\in \tilde{S}^t} \cup \{\frac{u_{i0}}{u_{i0} + \sum_{j^\prime\in \tilde{S}^t}u_{ij^\prime}}\}$\;
    $V^l_{\text{simulation}} \gets V^l_{\text{simulation}} + (r_{id^t} + c_{d^t})$\;
    
}
}
$V_{\text{simulation}} \gets \frac{1}{1000}\sum_{l=1}^{1000} V^l_{\text{simulation}}$\;
\end{algorithm}

\begin{table}[hbpt]
\fontsize{9pt}{10pt}\selectfont 
\renewcommand\tabcolsep{3pt}
\renewcommand{\arraystretch}{1.1}
\centering
\caption{Inventory simulation results for varying $u_{i0}$ and $T$.}
\label{tab:inventory-simulation}
\begin{tabular}{cccrrc}
\toprule
$(n,m)$ & $u_{i0}$ & $T$ & $V_\text{fluid}$ & $V_\text{simulation}$ & $\frac{V_\text{fluid} - V_\text{simulation}}{V_\text{fluid}} \times 100$ \\
\midrule
(100, 50) & 2  & 500  & 8168.83  & 7955.38  & 2.61  \\
(100, 50) & 2  & 1000 & 16337.66 & 16037.05 & 1.84  \\
(100, 50) & 2  & 2000 & 32675.32 & 32253.11 & 1.29  \\
(100, 50) & 5  & 500  & 7635.56  & 7379.38  & 3.36  \\
(100, 50) & 5  & 1000 & 15271.13 & 14905.74 & 2.39  \\
(100, 50) & 5  & 2000 & 30542.26 & 30032.24 & 1.67  \\
(100, 50) & 10 & 500  & 7103.47  & 6799.81  & 4.27  \\
(100, 50) & 10 & 1000 & 14206.94 & 13779.27 & 3.01  \\
(100, 50) & 10 & 2000 & 28413.87 & 27815.19 & 2.11  \\
\bottomrule
\end{tabular}
\end{table}

\rvtwo{\rvtwo{We conduct the simulation on instances randomly generated using the same data generation procedure as in Section \ref{subsec:num-setting}. Additionally, we vary the selling horizon as $T\in \{500, 1000, 2000\}$  and the number of sample paths is fixed at 1,000.} For the inventory cost, we set $c_j=1$ for all $j\in N$.}
\rvtwo{From Table~\ref{tab:inventory-simulation}, we observe that the relative gap between $V_{\text{fluid}}$ and $V_{\text{simulation}}$ decreases as $T$ increases. This pattern is consistent with the findings reported in \cite{zhang2024leveraging}. Additionally, when the preference weight of the no-purchase option for online segments increases, the actual purchase rate of each consumer declines. Consequently, the relative gap between $V_{\text{fluid}}$ and $V_{\text{simulation}}$ widens.}

\subsection{Comparison with Improved Revenue-Ordered Heuristics}\label{section:comparisonwithRO}

\rvtwo{In this subsection, we propose a heuristic that improves upon the two-step revenue-ordered policy introduced in Section~\ref{subsec:2S-RO}. In the first step, the offline assortment, which defines the common pool of products available to all customers, is constructed depending on whether a cardinality constraint is present. If such a constraint exists, the offline assortment is directly optimized using a standard MNL assortment model with a cardinality constraint. Otherwise, an enumeration-based approach is applied: products are ranked by their offline revenue potential, and candidate assortments are evaluated iteratively to identify the one that maximizes the combined offline and online revenue.}

\rvtwo{Given an offline assortment, the second step focuses on online personalization. For each consumer segment, the algorithm evaluates nested subsets of the offline set, ordered by segment-specific revenues, and selects the subset that yields the highest expected revenue. The consumer choice model is assumed to follow either the MNL or the 2SLM, depending on the setting.}

\rvtwo{Compared with the two-step revenue-ordered policy in Section~\ref{subsec:2S-RO}, Algorithm~\ref{alg:2SRO-heuristics} expands the search space while still leveraging revenue-based ordering to limit computational complexity. As a result, it has the potential to improve performance without sacrificing tractability. For implementation details, please refer to Algorithm~\ref{alg:2SRO-heuristics}.}

\rvtwo{
\begin{example}\label{example:improved_RO_is_bad}
 Before going into the numerical experiments, we first introduce a toy example with three products and two consumer types in the online channel. The arrival rates of the three types are $\alpha = (0.7, 0.2, 0.1)$, and the revenues of products are $(r_{01}, r_{02}, r_{03}) = (20, 14, 14)$, $(r_{11}, r_{12}, r_{13}) = (10,10,18) $, and $(r_{21}, r_{22}, r_{23}) =(10, 10, 20)$. The preference weight of products are $(u_{01}, u_{02}, u_{03}) = (100, 200, 1)$, $(u_{11}, u_{12}, u_{13}) = (1, 1, 100)$ and $(u_{21}, u_{22}, u_{23}) =(2, 2, 1)$, and the preference weights of the no-purchase option is $u_{00}=u_{10}=u_{20}=1$. The optimal gap of the \eqref{eq:2SNR} and the improved heuristic using Algorithm \ref{alg:2SRO-heuristics} are $\frac{18.39-15.53}{18.39} = 15.54\%$ and $\frac{18.39-15.72}{18.39} = 14.48\%$, respectively. \hfill \Halmos
\begin{table}[H]
\fontsize{9pt}{10pt}\selectfont 
\renewcommand\tabcolsep{1.5pt}
\renewcommand{\arraystretch}{1.1}
  \centering
  \caption{Assortment and revenue of heuristics and \eqref{eq:omni-mnl}}
  \resizebox{\textwidth}{!}{
    \begin{tabular}{cccccccccccc}
    \toprule
    \multirow{2}[4]{*}{\textbf{Channel}} & 
    \multirow{2}[4]{*}{\textbf{Type}} & 
    \multirow{2}[4]{*}{\textbf{Rate}} & 
    \multicolumn{3}{c}{\textbf{\ref{eq:2SNR}}} & 
    \multicolumn{3}{c}{\textbf{Algorithm \ref{alg:2SRO-heuristics}}} &
    \multicolumn{3}{c}{\textbf{\ref{eq:omni-mnl}}} \\
    \cmidrule{4-6}\cmidrule{7-9}\cmidrule{10-12} 
        &  &  & \textbf{Assort.} & \textbf{Exp.Rev} & \textbf{Total Rev.} & \textbf{Assort.} & \textbf{Exp.Rev} & \textbf{Total Rev.} & \textbf{Assort.} & \textbf{Exp.Rev} & \textbf{Total Rev.} \\
    \midrule
    Offline 
        & 0  & 0.7 & \{1\} & 19.82 & \multirow{3}{*}{15.53} & \{1,2,3\}& 15.94 &  \multirow{3}{*}{15.72} & \{1,3\}  & 19.75  & \multirow{3}{*}{18.39} \\
    \multirow{2}{*}{Online} 
        & 1  & 0.2 & \{1\} & 5.00 &      & \{3\}     & 17.82 &   &\{3\}   & 17.82   &    \\
        & 2  & 0.1 & \{1\} & 6.67 &      & \{1,2,3\} & 10.00 &   &\{1,3\} & 10.00  &    \\
    \bottomrule
    \end{tabular}}
  \label{table:toy_example_2_heuristics}
\end{table}
\end{example}
}

\rvtwo{Table~\ref{tab:heuristic_gap} reports the average (Ave) and maximum (Max) optimality gaps of Algorithm~\ref{alg:2SRO-heuristics} across 12 randomly generated instances. These instances are constructed using the same procedure described in Section~\ref{subsec:num-setting}. Specifically, we vary the problem sizes as \( (n, m) \in \{(20, 10), (50, 10), (50, 25), (100, 50)\} \), the offline arrival ratio as \( \alpha_0 \in \{0.1, 0.5\} \), and the online no-purchase utility weight as \( u_{i0} \in \{2, 5, 10\} \). We also introduce a cardinality constraint on the offline segment, denoted by \( K/n \), where the offline assortment must include no more than \( 0.1 \times n \) products when \( K/n = 0.1 \), and no constraint is applied when \( K/n = 1 \).}

\rvtwo{From the results in Table~\ref{tab:heuristic_gap}, we observe that the optimality gap tends to increase under the following conditions: when a cardinality constraint is applied, when 2SLM constraints are involved, when the offline segment represents a smaller proportion of total demand, and when both the product universe and the number of online consumer segments are reduced.}

\begin{algorithm}[hbpt]
\fontsize{9pt}{10pt}\selectfont 
\renewcommand{\baselinestretch}{1.2}\selectfont
\caption{2-Step Revenue-Order Heuristics}\label{alg:2SRO-heuristics}
\KwData{$u_{ij}$ for $i\in M^+, j\in N\cup\{0\}$; $r_{ij}$ for $i\in M^+, i\in N$; $G(N, E^i)$ for $i\in M$\;}
\KwResult{Offline offer set $S^{\text{off}}$; Online offer sets $S^{\text{on}} = \{S^{\text{on}}_1, \dots, S^{\text{on}}_M\}$\;}
\uIf{having cardinality constraint on offline}{
    $S^{\text{off}} \gets $ solve the MNL model for offline with cardinality constraint\;
}
\Else{
Sort $\{r_{0j}\}_{j\in N}$ as $r_{0\sigma_1}\ge r_{0\sigma_2} \ge \dots \ge r_{0\sigma_N}$\;
\For{$j$ in $N$}{
    $S \gets \{\sigma_1, \dots, \sigma_j\}$\;
    $z^{\text{off}}_j \gets \alpha_0 \times MNL(0, S)$; {\color{blue}\tcp{def: $MNL(i,S) \coloneqq \frac{\sum_{j\in S}r_{ij}u_{ij}}{u_{i0}+\sum_{j\in S}u_{ij}}$.}}
    $(S^{\text{on}}, z^{\text{on}}) \gets $ \texttt{EnumOnline}($S$)\; 
    $Z_j \gets z^{\text{off}}_j + z^{\text{on}}$\;
}
$j^\prime \gets \arg\max_{j\in N} Z_j$; $S^{\text{off}} \gets \{\sigma_1, \dots, \sigma_{j^\prime }\}$\; 
}
$(S^{\text{on}}, z^{\text{on}}) \gets $ \texttt{EnumOnline}($S^{\text{off}}$); {\color{blue}\tcp{finalize $S^{\text{on}}$.}} 
\SetKwFunction{onlineA}{EnumOnline}
\SetKwProg{Fn}{Function}{:}{}
\Fn{\onlineA{$S^{\textup{off}}$}}{
  \For{$i$ in $M$}{
      Sort  $\{r_{ij}\}_{j\in S^{\textup{off}}}$ as $r_{i\sigma_1}\ge r_{i\sigma_2 } \ge \dots \ge r_{i|S^{\textup{off}}|}$\;
      
    \eIf{no 2SLM}{
      $z_i \gets 0$\;
      \For{$j$ = $1:|S^{\textup{off}}|$}{
        $S \gets \{\sigma_1, \dots, \sigma_j\}$\;
            \eIf{$z_i< MNL(i, S)$}
            {$z_i \gets MNL(i, S)$; $S^{\text{on}}_i \gets S$\;
            }
            {break\;}
        }
      }
    {
        \For{$j$ = $1:|S^{\text{off}}|$}{
        $S \gets \{\sigma_1, \dots, \sigma_j\}$\;
        $\hat{S} \gets $ \texttt{applyLuce}($S$, $E^i$)\;
        $z_{ij} \gets MNL(i, \hat{S})$\;
      }
      $j^\prime \gets \arg\max_{1 \le j \le |S^{\text{off}}|} z_{ij}$\;
      $z_i \gets z_{ij^\prime}$; $S^{\text{on}}_i \gets \{\sigma_1, \dots, \sigma_{j^\prime}\}$
    }
  }
  \Return $S^{\text{on}}$, $\sum_{i\in M} \alpha_i\times z_i$\;
}
\SetKwFunction{applyLuce}{applyLuce}
\SetKwProg{Fn}{Function}{:}{}
\Fn{\applyLuce{$S$, $E$}}{
$D \gets \emptyset$\;
\For{$(j^\prime, j^{\prime\prime})$ in $E$ }{
    \If{$j^\prime \in S$ and $j^{\prime\prime} \in S$}{
    $D \gets D \cup \{j^{\prime\prime}\}$\;
    }
}
  \Return $S\backslash D$\;
}
\end{algorithm}

\begin{table}[H]
\footnotesize
\fontsize{9pt}{10pt}\selectfont 
\renewcommand\tabcolsep{3pt}
\renewcommand{\arraystretch}{1.1}
\centering
\caption{Optimality gap of the heuristics.}
\label{tab:heuristic_gap}
\begin{tabular}{cccccccc}
\toprule\multirow{3}{*}{$(n,m)$}  & \multirow{3}{*}{$\alpha_0$} & \multirow{3}{*}{$u_{i0}$} & \multirow{3}{*}{$K/n$} & \multicolumn{2}{c}{\textbf{without 2SLM}}   & \multicolumn{2}{c}{\textbf{with 2SLM}}   \\
\cmidrule{5-6} \cmidrule{7-8}
&  & & & Ave & Max & Ave & Max \\
\hline

(20, 10) & 0.1 & 2  & 0.1 & 7.59 & 14.12 & 7.97 & 16.04 \\
(20, 10) & 0.1 & 2  & 1   & 0.00 & 0.00  & 0.86 & 1.29  \\
(20, 10) & 0.1 & 5  & 0.1 & 6.78 & 13.73 & 7.03 & 15.31 \\
(20, 10) & 0.1 & 5  & 1   & 0.01 & 0.05  & 1.31 & 2.00  \\
(20, 10) & 0.1 & 10 & 0.1 & 4.66 & 11.09 & 5.01 & 12.39 \\
(20, 10) & 0.1 & 10 & 1   & 0.00 & 0.00  & 1.54 & 1.94  \\
(20, 10) & 0.5 & 2  & 0.1 & 1.20 & 5.10  & 1.44 & 6.03  \\
(20, 10) & 0.5 & 2  & 1   & 0.01 & 0.08  & 0.49 & 0.72  \\
(20, 10) & 0.5 & 5  & 0.1 & 0.57 & 3.44  & 0.70 & 4.00  \\
(20, 10) & 0.5 & 5  & 1   & 0.05 & 0.27  & 0.64 & 0.94  \\
(20, 10) & 0.5 & 10 & 0.1 & 0.21 & 1.89  & 0.31 & 2.23  \\
(20, 10) & 0.5 & 10 & 1   & 0.10 & 0.36  & 0.84 & 1.20  \\
\hline
(50, 10) & 0.1 & 2  & 0.1 & 5.26 & 10.66 & 6.10 & 12.39 \\
(50, 10) & 0.1 & 2  & 1   & 0.00 & 0.02  & 0.49 & 0.67  \\
(50, 10) & 0.1 & 5  & 0.1 & 6.83 & 14.71 & 7.97 & 17.41 \\
(50, 10) & 0.1 & 5  & 1   & 0.01 & 0.04  & 0.74 & 0.94  \\
(50, 10) & 0.1 & 10 & 0.1 & 6.70 & 15.32 & 7.86 & 17.82 \\
(50, 10) & 0.1 & 10 & 1   & 0.01 & 0.03  & 1.01 & 1.24  \\
(50, 10) & 0.5 & 2  & 0.1 & 0.79 & 3.10  & 1.26 & 4.05  \\
(50, 10) & 0.5 & 2  & 1   & 0.01 & 0.06  & 0.25 & 0.37  \\
(50, 10) & 0.5 & 5  & 0.1 & 0.79 & 3.25  & 1.16 & 3.83  \\
(50, 10) & 0.5 & 5  & 1   & 0.10 & 0.24  & 0.57 & 0.84  \\
(50, 10) & 0.5 & 10 & 0.1 & 0.54 & 2.18  & 0.81 & 2.65  \\
(50, 10) & 0.5 & 10 & 1   & 0.17 & 0.42  & 0.73 & 1.04  \\
\hline
(50, 25) & 0.1 & 2  & 0.1 & 2.65 & 9.30  & 3.28 & 9.27  \\
(50, 25) & 0.1 & 2  & 1   & 0.00 & 0.02  & 0.47 & 0.64  \\
(50, 25) & 0.1 & 5  & 0.1 & 2.99 & 10.05 & 3.50 & 10.82 \\
(50, 25) & 0.1 & 5  & 1   & 0.01 & 0.04  & 0.73 & 0.94  \\
(50, 25) & 0.1 & 10 & 0.1 & 2.75 & 8.91  & 3.36 & 9.10  \\
(50, 25) & 0.1 & 10 & 1   & 0.01 & 0.04  & 0.94 & 1.08  \\
(50, 25) & 0.5 & 2  & 0.1 & 0.47 & 2.28  & 0.74 & 2.55  \\
(50, 25) & 0.5 & 2  & 1   & 0.00 & 0.02  & 0.27 & 0.34  \\
(50, 25) & 0.5 & 5  & 0.1 & 0.35 & 1.83  & 0.61 & 2.06  \\
(50, 25) & 0.5 & 5  & 1   & 0.10 & 0.29  & 0.51 & 0.79  \\
(50, 25) & 0.5 & 10 & 0.1 & 0.19 & 1.08  & 0.38 & 1.26  \\
(50, 25) & 0.5 & 10 & 1   & 0.24 & 0.49  & 0.75 & 1.08  \\
\hline
(100, 50) & 0.1 & 2  & 0.1 & 1.80 & 3.41 & 2.49 & 4.05 \\
(100, 50) & 0.1 & 2  & 1   & 0.04 & 0.13 & 0.37 & 0.45 \\
(100, 50) & 0.1 & 5  & 0.1 & 2.08 & 4.05 & 3.26 & 5.21 \\
(100, 50) & 0.1 & 5  & 1   & 0.04 & 0.08 & 0.52 & 0.60 \\
(100, 50) & 0.1 & 10 & 0.1 & 2.21 & 4.30 & 3.70 & 5.63 \\
(100, 50) & 0.1 & 10 & 1   & 0.03 & 0.11 & 0.62 & 0.74 \\
(100, 50) & 0.5 & 2  & 0.1 & 0.39 & 0.79 & 0.76 & 1.18 \\
(100, 50) & 0.5 & 2  & 1   & 0.02 & 0.05 & 0.21 & 0.28 \\
(100, 50) & 0.5 & 5  & 0.1 & 0.36 & 0.69 & 0.78 & 1.12 \\
(100, 50) & 0.5 & 5  & 1   & 0.11 & 0.24 & 0.39 & 0.54 \\
(100, 50) & 0.5 & 10 & 0.1 & 0.27 & 0.55 & 0.71 & 1.03 \\
(100, 50) & 0.5 & 10 & 1   & 0.23 & 0.44 & 0.58 & 0.80 \\
\bottomrule
\end{tabular}
\end{table}

\end{document}